\documentclass[aop,MSNbibl,seceqn,dvips]{arximspdf}

\doi{10.1214/11-AOP736} 
\volume{41}
\issue{1}
\pubyear{2013}
\firstpage{1}
\lastpage{49}

\makeatletter

\newcommand{\rright}{\right}
\newcommand{\lleft}{\left}
\newcounter{Const} 
\def\Ct{\refstepcounter{Const}c_{\theConst}}

\newtheorem{thmm}{Theorem}[section]

\newtheorem{lem}[thmm]{Lemma}
\newtheorem{lemm}{Lemma}[section]

\newproclaim{dfn}{Definition}[section]
\newproclaim{ntt}[thmm]{Notation}
\newproclaim{cdn}[thmm]{Condition}
\newproclaim{ex}{Example}[section]
\newproclaim{prob}[thmm]{Problem}

\newproclaim{rem}{Remark}[section]

\newtheorem{Cora}{Corollary}[section] 

\newcommand\N{\mathbb{N}}
\newcommand\R{\mathbb{R}}
\newcommand\Z{\mathbb{Z}}
\newcommand\Rd{\R^d }
\newcommand{\map}[3]{#1\dvtx #2\to #3}
\newcommand\ts{ \times}
\newcommand{\limi}[1]{\lim_{#1\to\infty}}
\newcommand{\lsupi}[1]{\limsup_{#1\to\infty}}
\newcommand{\linfi}[1]{\liminf_{#1\to\infty}}
\newcommand{\PD}[2]{\frac{\partial #1}{\partial #2}}

\newcommand\SSSS{S}
\newcommand\SSS{\mathsf{S}}
\newcommand\nuN{\nu^{N}}
\newcommand\tCN{\tau_{r,rs }^N (\mathsf{x},\mathsf{s})}
\newcommand\tCNN{\tau_{r,rs }^N (\mathsf{x}',\mathsf{s})}
\newcommand\mA{\Lambda(d\mathsf{x})}
\newcommand\Aq{\mathsf{A}_q}
\newcommand\xxi{ x_i}
\newcommand\yyi{ y_i}

\newcommand\Sr{\SSSS_{r}}
\newcommand\Ss{\SSSS_{s}}
\newcommand\Srs{\SSSS_{rs}}
\newcommand\Sti{\tilde{\SSSS}}
\newcommand\Sri{\SSSS_{r\infty}}
\newcommand\Dr{\Sti_r}
\newcommand\Ds{\Sti_s}
\newcommand\Drs{\Sti_{rs}}
\newcommand\SZ{\SSS\bullet\SSSS}
\newcommand\Srhat{\tilde{\SSSS}_r}

\newcommand\wN{\varpi_{N}}
\newcommand\wM{\varpi_{M}}
\newcommand\wMxi{\wM(x_i)}
\newcommand\wNx{\wN(x)}
\newcommand\wNy{\wN(y)}
\newcommand\wNxi{\wN(x_i)}
\newcommand\wNyi{\wN(y_i)}
\newcommand\wNxl{\wN(x)^{\ell}}
\newcommand\wNwl{\wN(w)^{\ell}}

\newcommand\KN{K_{N}} 
\newcommand\KNstar{K_{N}^*}
\newcommand\KNover{\overline{K}_{N}}
\newcommand\KNstarover{\overline{K}_{N}^{*}}

\newcommand\NN{N}
\newcommand\ii{\mathrm{i}}

\newcommand\1{\Psi_r (\mathsf{x},\mathsf{s})}
\newcommand\2{\sup_{N\in\N} \sup_{r< s\in\N}\sup_{x\not=w\in\Sr}}
\newcommand\3{\frac{1}{({n_{N}})^d}\sum_{\xi\in\TN\cap({\Z^d}/{n_{N}})}}
\newcommand\4{\sum_{\yyi\in\Srs}}
\newcommand\6{\xi\in\TN\cap({\Z}/{n_{N}})}
\newcommand\7{x\not=w\in\Sr}
\newcommand\8{$ \murkm$-a.e. $ \mathsf{s} $}
\newcommand\9{\sum_{\yyi\in\Srs} \frac{1}{\bar{\varpi}_{N}(y_i)^{\ell} }}

\newcommand\V{\operatorname{Var}^{\mugN}}
\newcommand\Vz{\operatorname{Var}^{\mug}}
\newcommand\VV{\operatorname{Var}^{\mugNstar}}

\newcommand\xx{\mathsf{x}}

\newcommand\IN{\mathbb{I}_N}
\newcommand\INN{\mathbb{I}_{N+1}}
\newcommand\TN{\mathbb{T}_N}
\newcommand\T{\mathcal{T}_{N}}

\newcommand\HH{\mathsf{H}}
\newcommand\Hrk{\HH_{r,k}}
\newcommand\Hsl{\HH_{s,l}}

\newcommand\PhiN{\Phi^N }

\newcommand\PsiN{\Psi^{\NN}}
\newcommand\PsiNrs{\PsiN_{r,rs}} 
\newcommand\PsiNrstu{\PsiN_{rs,tu}} 
\newcommand\PsiNrst{\PsiN_{r,st}} 
\newcommand\PsiNrsi{\PsiN_{r,s\infty}} 
\newcommand\PsiNssi{\PsiN_{rs,s\infty}} 

\newcommand\ABN{_{r,k,\mathsf{s},rs}^{N,m}}
\newcommand\AB{_{r,k,\mathsf{s},rs}^{m}}
\newcommand\ABNs{_{r,k,\mathsf{s},rs}^{N,m}}
\newcommand\ABs{_{r,k,\mathsf{s},rs}^{m}}
\newcommand\ABss{_{r,k,\pi_{\Srs}(\mathsf{s}),rs}^{m}}

\newcommand\mN{\mathfrak{m}_{N}}
\newcommand\mukNm{\mu^{N,m}_{r,k}}
\newcommand\muNm{\mu^{N,m}_{r}}
\newcommand\murkm{\mu_{r,k}^{m}}

\newcommand\murky{\mu_{r,k,\mathsf{s}}^m} %
\newcommand\mury{\mu_{r,\mathsf{s}}^m} %
\newcommand\muN{\mu^{N}}
\newcommand\muA{\mu_{r,k,\mathsf{s}}^{ m }}
\newcommand\muAB{\mu\AB}
\newcommand\muABs{\mu\ABs}
\newcommand\muABN{\mu\ABN}

\newcommand\murmy{\mu_{r,k,\mathsf{s}}^m } %
\newcommand\murm{\mu_{r}^m }
\newcommand\murmch{\check{\mu}_{r,k,\mathsf{s}}^m}%

\newcommand\mug{\mu_{\mathrm{gin}}}
\newcommand\muNg{\mu_{\mathrm{gin}}^{N}}
\newcommand\muNgcheck{\check{\mu}_{\mathrm{gin}}^{N}}

\newcommand\KsinNb{\mathsf{K}^{N}_{\mathrm{dys},\beta}}
\newcommand\KsinNx{\mathsf{K}^{N}_{\mathrm{dys},1 }}
\newcommand\KsinNy{\mathsf{K}^{N}_{\mathrm{dys},2 }}
\newcommand\KsinNz{\mathsf{K}^{N}_{\mathrm{dys},4 }}
\newcommand\Ksinb{\mathsf{K}_{\mathrm{dys},\beta}}
\newcommand\Ksinx{\mathsf{K}_{\mathrm{dys},1 }}
\newcommand\Ksiny{\mathsf{K}_{\mathrm{dys},2 }}
\newcommand\Ksinz{\mathsf{K}_{\mathrm{dys},4 }}

\newcommand\mugN{\mu^{N}_{\mathrm{gin}}}
\newcommand\mugNstar{\mu^{N*}_{\mathrm{gin}}}
\newcommand\rg{\rho_{\mathrm{gin}}^n}
\newcommand\rgNx{\rho^{N,1}_{\mathrm{gin}}}
\newcommand\rgx{\rho_{\mathrm{gin}}^{1}}
\newcommand\kg{\mathsf{K}_{\mathrm{gin}}}

\newcommand\rN{\rho_N^{n}}
\newcommand\rNone{\rho_N^{1}}
\newcommand\rNtwo{\rho_N^{2}}
\newcommand\rbN{\rho_{N}^{n}}
\newcommand\rgN{\rho_{\mathrm{gin}}^{N,n}}
\newcommand\kgN{\mathsf{K}^{N}_{\mathrm{gin}}}

\newcommand\LAB{L^{\infty}(\SSS, \Lambda)}
\newcommand\LABone{L^{1}(\SSS, \Lambda)}

\newcommand\skAB{\sigma\AB}
\newcommand\skABs{\sigma\ABs}
\newcommand\skABNs{\sigma\ABNs}

\newcommand\Lm{L^2(\SSS, \mu)}
\newcommand\Ln{L^2(\SSS, \nu)}
\newcommand\LmNone{L^1(\SSS, \muN)}
\newcommand\Lmk{L^2(\SSS, \murkm)}

\newcommand\dom{\mathcal{D}^{\mu}}
\newcommand\domi{\mathcal{D}^{\mu}_{\infty}}
\newcommand\domai{\mathcal{D}^{ a ,\mu}_{\infty}}
\newcommand\domaz{\mathcal{D}^{ a_0I ,\mu}_{\infty}}
\newcommand\doma{\mathcal{D}^{ a ,\mu}}

\newcommand\Em{\mathcal{E}^{\mu}}
\newcommand\Ema{\mathcal{E}^{a,\mu}}
\newcommand\Erm{\mathcal{E}^{m, a,\mu}_r}
\newcommand\Ermk{\mathcal{E}^{m, a,\mu}_{r,k}}
\newcommand\Ermky{\mathcal{E}^{m, a,\mu}_{r,k,\mathsf{s}}}
\newcommand\Eaz{\mathcal{E}^{m, a_0I,\mu}_{r,k,\mathsf{s}}}

\newcommand\Dsp{(\Em,\dom,\Lm)}
\newcommand\Dspa{(\Ema,\doma,\Lm)}

\newcommand\ur{\mathsf{w}^{j}_{r }}
\newcommand\urr{\mathsf{w}^{j-1}_{r }}
\newcommand\up{\mathsf{w}^{j}_{q }}
\makeatother

\begin{document}
\begin{frontmatter}

\title{Interacting Brownian motions in infinite dimensions with
logarithmic interaction potentials\thanksref{T1}}
\thankstext{T1}{Supported by the Grant-in-Aid for Scientific Research
(B) 21340031 (Japan).}
\runtitle{Interacting Brownian motions with log potentials}

\begin{aug}
\author[A]{\fnms{Hirofumi} \snm{Osada}\corref{}\ead[label=e1]{osada@math.kyushu-u.ac.jp}}
\runauthor{H. Osada}
\affiliation{Kyushu University}
\address[A]{Faculty of Mathematics\\
Kyushu University\\
Fukuoka 819-0395\\
Japan \\
\printead{e1}}
\end{aug}

\received{\smonth{6} \syear{2009}}
\revised{\smonth{1} \syear{2011}}

%
\begin{abstract}
We investigate the construction of diffusions consisting of infinitely
numerous Brownian particles
moving in $ \mathbb{R}^d $ and interacting via logarithmic functions
(two-dimensional Coulomb potentials).
These potentials are very strong and act over a long range in nature.
The associated equilibrium states are no longer Gibbs measures.

We present general results for the construction of such diffusions and,
as applications thereof, construct two typical interacting Brownian
motions with logarithmic interaction potentials, namely the Dyson model
in infinite dimensions and Ginibre interacting Brownian motions.
The former is a particle system in $ \mathbb{R} $, while the latter
is in $ \mathbb{R}^2 $. Both models are translation and rotation
invariant in space,
and as such, are prototypes of dimensions $ d = 1,2 $, respectively.
The equilibrium states of the former diffusion model are
determinantal or Pfaffian random point fields with sine kernels.
They appear in the thermodynamical limits of the spectrum of
the ensembles of Gaussian random matrices such as GOE, GUE
and GSE. The equilibrium states of the~latter diffusion model are the
thermodynamical limits of the spectrum of the ensemble of complex
non-Hermitian Gaussian random matrices known as the Ginibre ensemble.
\end{abstract}

%


\end{frontmatter}

\section{Introduction} \label{s1}
Interacting Brownian motions (IBMs) in infinite dimensions
are diffusions $ \mathbf{X}_t = (X_t^i)_{i\in\Z} $ consisting of
infinitely many particles moving in $ \Rd$ with the effect of the
external force coming from a self-potential $ \map{\Phi}{\Rd}{\R
\cup\{ \infty\}}$ and that of the mutual interaction coming from an
interacting potential
$ \map{\Psi}{\Rd\ts\Rd}{\R\cup\{ \infty\}} $ such that $ \Psi
(x,y) = \Psi(y,x)$.

Intuitively, an IBM is described by the infinitely dimensional
stochastic differential equation (SDE) of the form
%
\begin{equation}
\label{11} \qquad dX_t^i = dB_t^i
- \frac{1}{2} \nabla\Phi\bigl(X_t^i\bigr) \,dt -
\frac{1}{2} \sum_{j \in\Z, j\not= i } \nabla\Psi
\bigl(X_t^i , X_t^j \bigr)
\,dt \qquad(i \in\Z) .
\end{equation}
The state space of the process
$\mathbf{X}_t = (X_t^i)_{i\in\Z}$ is $(\Rd)^{\Z}$ by construction.
Let $\mathsf{X}$ be the configuration-valued process given by
%
\begin{equation}
\label{12}  \mathsf{X}_t = \sum_{i\in\Z}
\delta_{X_t^i} .
\end{equation}
Here $ \delta_a $ denotes the delta measure at $ a $ and
a configuration is a Radon measure consisting of
a sum of delta measures.
We call $ \mathbf{X}$ the labeled dynamics and
$ \mathsf{X}$ the unlabeled dynamics.

The SDE (\ref{11}) was initiated by Lang \cite{La1,La2}.
He studied the case $ \Phi= 0 $, and
$ \Psi(x,y) = \Psi(x-y) $, where $ \Psi$
is of $ C^3_0(\Rd) $, superstable and regular
according to Ruelle \cite{ruelle2}. With the last two assumptions,
the corresponding unlabeled dynamics~$ \mathsf{X}$ has Gibbsian
equilibrium states.
See \cite{shiga,Fr} and \cite{T2} for other works concerning
the SDE (\ref{11}).

In \cite{odfa} the unlabeled diffusion was constructed using the
Dirichlet form.
The advantage of this method is that it gives a general and simple
proof of construction,
and more significantly, it allows us to apply singular interaction
potentials, which are particularly of interest, such as the
Lennard--Jones 6--12 potential and hard core potential.
We note that all these potentials were excluded in the SDE approach.
See \cite{yoshida,ark,taneudf} and \cite{yuu05}
for other works on applying the Dirichlet form approach to IBMs.

We remark that in all these works, except some parts of \cite{odfa},
the equilibrium states are supposed to be Gibbs measures
with Ruelle's class interaction potentials~$ \Psi$.
Thus, the equilibrium states are described
by the Dobrushin--Lanford--Ruelle (DLR) equations [see (\ref{qg5})],
the usage of which plays a pivotal role in the previous works.

The purpose of this paper is to construct
unlabeled IBMs in infinite dimensions
with the logarithmic interaction potentials
%
%
\begin{equation}
\label{13}  \Psi(x,y) = -\beta\log|x-y| .
\end{equation}
We present a sequence of general theorems to construct IBMs
and apply these to logarithmic potentials.
We remark that the equilibrium states
are not Gibbs measures because
the logarithmic interaction potentials are unbounded at infinity.

The above potential $ \Psi$ in (\ref{13}) is known
to be the two-dimensional Coulomb potential.
In practice, such systems are regarded as one-component plasma
consisting of equally charged particles.
To prevent the particles from all repelling to explode,
a neutralizing background charge is imposed.
The self-potential $ \Phi$ denotes this particle--background
interaction; see \cite{for}.

We study two typical examples, namely
Dyson's model (Section~\ref{s21}) and Ginibre IBMs (Section~\ref{s22}).
In the first example, we take $ d = 1 $, $ \Phi= 0 $ and
$ \Psi(x,y) = -\beta\log|x-y| $ ($ \beta= 1,2,4 $),
while in the second $ d = 2 $, $ \Phi(z) = |z|^2 $, and
$ \Psi(x,y) = -2 \log|x-y| $.

For the special values $ \beta= 1,2,4 $ and particular self-potentials
$ \Phi$,
the associated equilibrium states are limits of the spectrum of random matrices.
Recently, much intensive research has been carried out on random point fields
related to random matrices. Our purpose in this paper is
a rather more dynamical one; that is,
we construct diffusions, the equilibrium states
of which are these random point fields related to random matrices.

The labeled dynamics of the Dyson model in infinite dimensions
is represented by the following SDE:
%
%
\begin{equation}
\label{14}  dX_t^i = dB_t^i
+\frac{\beta}{2} \limi{R} \sum_{|X_t^j|\le R , j \in\Z, j\not= i }
\frac{1}{X_t^i - X_t^j } \,dt\qquad (i \in\Z) .
\end{equation}
Here $ \beta= 1 , 2, 4 $,
corresponding to the Gaussian orthogonal ensemble (GOE), the Gaussian
unitary ensemble (GUE) and
the Gaussian symplectic ensemble (GSE), respectively.
The invariant probability measures $ \mu_{\mathrm{dys}, \beta}$ of
the (unlabeled) Dyson models are translation invariant.
Hence, if the distribution of $ \mathsf{X}_0 $ equals $ \mu_{\mathrm
{dys}, \beta}$, then
for all $ t $,
%
\begin{equation}
\label{15} \sum_{ j \in\Z, j\not= i } \frac{1}{|X_t^i - X_t^j| } \,dt =
\infty\qquad\mbox{a.s} .
\end{equation}
This means that only conditional convergence is possible
in the summation of the drift term in (\ref{14}), which
is the cause of the difficulty in dealing with the Dyson model.
It is well known that the equilibrium states are
the thermodynamic limits of the distribution of
the spectrum of Gaussian random matrices at the bulk
\cite{so-,for,mehta}.

The labeled dynamics of Ginibre IBMs is represented by the following SDE.
For convenience, we regard $ \SSSS$ as $ \mathbb{C} $ rather than $
\R^2 $.
%
\begin{equation}
\label{17} \qquad dZ_t^i = dB_t^i
- Z_t^i \,dt + \limi{R } \sum
_{|Z_t^j | \le R , j \in\Z, j\not= i } \frac{ Z_t^i - Z_t^j }{| Z_t^i - Z_t^j |^2 } \,dt\qquad  (i \in\Z) .
\end{equation}
Here $ Z_t^i = X_t^i + \ii Y_t^i \in\mathbb{C} $, where
$ \ii= \sqrt{-1}$, and $ \{ B_t^i \}_{i \in\Z}$
are independent complex Brownian motions. That is,
$ B_t^i = B_t^{i, \operatorname{Re}} + \ii B_t^{i, \operatorname{Im}}$,
where
$ \{ B_t^{i, \operatorname{Re}}, B_t^{i, \operatorname{Im}} \}_{i \in
\Z}$
is a system of independent one-dimensional Brownian motions.
The stationary measure~$ \mug$ of the unlabeled dynamics is the
thermodynamic limit of the distribution of the spectrum of random
Gaussian matrices called the Ginibre ensemble; cf. \cite{so-}.
$ \mug$ is a random point field with logarithmic interaction potential
and is known to be translation invariant.
If Ginibre IBMs
$ \mathsf{Z} = \{ \mathsf{Z}_t \} = \{ \sum_i\delta_{Z_t^i} \}$
start from the stationary measure $ \mug$,
then $ \mathsf{Z}$ is also
translation invariant in space. Moreover, Ginibre IBMs
$ \mathsf{Z}$ satisfy the SDE of
the translation invariant form
%
\begin{equation}
\label{214b0}  dZ_t^i = dB_t^i
+ \limi{R } \sum_{|Z_t^i -Z_t^j | \le R , j \in\Z, j\not= i } \frac{ Z_t^i - Z_t^j }{| Z_t^i - Z_t^j |^2 } \,dt\qquad  (i \in
\Z) .
\end{equation}
This variety of SDE representations of
Ginibre IBMs is a result
of the strength of the interaction potential.

A diffusion
$ (\mathsf{X},\mathsf{P})$ is a family of probability measures
$ \mathsf{P} = \{ \mathsf{P}_{\mathsf{x}} \} $
with continuous sample path $ \mathsf{X} = \{ \mathsf{X}_t\}$
starting at each point $ \mathsf{x} $ of the state space
with a strong Markov property; see \cite{fot}.
We emphasize that we construct not only a Markov semi-group or a
stationary Markov process, but also a diffusion in the above sense, and
also that, to apply stochastic analysis effectively, we require the
construction of diffusions.

In \cite{oisde}, we give another general result for the SDE
representation of unlabeled diffusions constructed in this paper.
The SDEs (\ref{14}), (\ref{17}) and (\ref{214b0}) of the labeled
dynamics are solved there using the main results Theorems~\ref{l22}
and~\ref{l23} in the present paper.
These SDEs provide a clear trajectory level description of the
diffusions obtained in the present paper. We also note that in \cite
{oisde} the fully labeled dynamics $ \mathbf{X}_t$ is a diffusion on
$ \R^{\Z}$ (Dyson's model) and $ (\R^2)^{\Z}$ (Ginibre IBMs).

Because of the long range nature of the logarithmic interaction,
the diffusion has not yet been constructed.
The only exception is the Dyson model with $ \beta= 2 $.
In \cite{sp2} Spohn proved the closability of
the Dirichlet form associated with (\ref{11}) for this model.
This implies the construction of the unlabeled dynamics~(\ref{12}) in the sense of an $ L^2 $-Markovian semigroup.
An associated diffusion was constructed in~\cite{odfa}
by combining Spohn's result with the result from \cite{odfa}, Theorem~0.1,
for the quasi-regularity of Dirichlet forms.

In one space dimension, some explicit computations of
space--time correlation functions of infinite particle systems
related to random matrices have been obtained.
Indeed, Katori and Tanemura \cite{kt07} recently studied the
thermodynamic limit of the space--time correlation functions related to
the Dyson model and Airy process.
Their limit space--time correlation functions define a stochastic process
starting from a limited set of initial distributions.
However, the Markov (semi-group) property of the process
has not yet been proved.
They also proved that, if their process is Markovian, the associated
Dirichlet form is the same as
the one obtained in this paper and their processes coincide with the
processes constructed here. It is an interesting open problem to prove
the Markov property of their processes
and identify these two processes.\setcounter{footnote}{1}\footnote{Recently, the Markov
property of the processes in \cite{kt07} has been proved in \cite
{kt11}. The domains of Dirichlet forms in \cite{kt11} include ones
in the present paper, but the identification of these two kinds of
Markov processes is still open.}

We also refer to
\cite{johansson03,knt04,kt07-ptrf}
and \cite{p-spohn02} 
for stochastic processes of one-dimensional infinite particle systems
related to random matrices.

As for two-dimensional infinite systems with logarithmic interactions,
the construction of stochastic processes based on the explicit
computation of space--time correlation functions
has not been done. Techniques useful in one-dimension,
such as applying the Karlin--McGregor formula, are no longer valid in
two dimensions.

Let us briefly explain the main idea.
We introduce the notion of quasi-Gibbs measures as a substitution for
Gibbs measures. These measures satisfy
inequality~(\ref{qg2})
involving a~(finite volume) Hamiltonian.
Inequality (\ref{qg2}) is sufficient for the closability of the
Dirichlet forms and the construction of the diffusions.

To obtain the above-mentioned inequality
we control the difference of the infinite volume Hamiltonians
instead of the Hamiltonian, itself.
The key point of the control is the usage of the geometric property of
the random point fields behind the dynamics.
Indeed, although the difference still diverges for Poisson random
fields and Gibbs measures
with translation invariance, it becomes finite for random point fields
such as
Dyson random point fields and Ginibre random point fields.
For these random point fields the fluctuations of particles are
extremely suppressed because the logarithmic potentials are quite strong.
This cancels the sum of the difference of the infinite-volume Hamiltonians.

The organization of the paper is as follows.
In Section~\ref{s2}, we describe the set-up and state the main results
(Theorems~\ref{l22} and~\ref{l23}).
We first introduce the notion of quasi-Gibbs measures and
give a general result (Lemma~\ref{l21}) concerning the closability of
bilinear forms.
As applications, we then construct the diffusions of the Dyson model
and the Ginibre IBMs cited above in Theorems~\ref{l22} and~\ref{l23}, respectively.
Section~\ref{s3} is devoted to preparation from the Dirichlet form theory
and the proof of Lemma~\ref{l21}.
The most crucial assumption of Lemma~\ref{l21} is the quasi-Gibbs property.
In Section~\ref{s4}, we introduce Theorem~\ref{l41},
which gives a pair of sufficient conditions (\hyperlink{A4}{A.4}) and (\hyperlink{A5}{A.5})
for the quasi-Gibbs property. We also explain the strategy of the
proof of Theorem~\ref{l41}.
In Section~\ref{s5}, we prove Theorem~\ref{l41}.
In Section~\ref{s6}, we prove Theorem~\ref{l63},
which allows us to deduce (\hyperlink{A5}{A.5}) from the new condition (\hyperlink{A6}{A.6}).
In Section~\ref{s7}, we give a sufficient condition of (\hyperlink{A6}{A.6}),
directly used in the proof of Theorems~\ref{l22} and~\ref{l23}.
In Section~\ref{s8}, we give a representation of the $ L^2$-norm of linear
statistics in terms of Fourier series when random fields are periodic,
which is a preparation of the proof of Theorem~\ref{l22}.
In Section~\ref{s9}, we prove Theorem~\ref{l22}.
In Section~\ref{s10}, we prove Theorem~\ref{l23}.
In Appendix~\ref{sA1} we prove Lemmas~\ref{l34} and~\ref{l35},
in Appendix~\ref{sA4} we prove Lemma~\ref{l111}.

\section{Set up and main results}\label{s2}
Let $\SSSS$ be a closed set in $ \Rd$ such that $ 0 \in\SSSS$ and $
\overline{\SSSS}^{\mathrm{int}}  = \SSSS$,
where $ \SSSS^{\mathrm{int}} $ means the interior of $ \SSSS$.
Let
$ \SSS= \{ \mathsf{s} = \sum_i \delta_{s_i} ; \mathsf{s}(K )<
\infty
\mbox{ for any compact set } K \} $, where $ \{ s_i \} $ is a sequence
in $ \SSSS$.
Then $ \SSS$ is the set of configurations on $ \SSSS$ by definition.
We endow $ \SSS$ with vague topology,
under which $ \SSS$ is a Polish space.

Let $ \mu$ be a probability measure on
$ (\SSS, \mathcal{B}(\SSS) ) $.
We construct $ \mu$-reversible diffusions
$(\mathsf{X},\mathsf{P})$
with state space $ \SSS$ using the Dirichlet form theory.
Hence, we begin by introducing Dirichlet forms in the following.

For a subset $ A \subset\SSSS$, we define the map
$ \map{\pi_{A }}{\SSS}{\SSS} $ by
$ \pi_{A } (\mathsf{s}) = \mathsf{s}( A \cap\cdot) $.
We say a function $ \map{f}{\SSS}{\R} $
is local if $ f $ is $ \sigma[\pi_{ A }]$-measurable
for some bounded Borel set~$ A $.
We say $ f $ is smooth if $ \tilde{f} $ is smooth,
where $ \tilde{f}((s_i)) $ is the permutation invariant function in $
(s_i) $ such that
$ f (\mathsf{s}) = \tilde{f} ((s_i)) $ for
$ \mathsf{s} = \sum_i \delta_{s_i} $.

Let
$ \SZ= \{ (\mathsf{s}, s )\in\SSS\ts\SSSS ;
\mathsf{s}(\{ s \})\ge1 \}
$.
Let $ \map{a = (a_{kl})}{\SZ}{\R^{d^2}} $ be such that
$ a_{kl} = a_{lk} $ and $ (a_{kl}(\mathsf{s}, s )) $
is nonnegative definite. Set
%
%
\begin{equation}
\label{A31}  \mathbb{D}^{a} [f,g] (\mathsf{s} ) = \frac{1}{2}
\sum_i \sum_{k,l = 1}^{ d }
a_{kl}(\mathsf{s} , s_i ) \PD{\widetilde{f}}
{s_{ik}}\cdot\PD{\widetilde{g}} {s_{il}} .
\end{equation}
Here $ s_i = (s_{i1},\ldots,s_{id}) \in\SSSS$ and
$ \mathsf{s} = \sum_i \delta_{s_i} $.
For given $ f $ and $ g $, it is easy to see that the right-hand side
depends only on $ \mathsf{s}$.
Therefore, the square field $ \mathbb{D}^{a}[f,g] $ is well defined.
We assume $ \map{\mathbb{D}^{a} [f,g] }{\SSS}{\R} $ is
$ \mathcal{B}(\SSS) $-measurable for each of the local, smooth
functions $ f $ and $ g $.

For $ a $ and $ \mu$, we consider the bilinear form $ (\Ema, \domai)
$ defined by
%
\begin{eqnarray}
\label{A30} \qquad\Ema(f,g) &=& \int_{\SSS} \mathbb{D}^{a}[f,g]
\,d\mu,
\nonumber
\\[-8pt]
\\[-8pt]
\nonumber
\domai&=& \bigl\{ f \in\Lm; f \mbox{ is local and smooth, } \Ema(f,f)<
\infty\bigr\} .
\end{eqnarray}
When $ a_{kl} = \delta_{kl} $ ($ \delta_{kl} $ is the Kronecker
delta), we write
$ \mathbb{D}^{a} = \mathbb{D} $,
$ \Ema= \mathcal{E}^{\mu} $, and $ \domai= \domi$.

All examples in this paper satisfy
$ a_{kl} = \delta_{kl} $. We, however, state the assumption in a
general framework.
We assume the coefficients $ \{ a_{kl} \}$
satisfy the following:

(A.0)\hypertarget{A0}\ There exists a nonnegative,
bounded, lower semicontinuous function
$ \map{a_0}{\SZ}{[0,\infty)} $ and
a constant
$ \Ct\label{;a33} $ $ \ge1 $ such that
%
%
\begin{equation}
\label{A33}  c_{\scriptsize\ref{;a33}}^{-1} a_0 (\mathsf{s},s )
|x|^2 \le \sum_{k,l = 1}^d
a_{kl} (\mathsf{s},s ) x_{k} x_{l} \le
c_{\scriptsize\ref{;a33}} a_0 (\mathsf{s},s ) |x|^2
\end{equation}
for all $ x = (x_1,\ldots,x_d)\in\Rd$, $(\mathsf{s}, s)
\in\SZ$.\vadjust{\goodbreak}

We call a function $ \rho^n $ the $ n $-correlation function of $ \mu$
with respect to (w.r.t.) the Lebesgue measure
if $ \map{\rho^n }{\SSSS^n}{\R} $ is a permutation invariant
function such that
%
\begin{equation}
\label{28} \int_{A_1^{k_1}\ts\cdots\ts A_m^{k_m}} \rho^n (x_1,
\ldots,x_n) \,dx_1\cdots dx_n = \int
_{\SSS} \prod_{i = 1}^{m}
\frac{\mathsf{s} (A_i) ! } {
(\mathsf{s} (A_i) - k_i )!} \,d\mu
\end{equation}
for any sequence of disjoint bounded measurable subsets
$ A_1,\ldots,A_m \subset\SSSS$ and a sequence of natural numbers
$ k_1,\ldots,k_m $ satisfying $ k_1+\cdots+ k_m = n $.
It is well known \cite{so-} that under a mild condition, the
correlation functions
$ \{ \rho^n \}_{n \in\N}$ determine the measure $ \mu$.

We assume $ \mu$ satisfies the following.

(A.1)\hypertarget{A1}\ The measure $ \mu$ has
a locally bounded, $ n $-correlation function
$ \rho^n $ for each $ n \in\N$.

We introduce a Hamiltonian on a bounded Borel set $ A $ as follows.
For Borel measurable functions
$ \map{\Phi}{\SSSS}{\R\cup\{\infty\}} $
and
$ \map{\Psi}{\SSSS\ts\SSSS}{\R\cup\{\infty\}} $
with $ \Psi(x,y) = \Psi(y,x) $, let
%
%
\begin{equation}
\label{2y}  \quad\mathcal{H}_{A }^{\Phi, \Psi} (\mathsf{x}) = \sum
_{x_i\in A } \Phi( x_i ) + \sum
_{x_i, x_j\in A , i < j } \Psi( x_i, x_j) \qquad\mbox{where }
\mathsf{x} = \sum_i \delta_{x_i} .
\end{equation}
We assume $ \Phi< \infty$ a.e. to avoid triviality.

For two measures $ \nu_1,\nu_2 $ on a measurable space
$ (\Omega, \mathcal{B})$ we write $ \nu_1 \le\nu_2 $
if $ \nu_1(A)\le\nu_2(A)$ for all $ A\in\mathcal{B}$.
We say a sequence of finite Radon measures $ \{ \nuN\} $
on a Polish space $ \Omega$ converge weakly to
a finite Radon measure $ \nu$
if $ \lim_{N\to\infty} \int fd\nuN= \int f \,d\nu$
for all $ f \in C_b(\Omega) $.

Throughout this paper, $ \{ b_r \}$ denotes
an increasing sequence of natural numbers. We set
%
\begin{equation}
\label{qg0}\Sr= \bigl\{ s \in\SSSS; |s| < b_r \bigr\},\qquad \SSS_r^{m}=
\bigl\{\mathsf{s}\in\SSS;\mathsf{s}(\Sr) = m \bigr\} .
\end{equation}

%
\begin{dfn}\label{dfn1}
A probability measure $ \mu$ is said to be
a $ (\Phi, \Psi)$-quasi Gibbs measure if
there exists an increasing sequence $ \{ b_r \} $
of natural numbers and measures $ \{ \murkm\} $
such that, for each $ r ,m \in\N$,
$ \murkm$ and $ \murm: = \mu(\cdot\cap\SSS_r^{m}) $
satisfy
%
%
\begin{equation}
\label{qg1} \murkm\le\mu_{r,k+1}^{m}\qquad \mbox{for all } k ,\qquad
\limi{k} \murkm= \murm\qquad \mbox{weakly,}
\end{equation}
and that, for all $ r,m,k \in\N$ and for \8 $\in\SSS$,
%
%
\begin{equation}
\label{qg2} c_{\scriptsize\ref{;2y}}^{-1} e^{-\mathcal{H}_{r}(\mathsf{x})} 1_{\SSS_r^{m}}(
\mathsf{x}) \Lambda(d\mathsf{x})\le \murky(d\mathsf{x}) \le c_{\scriptsize\ref{;2y}}
e^{-\mathcal{H}_{r}(\mathsf{x})} 1_{\SSS_r^{m}}(\mathsf{x}) \Lambda(d\mathsf{x}) .
\end{equation}
Here $ \mathcal{H}_{r}(\mathsf{x}) = \mathcal{H}_{\Sr}^{\Phi,\Psi
}(\mathsf{x})$,
$\Ct\label{;2y} = c_{\scriptsize\ref{;2y}}(r,m,k,\pi_{\Sr^c}(\mathsf{s}))$
is a positive constant, $ \Lambda$ is the Poisson random point field
whose intensity is the Lebesgue measure on $ \SSSS$ and
$ \murky$ is the conditional probability measure of
$ \murkm$ defined by
%
\begin{equation}
\label{qg4}  \murky(d\mathsf{x}) = \murkm\bigl(\pi_{ \Sr} \in d
\mathsf{x} | \pi_{ \Sr^c }(\mathsf{s})\bigr) .
\end{equation}
\end{dfn}

We call $ \Phi$ (resp., $ \Psi$) a free (interaction) potential.
When $ \Psi$ is an interaction potential, we implicitly assume that
$ \Psi(x,y)= \Psi(y,x)$.
Our second assumption is as follows.

(A.2)\hypertarget{A2}\ $ \mu$ is a $ (\Phi, \Psi) $-quasi Gibbs measure.

%
\begin{rem}\label{r1}
(1) By definition, $ \murkm((\SSS_r^{m})^c) = 0 $.
Since $ \murky$ is
$ \sigma[\pi_{ \Sr^c }]$-measurable in
$ \mathsf{s}$,
we have the disintegration of the measure $\murkm$
%
\begin{equation}
\label{qg6} \murkm\circ\pi_{ \Sr}^{-1}(d\mathsf{x}) = \int
_{\SSS} \murky(d\mathsf{x})\murkm(d \mathsf{s}) .\vspace*{-6pt}
\end{equation}
\begin{longlist}[(2)]
\item[(2)]
Let
$ \mury(d\mathsf{x}) = \murm(\pi_{ \Sr} (\mathsf{s})
\in d\mathsf{x} | \pi_{ \Sr^c }(\mathsf{s})) $.
Recall that a probability measure~$\mu$ is said to be
a $ (\Phi, \Psi) $-canonical Gibbs measure if
$ \mu$ satisfies the DLR equation~(\ref{qg5}), that is,
for each $ r , m \in\N$, the conditional probability $ \mury$ satisfies
%
%
\begin{equation}
\label{qg5} \mury(d\mathsf{x}) = \frac{1}{c_{\scriptsize\ref{;2yy}}} e^{-\mathcal{H}_{r}(\mathsf{x}) -
\1
}
1_{\SSS_r^{m}}(\mathsf{x}) \Lambda(d\mathsf{x}) \qquad\mbox{for $ \murm$-a.e. $
\mathsf{s} $.}
\end{equation}
Here $ 0 < \Ct\label{;2yy}< \infty$ is the normalization and,
for $ \mathsf{x} = \sum_i \delta_{x_i}$ and
$ \mathsf{s} = \sum_j \delta_{s_j}$, we set
%
%
\begin{equation}
\label{qg8} \1 = \sum_{x_i\in\Sr, s_j \in\Sr^c } \Psi( x_i,
s_j) .
\end{equation}
\end{longlist}

We remark that $ (\Phi, \Psi) $-canonical Gibbs
measures are $ (\Phi, \Psi) $-quasi Gibbs measures.
The converse is, however, not true.
When $ \Psi(x,y) = - \beta\log|x-y|$ and $ \mu$ are translation
invariant, $ \mu$ are not $ (\Phi, \Psi) $-canonical Gibbs measures.
This is because the DLR equation does not make sense.
Indeed, $ |
\1
| = \infty$ for $ \mu$-a.s. $ \mathsf{s}$.
The point is that one can expect a cancellation between
$ c_{\scriptsize\ref{;2yy}} $ and
$ e^{-
\1 
} $,
even if
$ |
\1
| = \infty$.
\end{rem}

(A.3)\hypertarget{A3}\ There exist
upper semicontinuous functions
$ \map{\Phi_0, \Psi_0 }{\SSSS}{\R\cup\{ \infty\}}$ and positive constants
$ \Ct\label{;a35} $ and $ \Ct\label{;a3x} $
such that
%
\begin{eqnarray}
\label{A34}  c_{\scriptsize\ref{;a35}}^{-1} \Phi_0 (s) &\le&
\Phi(s) \le c_{\scriptsize\ref{;a35}} \Phi_0 (s)
\\
\label{Ax} \qquad\quad c_{\scriptsize\ref{;a3x}}^{-1} \Psi_0 (s-t) &\le&
\Psi(s,t) \le c_{\scriptsize\ref{;a3x}} \Psi_0 (s-t),\qquad \Psi_0
(s) = \Psi_0 (-s) \qquad (\forall s) .
\end{eqnarray}
Moreover, $ \Phi_0 $ and $ \Psi_0 $ are locally bounded from below, and
$ \mathsf{\Gamma}: = \{ s ; \Psi_0 (s ) = \infty\}$ is a
compact set.

We use the following result obtained in \cite{odfa} and \cite{om}.

%
\begin{lem}[(\cite{odfa,om})] \label{l21}
Assume \textup{(\hyperlink{A0}{A.0})--(\hyperlink{A3}{A.3})}.
Then $ (\Ema, \domai,\Lm) $ is closable,
and its closure $ \Dspa$ is a local, quasi-regular Dirichlet space.
\end{lem}

See Section~\ref{s3} for the definition of
``a local, quasi-regular Dirichlet space'' and
necessary notions of the Dirichlet form theory.\vadjust{\goodbreak}
Combining Lemma~\ref{l21} with the Dirichlet form theory developed
in \cite{fot} and \cite{mr}, we obtain the following.

%
\begin{Cora}
Assume \textup{(\hyperlink{A0}{A.0})--(\hyperlink{A3}{A.3})}.
Then there exists a diffusion $(\mathsf{X},\mathsf{P})$ associated
with $\Dspa$.
Moreover, the diffusion $(\mathsf{X},\mathsf{P})$ is $ \mu$-reversible.
\end{Cora}
We say a diffusion $(\mathsf{X},\mathsf{P})$
is associated with the Dirichlet form
$ (\Ema,\doma) $ on $ \Lm$ 
if $ E_{\mathsf{x}}[f(\mathsf{X}_t )] = T_tf(\mathsf{x})$
$ \mu$-a.e. $ \mathsf{x}$ for all $ f\in\Lm$.
Here $ T_t $ is the $ L^2 $-semi group associated with
the Dirichlet space $ \Dspa$. Moreover,
$(\mathsf{X},\mathsf{P})$ is called $ \mu$-reversible if
$(\mathsf{X},\mathsf{P})$
is $ \mu$-symmetric, and $ \mu$ is an invariant probability measure
of $(\mathsf{X},\mathsf{P})$.
\subsection{The Dyson model in infinite dimensions (Dyson IBMs)}
\label{s21}
Let $ S = \R$. Let $ \mu_{\mathrm{dys}, \beta}$ $ (\beta= 1,2,4)$
be the probability measure on $ \SSS$
whose $ n $-correlation function $ \rho^n$ is given by
%
\begin{equation}
\label{224}  \rho^n(x_1,\ldots,x_n) = \det
\bigl[\Ksinb(x_i,x_j)\bigr]_{1 \le i,j \le n} .
\end{equation}
Here for $ \beta= 2 $, we take $ \Ksiny(x,y) = {\sin( \pi(x-y)
)}/{\pi(x-y)} $.
$ \Ksiny$ is called the sine kernel. We remark that $ \Ksiny(x,y) =
\frac{1}{2\pi}
\int_{|k| \le\pi} e^{\ii k (x-y)} \,dk $
and $ 0 \le\Ksiny\le\operatorname{Id}$ as an operator on $ L^2 (\R
) $.
It is known that $ \Ksiny$ generates a determinantal random point
field \cite{so-}.
The definition of $ \Ksinb$ for $ \beta= 1,4 $ is given by~(\ref{"91s}) and (\ref{"91t}).
We use quaternions to denote the kernel $ \Ksinb$ for $ \beta= 1,4 $.
The precise meaning of the determinant of (\ref{224}) for $ \beta=
1,4 $
is given by (\ref{"91p}).
%
%
\begin{thmm}\label{l22}
Let $ \Phi(x) = 0 $ and $ \Psi(x,y) = -\beta\log|x-y| $.
Then $ \mu_{\mathrm{dys}, \beta}$ is a~quasi-Gibbs measure with
potentials $ (\Phi,\Psi)$.
\end{thmm}

From Corollary~\ref{l21} and Theorem~\ref{l22}, we obtain:
\begin{Cora}
Let $ \Dsp$ be as in Lemma~\ref{l21}
with $ a = (\delta_{kl}) $ and $ \mu= \mu_{\mathrm{dys}, \beta}$.
Then there exists a $ \mu$-reversible diffusion
$(\mathsf{X},\mathsf{P})$ associated with $ \Dsp$.
\end{Cora}

%
\begin{rem}\label{r25}
(1) We write\vspace*{-1pt}
$ \mathsf{X}_t = \sum_{i\in\Z} \delta_{X_t^i}$.
Here $ \mathbf{X}_t = (X_t^i)_{i\in\Z}$
is the associated labeled dynamics. It is known \cite{ocol} that
particles $ X_t^i $ never collide with each other.
Moreover, in \cite{oisde}, we prove that
the associated labeled dynamics
$ (X_t^i)_{i\in\Z} $ is a solution of the SDE
%
\begin{equation}
\label{214}  dX_t^i = dB_t^i
+\frac{\beta}{2} \limi{R} \sum_{|X_t^j|\le R , j \in\Z, j\not= i }
\frac{1}{X_t^i - X_t^j } \,dt\qquad (i \in\Z)
\end{equation}
with $ (X_0^i) = (x_i) $ for $ \mu_{\mathrm{dys}, \beta}$-a.s. $
\mathsf{x} = \sum_i\delta_{x_i}$.\vspace*{-6pt}
\begin{longlist}[(2)]
\item[(2)]
We remark that $ \mu_{\mathrm{dys}, \beta}$ is translation invariant.
The dynamics $ \mathsf{X}_t $ inherits the translation invariance from
the equilibrium state $ \mu_{\mathrm{dys}, \beta}$.
Indeed, if $ \mathsf{X}_t $ starts from the distribution $ \mu_{\mathrm{dys}, \beta}$, then the distribution of $ \mathsf{X}_t $
becomes translation invariant in time and space.

\item[(3)] One can easily see that
$ \rho^1(x) = 1 $.
By scaling in space, we can treat $ \mu_{\mathrm{dys}, \beta}$ with intensity
$ \rho^1(x) = \bar{\rho}$ for any
$ 0 < \bar{\rho}< \infty$.
\end{longlist}
\end{rem}

\subsection{Ginibre interacting Brownian motions}\label{s22}
Next we proceed with the Ginibre IBMs. For this purpose, we first introduce
a Ginibre random point field, which is a stationary probability measure
for a Ginibre IBM.

Let the state space $ \SSSS$ of particles be $ \mathbb{C}$. Let
%
\begin{equation}
 \label{27a} \kg(z_1,z_2) = \frac{1}{\pi} \exp
\biggl( - \frac{|z_1|^2}{2}-\frac{|z_2|^2}{2} + z_1 \cdot
\bar{z}_2 \biggr) .
\end{equation}
Here $ z_1, z_2 \in\mathbb{C} $ and $ \bar{z}$ denotes the complex
conjugate of $ z \in\mathbb{C}$.
Let $ \mug$ be the probability measure
whose $ n $-correlation $ \rg$ is given by
%
\begin{equation}
 \label{27b} \rg(z_1,\ldots,z_n) = \det\bigl[
\kg(z_i,z_j)\bigr]_{1 \le i,j \le n} .
\end{equation}
We call $ \mug$ the Ginibre random point field.
It is well known \cite{mehta} that
$ \mug$ is the thermodynamic limit of the distribution of
the spectrum of the random Gaussian matrix called
the Ginibre ensemble (cf. \cite{so-}),
which is the ensemble of complex non-Hermitian random
$ N \ts N $ matrices whose 2$ N^2 $ parameters
are independent Gaussian random variables
with mean zero and variance~$1/2$.

%
%
\begin{thmm} \label{l23}
Let $ \Phi(z) = |z|^2 $ and $ \Psi(z_1,z_2) = -2 \log|z_1-z_2| $.
Then $ \mug$ is a quasi-Gibbs measure with potential $ (\Phi,\Psi)$.
\end{thmm}

From Corollary~\ref{l21} and Theorem~\ref{l23}, we obtain:
\begin{Cora}
Let $ \Dsp$ be as in Lemma~\ref{l21} with $ a = (\delta_{kl}) $ and $
\mu= \mug$.
Then there exists a $ \mu$-reversible diffusion
$(\mathsf{Z},\mathsf{P})$ associated with $ \Dsp$.
\end{Cora}
We write $ \mathsf{Z}_t = \sum_{i\in\Z} \delta_{Z_t^i}$.
In \cite{oisde}, we prove that the associated labeled dynamics $
(Z_t^i)_{i\in\Z} $ is a solution of the SDE
%
%
\begin{equation}
\label{1g} \qquad dZ_t^i = dB_t^i
- Z_t^i \,dt + \limi{R } \sum
_{|Z_t^j | \le R , j \in\Z, j\not= i } \frac{ Z_t^i - Z_t^j }{| Z_t^i - Z_t^j |^2 } \,dt \qquad (i \in\Z) .
\end{equation}
Here $ Z_t^i \in\mathbb{C} $ and
$ \{ B_t^i \}_{i \in\Z}$ are independent complex Brownian motions.

We remark that the kernel $ \mathsf{K}_{\mathrm{gin}}$
is \textit{not} translation invariant.
The measure $ \mug$ is, however, rotation and
translation invariant. Such invariance is inherited by
the unlabeled diffusion
$ \mathsf{Z}_t = \sum_{i\in\Z} \delta_{Z_t^i}$.
This may be surprising because\vspace*{-1pt} SDE (\ref{1g}) is not
translation invariant at first glance.
In \cite{oisde}, we prove that
$ (Z_t^i)_{i\in\Z} $ satisfies the following SDE
%
\begin{equation}
\label{214b}  dZ_t^i = dB_t^i
+ \limi{R } \sum_{|Z_t^i -Z_t^j | \le R , j \in\Z, j\not= i } \frac{ Z_t^i - Z_t^j }{| Z_t^i - Z_t^j |^2 } \,dt \qquad(i \in
\Z)
\end{equation}
if $ \mathsf{Z}_t $ starts from the distribution $ \mug$.
The passage from (\ref{1g}) to (\ref{214b}) is a result of
the cancellation between the repulsion of the mutual interaction
of the particles and the neutralizing background charge.

\section{Preliminaries from the Dirichlet form theory} \label{s3}
In this section, we prepare some results from the Dirichlet form theory
and give a proof of Lemma~\ref{l21}.
The proof of Lemma~\ref{l21} is essentially the same as that in
\cite{odfa} and \cite{om} although the notion of quasi-Gibbs
measures was not introduced in those papers and the statement was
different to Lemma~\ref{l21}.
For the reader's convenience, we present the proof here.

We begin by recalling the definition of Dirichlet forms and related notions
according to \cite{fot} and \cite{mr}.
Let $ X $ be a Polish space
and $ m $ be a $ \sigma$-finite Borel measure on $ X $ whose
topological support equals $ X $. Let $ \mathcal{F}$ be a dense
subspace of $ L^2(X,m) $ and $ \mathcal{E}$ be a nonnegative bilinear
form defined on $ \mathcal{F}$.
We call $ (\mathcal{E},\mathcal{F}) $ a Dirichlet form on
$ L^2(X,m) $ if $ (\mathcal{E},\mathcal{F}) $ is closed and
Markovian. Here we say $ (\mathcal{E},\mathcal{F}) $ is Markovian if
$ \bar{u} : = \min\{ \max\{ u,0 \},1 \} \in\mathcal{F} $ and
$ \mathcal{E}(\bar{u},\bar{u})\le\mathcal{E}(u,u) $ for any
$ u \in\mathcal{F} $. The triplet
$(\mathcal{E},\mathcal{F},L^2(X,m))$ is called a Dirichlet space. We
say $(\mathcal{E},\mathcal{F},L^2(X,m))$ is local
if $ \mathcal{E}(u,v) = 0 $ for any $ u,v \in\mathcal{F}$
with disjoint compact supports.
Here a support of $ u \in\mathcal{F}$ is the topological support of
the signed measure $ u\,dm $; see \cite{fot}.

For a given Dirichlet space, there exists an $ L^2 $-Markovian
semi-group associated with the Dirichlet space. If the Dirichlet space
satisfies the quasi-regularity explained below, then there exists a
Hunt process associated with the Dirichlet space. Moreover, if the
Dirichlet form is local, then the Hunt process becomes a diffusion;
that is, a strong Markov process with continuous sample paths.

We say a Dirichlet space $(\mathcal{E},\mathcal{F},L^2(X,m))$
is quasi-regular if:
\begin{longlist}[(Q1)]
\item[(Q1)] There exists an increasing sequence of compact sets
$ \{ K_n \}$ such that
$ \bigcup_{n}\mathcal{F}(K_n) $
is dense in $ \mathcal{F} $ w.r.t. $ \mathcal{E}_{1}^{1/2} $-norm.
Here $ \mathcal{F}(K_n) = \{ f \in\mathcal{F} ; f = 0\
m\mbox{-a.e.} \mbox{on }K_n^c\}$, and
$ \mathcal{E}_{1}^{1/2}(f) = \mathcal{E}(f,f)^{1/2}+ \|f\|_{L^2(E
,m)} $.

\item[(Q2)] There exists a $ \mathcal{E}_{1}^{1/2} $-dense subset
of $ \mathcal{F}$ whose elements have $ \mathcal{E}$-quasi continuous
$ m $-version.

\item[(Q3)]
There exists a countable set $ \{ u_n \}_{n\in\N}$
having $ \mathcal{E}$-quasi continuous
$ m $-version~$ \tilde{u}_n $,
and an exceptional set $ \mathcal{N} $ such that
$ \{ \tilde{u}_n \}_{n\in\N}$ separates the points of
$ E \setminus\mathcal{N}$.
\end{longlist}

%
\begin{lem}\label{l31}
\textup{(1)} Assume \textup{(\hyperlink{A1}{A.1})}.
Let $(\Em, \domi)$ be as in (\ref{A30})
with $ a_{kl} = \delta_{kl} $.
Assume $(\Em, \domi)$ is closable on $ \Lm$.
Then its closure $(\Em, \dom)$ on $ \Lm$ is
a local, quasi-regular Dirichlet form.\vspace*{-6pt}
\begin{longlist}
\item[(2)] In addition, assume \textup{(\hyperlink{A0}{A.0})} and that
$(\Ema, \domai)$ is closable $ \Lm$.
Then its closure $(\Ema, \domai)$ on $ \Lm$ is
a local, quasi-regular Dirichlet form.
\end{longlist}
\end{lem}
\begin{pf}
(1) follows from \cite{odfa}, Theorem 1,
in which we suppose that
the density functions are locally bounded and
$ \sum_{m = 1}^{\infty} m \mu(\SSS_r^m) < \infty$.
We remark that these assumptions follow immediately from (\hyperlink{A1}{A.1}).
We have thus obtained (1).

Let $ \Ct\label{;31} = c_{\scriptsize\ref{;a33}}\sup|a_0(\mathsf{s},s )| $.
Then by (\hyperlink{A0}{A.0}), we see that $ c_{\scriptsize\ref{;31}}< \infty$ and
\[
 \doma\supset\dom, \qquad\Ema(f,f) \le c_{\scriptsize\ref{;31}} \Em(f,f) \qquad\mbox{for all } f
\in\dom .
\]
Hence, (2) follows from (1).
\end{pf}

We now proceed with the proof of closability.\vspace*{1pt}
Let $ \murm$ be as in Definition \ref{dfn1}.
We remark that $ \sum_{m=0}^{\infty} \murm= \mu$ by construction.
Let $ \Erm$ be the bilinear form defined by
%
%
\begin{equation}
\label{31} \Erm(f,g) = \int\mathbb{D}^{a}[f,g] \,d\murm .
\end{equation}
Then we have $ \Ema= \sum_{m = 1}^{\infty} \Erm$ for each $ r \in
\N$,
where $\Ema$ is the bilinear form given by (\ref{A30}).
We now quote a result from \cite{odfa}.
%
%
\begin{lem}[(Theorem 2 in \cite{odfa})] \label{l32}
Assume $ (\Erm, \domai)$ is closable on $ \Lm$ for all $ r,m \in\N$.
Then $ (\Ema, \domai) $ is closable on $ \Lm$.
\end{lem}
\begin{pf}
When $ b_r = r $ and the coefficient is the unit matrix,
Lemma~\ref{l32} was proved in Theorem 2 in \cite{odfa}. The
generalization to the present case is trivial.
\end{pf}
Let $ \murkm$ be as in Definition \ref{dfn1}.
Define the bilinear form $ \Ermk$ by
%
%
\begin{equation}
\label{31a} \Ermk(f,g) = \int\mathbb{D}^{a}[f,g] \,d\murkm .
\end{equation}
%

\begin{lem} \label{l33}
Assume $ (\Ermk, \domai)$ is closable on $ \Lmk$ for all $ k $.
Then $ (\Erm, \domai) $ is closable on $ \Lm$.
\end{lem}
\begin{pf}
By (\ref{qg1}), we have $\murkm\le\mu$.
This implies $(\Ermk, \domai)$ is closable,
not only on $ \Lmk$, but also on $ \Lm$. We deduce from (\ref{qg1}) that
the forms $ \{(\Ermk, \domai)\}$ are nondecreasing in $ k $ and
converge to
$(\Erm, \domai)$ as $ k\to\infty$.
Hence, $(\Erm,\domai)$ is closable on $\Lm$
according to the monotone convergence theorem of
closable bilinear forms.~%
\end{pf}

Let $ \murky$ be as in (\ref{qg4}). Let
$\Ermky(f,g) =
\int_{\SSS} \mathbb{D}^{a}[f,g] \,d\murmy$.
By (\ref{qg6}) and~(\ref{31})
%
\begin{eqnarray}
\label{32a} \Ermk(f,g) &=& \int_{\SSS} \Ermky(f,g) \murkm(d
\mathsf{s}),
\\
\label{32b}  \| f \|^2_{L^2(\SSS_r^{m}, \murkm)} &=& \int_{\SSS}
\| f \|^2_{L^2(\SSS_r^{m},\murmy)} \murkm(d\mathsf{s}) .
\end{eqnarray}

%
\begin{lem} \label{l34}
Assume $ (\Ermky, \domai)$ is closable on $ L^2(\SSS_r^{m}, \murmy) $
for $ \murkm$-a.s. $ \mathsf{s}$. Then $ (\Ermk, \domai)$ is closable
on $ \Lmk$.
\end{lem}
%
%
\begin{lem} \label{l35}
Assume \textup{(\hyperlink{A0}{A.0}), (\hyperlink{A2}{A.2})} and \textup{(\hyperlink{A3}{A.3})}.
Then $ (\Ermky, \domai)$ is closable on $ L^2(\SSS_r^{m}, \murmy) $
for $ \murkm$-a.s. $ \mathsf{s}$.
\end{lem}
Although the proof of Lemma~\ref{l34} is the same as that of Theorem~4 in~\cite{odfa},
we present it in Appendix \ref{sA1} for the reader's
convenience.
We also give the proof of Lemma~\ref{l35} in Appendix~\ref{sA1}.
We are now ready to prove the closability of $ (\Ema, \domai, \Lm) $.

%
\begin{lem}\label{l36}
Assume \textup{(\hyperlink{A0}{A.0}), (\hyperlink{A2}{A.2})} and \textup{(\hyperlink{A3}{A.3})}.
Then $ (\Ema, \domai, \Lm) $ is closable.
\end{lem}
\begin{pf}
By Lemmas~\ref{l32}--\ref{l35}, we conclude Lemma~\ref{l36}.
\end{pf}
\begin{pf*}{Proof of Lemma~\ref{l21}}
Lemma~\ref{l21} follows immediately from Lemmas~\ref{l31} and~\ref{l36}.
\end{pf*}

\begin{pf*}{Proof of Corollary \protect\ref{l21}}
By Lemma~\ref{l21} and \cite{fot}, Theorems~4.5.1, there exists
a $ \mu$-symmetric diffusion whose Dirichlet space
is $(\Ema,\doma,\break\Lm)$. Since $ 1 \in\doma$,
the diffusion is conservative,
which completes the proof.
\end{pf*}

\section{A sufficient condition of the quasi-Gibbs property}
\label{s4}
The most crucial assumption in Lemma~\ref{l21} is that of the quasi-Gibbs
property (\hyperlink{A2}{A.2}).
In this section, we introduce assumptions (\hyperlink{A4}{A.4}) and (\hyperlink{A5}{A.5}) below
to obtain a sufficient condition of (\hyperlink{A2}{A.2}).
These conditions guarantee that $ \mu$ has
a good finite-particle approximation
$\{\muN\}_{N\in\N}$ that enables us to prove the quasi-Gibbs property.
We set $ \Srhat= \{ x \in\SSSS ; |x|< r \} $ and $ \Dr^n =
\prod_{m=1}^{n} \{ |x_m|< r \} $.

%
(A.4)\hypertarget{A4}\
There exists a sequence of probability measures $\{ \muN\}_{N\in\N}$ on
$ \SSS$ satisfying the following.

\begin{longlist}[(1)]
\item[(1)]
The $ n $-correlation functions $ \rN$ of $ \muN$ satisfy
%
\begin{eqnarray}
\label{41a}\quad  \lim_{N \to\infty} \rN(x_1,\ldots,x_n) &=&
\rho^n (x_1,\ldots,x_n)\qquad \mbox{ a.e. for all $ n \in\N$,} %
\\
\label{41b}\quad\qquad  \sup_{N\in\N} \sup_{(x_1,\ldots,x_n) \in\tilde{\SSSS}_r^n } \rN(x_1,
\ldots,x_n) &\le&\bigl\{ c_{\scriptsize\ref{;70}} n^{\delta}\bigr
\}^n\qquad \mbox{for all $ n,r\in\N$} ,
\end{eqnarray}
where $ \Ct\label{;70}=c_{\scriptsize\ref{;70}}(r) >0$ and $ \delta= \delta
(r) <
1 $
are constants depending on $ r \in\N$.

\item[(2)] $ \muN(\mathsf{s}(\SSSS) \le n_{N}) = 1 $
for some $ n_{N}\in\N$.

\item[(3)]
$ \muN$ is a $ (\PhiN,\PsiN) $-canonical Gibbs measure.

\item[(4)]
The potentials
$ \map{\PhiN}{\SSSS}{\R\cup\{ \infty\}} $ and
$ \map{\PsiN}{\SSSS\ts\SSSS}{\R\cup\{ \infty\}} $
satisfy the following:
%
%
\begin{eqnarray}
\label{41c} \limi{N} \PhiN(x) &=& \Phi(x)\qquad \mbox{for a.e. $ x $,}\qquad
\inf_{N\in\N} \inf_{x \in\SSSS}\PhiN(x) > -\infty,
\\
\label{41d} \limi{N} \PsiN&=& \Psi\qquad\mbox{compact uniformly in } C^1
\bigl(\SSSS\ts\SSSS\setminus\{ x=y \}\bigr) ,
\nonumber
\\[-8pt]
\\[-8pt]
\nonumber
\qquad\inf_{N\in\N} \inf_{ x,y \in\Sr}\PsiN(x,y) &>& -
\infty\qquad \mbox{for all }r\in\N .
\end{eqnarray}
\end{longlist}

%
\begin{rem}\label{r22}
(1)
By (\ref{41a}) and (\ref{41b}), we see that
$ \limi{N} \muN= \mu$ weakly in~$ \SSS$ (see Lemma~\ref{l111}).
By $ \muN(\mathsf{s}(\SSSS) \le n_{N}) = 1 $, the DLR equation
(\ref{qg5}) makes sense even if $ \PsiN$ is a logarithmic function.
(\ref{41d}) implies the core $ \mathsf{\Gamma}$ in (\hyperlink{A3}{A.3}) becomes
$ \mathsf{\Gamma}= \{ 0 \} $ or $ \varnothing$.\vspace*{-6pt}
\begin{longlist}[(2)]
\item[(2)] By assumption, for each $ r \in\N$,
$ \PsiN\in C^1( \Srhat\ts\Srhat\setminus\{ x=y \}) $
for all sufficiently large $\NN$, and $ \Psi\in C^1(\SSSS\ts\SSSS
\setminus\{ x=y \})$.
We note that $ \PsiN$ is not necessarily in $ C^1(\SSSS\ts\SSSS
\setminus\{ x=y \})$.
\end{longlist}
\end{rem}

The difficulty in treating the logarithmic interaction
is the unboundedness at infinity.
Indeed, the DLR equation does not make sense for infinite volume.
The key issue in overcoming this difficulty
is the fact that the logarithmic functions have
small variations at infinity.
With this property, we can control the difference of interactions
rather than the interactions themselves.
Bearing this in mind, we introduce the set $ \Hrk$ in (\ref{41f})
and the assumption (\hyperlink{A5}{A.5}) below.

For $ \{\Sr\}$ in (\ref{qg0}), we set
$ \Srs= \Ss\setminus\Sr$ and $ \SSSS_{r\infty}= \Sr^c $.
For $r<s \le t < u \le\infty$ and
$ \mathsf{x} = \sum\delta_{x_i}, \mathsf{y} = \sum\delta_{y_j}\in\SSS$, we set
%
%
\begin{equation}
\label{41e} \PsiNrstu(\mathsf{x},\mathsf{y}) = \sum
_{x_i\in\SSSS_{rs} , y_j \in\SSSS_{tu}} \PsiN(x_i,y_j) 
.
\end{equation}
We write $\PsiN_{r,st} = \PsiN_{0r,st}$ and
$ \PsiN_{r,rs}(\mathsf{x},\mathsf{y}) = \PsiN_{r,rs}(x,\mathsf{y})
$ if
$ \mathsf{x}=\delta_{x}$.
We set $ \Hrk$ by
%
\begin{equation}
\label{41f}  \Hrk= \biggl\{ \mathsf{y}\in\mathsf{S} ; \2 \frac{|\PsiN_{r,rs}(x,\mathsf{y}) -\PsiN_{r,rs}(w,\mathsf{y})|} {
|x-w|}
\le{k}\biggr\} .
\end{equation}
The following is a tightness condition on $\{\muN\}$
according to 
$ \PsiN$.\vadjust{\goodbreak}

(A.5)\hypertarget{A5}\ The measures $ \{ \muN\} $ satisfy the following:
%
\begin{equation}
\label{41g}  \limi{k} \limsup_{ \NN\to\infty} \muN\bigl(\Hrk^c
\bigr) = 0\qquad \mbox{for all $ r \in\N$.}
\end{equation}

%
\begin{thmm} \label{l41}
Assume \textup{(\hyperlink{A4}{A.4})} and \textup{(\hyperlink{A5}{A.5})}.
Then $ \mu$ is a $ (\Phi, \Psi) $-quasi Gibbs measure.
\end{thmm}
We will prove Theorem~\ref{l41} in Section~\ref{s5}.

\begin{Cora}
Assume \textup{(\hyperlink{A0}{A.0})}, \textup{(\hyperlink{A1}{A.1})} and
\textup{(\hyperlink{A3}{A.3})--(\hyperlink{A5}{A.5})}.
Then we have the following:
\begin{longlist}[(1)]
\item[(1)]
$ (\Ema, \domai, \Lm) $ is closable, and its closure
$ \Dspa$ is a local, quasi-regular Dirichlet space.
\item[(2)]There exists a $ \mu$-reversible diffusion
$(\mathsf{X},\mathsf{P})$ associated with the Dirichlet space $ \Dspa$.
\end{longlist}
\end{Cora}

Let $ \SSS_r^{m} $ be as in (\ref{qg0}).
Using the set $ \Hrk$, we introduce cut-off measures $ \mukNm$,
%
\begin{equation}
\label{42a} \mukNm= \muN\bigl(\cdot\cap\SSS_r^{m} \cap
\Hrk\bigr) .
\end{equation}

We will prove Theorem~\ref{l41} along this sequence $ \{ \mukNm\} $.
For this, we first note the following.

%
\begin{lem} \label{l42}
There exists a weak convergent subsequence of $ \{ \mukNm\}$,
denoted by the same symbol, with limit measures $ \{ \murkm\}$
satisfying (\ref{qg1}) for all $ r,k,m $.
\end{lem}
\begin{pf}
Recall that $ \{ \muN\}$ is a weak convergent sequence.
This combined with $ \mukNm\le\muN$
shows that $ \{\mukNm\}$\vspace*{-1pt} is relatively compact for each $ r,k,m \in\N$.
Hence, we can choose a convergent subsequence\vspace*{-1pt} $ \{\mu^{n_N(r,k),m}_{r,k} \}$
from any subsequence of $ \{\mukNm\}$ for each $ r,k,m $.
Then by diagonal argument, we obtain a weak convergent subsequence\vspace*{-1pt} with limit
$ \{ \murkm\}$.

Since $\Hrk\subset\HH_{r,k+1}$, we have $\mukNm\le\mu^{N,m}_{r,k+1} $
by (\ref{42a}).
This allows us to deduce $\murkm\le\mu_{r,k+1}^{m}$,
which is the first claim of (\ref{qg1}).
Because of the weak convergence, we see that for $ f \in C_b(\SSS) $,
\begin{eqnarray*}
&& \biggl|\int f \,d \murkm- \int f \,d \murm\biggr|
\\
&&\qquad\le\limi{N} \biggl|\int f \,d \murkm- \int f \,d \mukNm\biggr| +
\limsup_{N\to\infty}\biggl |\int f \,d \mukNm- \int f \,d \muNm\biggr|
\\
&&\qquad\quad{} + \limi{N} \biggl|\int f \,d \muNm- \int f \,d \murm\biggr|
\\
&&\qquad= \limsup_{N\to\infty} \biggl|\int f \,d \mukNm- \int f \,d
\muNm\biggr|
\\
&&\qquad\le\Bigl\{\sup_{\mathsf{s}}\bigl|f(\mathsf{s})\bigr|\Bigr\}\cdot
\limsup_{N\to\infty} \muNm\bigl(\{ \Hrk\}^c \bigr) .
\end{eqnarray*}
By (\ref{41g}) we deduce that the right-hand side converges to zero
as $ k \to\infty$, which is the second claim of (\ref{qg1}).
We thus see that the limit measures $ \{ \murkm\}$ satisfy (\ref{qg1}).
\end{pf}

Let $ \muABN$ denote the conditional probability of
$ \mukNm$ defined by
\[
\muABN(d\mathsf{x}) = \mukNm\bigl(\pi_{ \Sr} \in d \mathsf{x} |
\pi_{ \Srs} (\mathsf{s})\bigr) .
\]
We note that, although\vspace*{-1pt} $ \mukNm$ is not necessarily a probability measure,
we take the normalizing in such a way that the conditional measure\vspace*{-1pt} $
\muABN$
to be a probability measure. As a result, we have $ \muABN(\SSS) = 1$ and
%
\begin{equation}
\label{42e} \mukNm\circ\pi_{ \Sr}^{-1}(d\mathsf{x}) = \int
_{\SSS}\muABN(d\mathsf{x}) \mukNm\circ\pi_{ \Srs
}^{-1}(d
\mathsf{s}) .
\end{equation}
Recall that by (\hyperlink{A4}{A.4}), $ \muN$ is a $ (\PhiN,\PsiN)
$-canonical Gibbs measure.
Then $ \muN$ satisfies the DLR equation (\ref{qg5}).
Hence, $ \muABN$ is absolutely continuous w.r.t. $ e^{-\mathcal
{H}_{r}^{N}(\mathsf{x})}\mA$.
Therefore, we denote its density by $ \sigma\ABN$.
Then by definition, we have for $\mukNm$-a.e. $\mathsf{s}$,
%
%
\begin{equation}
\label{42g}\qquad \sigma\ABN(\mathsf{x}) e^{-\mathcal{H}_{r}^{N}(\mathsf
{x})}\mA = \muABN(d\mathsf{x})
\qquad\mbox{where } \mathcal{H}_{r}^{N} = \mathcal{H}_{\Sr}^{\Phi^{N},\PsiN}
.
\end{equation}

The quasi-Gibbs property consists of two conditions:
(\ref{qg1}) and (\ref{qg2}). We have already proved (\ref{qg1})
by Lemma~\ref{l42}.
Therefore, it only remains to prove (\ref{qg2}).
This task is the most difficult part of the proof, and it is carried
out in the next section.
In the rest of this section, we explain
the strategy of the proof of (\ref{qg2}).

By taking the representation (\ref{42e}) into account,
the proof consists of two kinds of limit procedures:
(\ref{43a}) $ N\to\infty$ and then (\ref{43b}) $ s \to\infty$,
which involve the following convergence:
%
\begin{eqnarray}
 \label{43a} \limi{N} \muABN&= &\muAB, \qquad\limi{N} \mukNm\circ\pi_{ \Srs}^{-1}
= \murkm\circ\pi_{ \Srs}^{-1} ,
\\
 \label{43b} \limi{s} \muAB&=& \muA .
\end{eqnarray}
Note that two of these are the convergence of the \textit{conditional}
measures.
Comparing with the weak convergence of $ \{ \mukNm\}$ in Lemma~\ref{l42},
it is noted that
the convergence of conditional measures is much more delicate.
It involves a kind of strong convergence of the conditioned variable $
\mathsf{s}$.

In each step, we prove the bounds of the densities being uniform in $
N, s$
[(\ref{52c}) and (\ref{57j})]
and the related quantities as well as the convergence of measures as above.
The uniformity of the bounds is the crucial point of the proof.
We emphasize that we can carry out the proof because we treat the
cut-off measures
$ \{ \mukNm\} $ defined by (\ref{42a}). This cut-off is done by the
set $ \Hrk$.
Therefore, assumption (\hyperlink{A5}{A.5}) plays a significant role in the
proof of Theorem~\ref{l41}.

The first step consists of three lemmas.
Recall expressions (\ref{42e}) and (\ref{42g}).
We prove the uniform bounds\vspace*{-1pt} of
$ \int_{\SSS_r^{m} } e^{-\mathcal{H}_{r}^{N}(\mathsf{x})}\mA$
(Lemma~\ref{l51}) and $ \sigma\ABN$ (Lemma~\ref{l55}).\vspace*{-1pt}
We then prove the weak convergence
$ \limi{N} \mukNm\circ\pi_{ \Srs}^{-1} = \murkm\circ\pi_{ \Srs}^{-1}$
and the $ L^1$ convergence of their densities (Lemma~\ref{l56}).

The second step consists of two lemmas. In Lemma~\ref{l57}, we prove the
absolute continuity of the measures $ \muABs$ and the uniform bound
(\ref{57j}) of their densities $ \skABs(\mathsf{x})$.
Finally, in Lemma~\ref{l58} we prove the convergence of $ \skABs
(\mathsf
{x}) $ as $ s\to\infty$
using martingale convergence theorems to complete the proof of the
quasi-Gibbs property.

\section{\texorpdfstring{Proof of Theorem~\protect\ref{l41}}{Proof of Theorem 4.1}} \label{s5}
In this section, we prove (\ref{qg2}) to complete the proof of
Theorem~\ref{l41}.
We fix $ r , m \in\N$ throughout this section.
We divide this section into two parts.
In Section~\ref{s51}, we prove the first step (\ref{43a}), and
in Section~\ref{s52}, we prove the second step (\ref{43b}).

\subsection{Proof of the first step}\label{s51}
%
\begin{lem} \label{l51}
Set
\[
\Ct(n)\label{;43}= \sup_{n\le N\in\N} \max \biggl\{ \int_{\SSS_r^{m} }
e^{-\mathcal{H}_{r}^{N}(\mathsf{x})}\mA, \biggl[\int_{\SSS_r^{m} } e^{-\mathcal{H}_{r}^{N}(\mathsf{x})}\mA
\biggr]^{-1}\biggr\} .
\]
Then there exists an $ N_0 $ such that $ c_{\scriptsize\ref{;43}}(N_0) < \infty$.
\end{lem}
\begin{pf}
By (\hyperlink{A4}{A.4}), we see that
$ \sup\{e^{ -\mathcal{H}_{r}^{N}(\mathsf{x})} ; N\in\N, \mathsf
{x}\in\SSS_r^{m} \} <\infty$.
Hence, by~(\ref{41c}), (\ref{41d}) and the bounded convergence theorem,
we deduce that
\[
\lim_{N\to\infty} \int_{\SSS_r^{m} } e^{-\mathcal{H}_{r}^{N}(\mathsf{x})}
\mA= \int_{\SSS_r^{m}} e^{-\mathcal{H}_{r}(\mathsf{x})}\mA < \infty .
\]
Recall that $ \Phi(x) < \infty$ a.e. by assumption
[see the line after (\ref{2y})]
and $ \Psi(x,y)< \infty$ a.e. by the first assumption of (\ref{41d}).
Therefore, $ \mathcal{H}_{r}(\mathsf{x})<\infty$ a.e.
Hence, $ \int_{\SSS_r^{m}}e^{-\mathcal{H}_{r}(\mathsf{x})}\mA> 0 $.
Combining these completes the proof.
\end{pf}

We next consider\vspace*{-1pt} a decomposition of $ \sigma\ABN$ in (\ref{42g}).
By the DLR equation and (\ref{42a}), we deduce that
for $\mukNm$-a.e. $\mathsf{s} $,\vadjust{\goodbreak}
the density $\sigma\ABN$ is expressed in such a way that
%
%
\begin{equation}
\label{42h} \sigma\ABN(\mathsf{x}) = e^{ -\PsiNrs(\mathsf{x},\mathsf{s})} {\tCN}/{c_{\scriptsize\ref{;42}}^{N}(
\mathsf{s})} .
\end{equation}
Here $ \PsiNrs$ is given by (\ref{41e}).
We define $\tCN$ and $ \Ct^{N}\label{;42}(\mathsf{s})$ by
%
%
\begin{eqnarray}
\label{42i} \tCN&=&
1_{\SSS_r^{m} }(\mathsf{x}) \int_{\SSS}
1_{\Hrk}\bigl(\pi_{\Srs}(\mathsf{s})+ \mathsf{z}\bigr)
\nonumber
\\[-8pt]
\\[-8pt]
\nonumber
&&\hspace*{43pt}{}\times e^{ -\PsiNrsi(\mathsf{x},\mathsf{z})- \PsiNssi(\mathsf{s},\mathsf{z})} \mukNm\circ\pi_{\SSSS_{s\infty} }^{-1} (d\mathsf{z}) ,
\\
\label{42j} c_{\scriptsize\ref{;42}}^{N}(\mathsf{s}) &=& \int
_{\SSS} e^{ -\PsiNrs(\mathsf{x},\mathsf{s})} \tCN e^{-\mathcal{H}_{r}^{N}(\mathsf{x})} \mA .
\end{eqnarray}
We remark that, since $ \muN(\mathsf{s} (\SSSS)\le n_{N})= 1$,
$ \PsiNrsi$ and $ \PsiNssi$ are well defined
for $ \muN$-a.s. $ \mathsf{s}$.

Set $\Ct\label{;44} (k) = mk \cdot\operatorname{diam} (\Sr) $.
Then from (\ref{41e}) and (\ref{41f}), we deduce that
%
\begin{equation}
\label{52a}
\sup_{N\in\N} \sup_{r\le s < t \in\N} \sup_{\mathsf{x}, \mathsf{x}' \in\SSS_r^{m} }
\sup_{\mathsf{s}\in\Hrk} \bigl| \PsiNrst(\mathsf{x},\mathsf{s}) -\PsiNrst\bigl(
\mathsf{x}',\mathsf {s}\bigr) \bigr| \le c_{\scriptsize\ref{;44}}
\end{equation}
for each $ k\in\N$. Let
$ \mathsf{S}^{n}_{rs} =
\{\mathsf{x}\in\SSS; \mathsf{x}(\Srs) = n \}$.
Then from (\ref{41f}) and $ \Srs\subset\Ss$,
we deduce that
%
\begin{equation}
\label{52b}  \sup_{N\in\N} \sup_{r\le s < t \in\N} \sup_{\mathsf{y},\mathsf{y}'\in\mathsf{S}^{n}_{rs} }
\sup_{\mathsf{s}\in\Hsl} \biggl\{ \frac
{| \PsiN_{rs,st} (\mathsf{y},\mathsf{s}) -
\PsiN_{rs,st} (\mathsf{y}',\mathsf{s}) |} {
d_{\mathsf{S}^{n}_{rs}} (\mathsf{y},\mathsf{y}')} \biggr\} \le l
\end{equation}
for each $ n ,l \in\N$.
Here for $ \mathsf{s}, \mathsf{t}\in\mathsf{S}^{n}_{rs} $, we set
$ d_{\mathsf{S}^{n}_{rs}}(\mathsf{s},\mathsf{t}) = \min\sum_{i=1}^n |s_i-t_i|$,
where the minimum is taken over the labeling such that
$ \pi_{\Srs}(\mathsf{s}) = \sum_{i=1}^n \delta_{s_i} $ and
$ \pi_{\Srs}(\mathsf{t}) = \sum_{i=1}^n \delta_{t_i} $. %
Moreover, we used the inequality
$ \{ a_1+\cdots+a_n \}/\{ b_1+\cdots+b_n \}\le
\max\{ a_m/b_m ;m=1,\ldots,n \}$ for $ a_i \ge0 $ and $ b_j > 0 $.

%
\begin{lem} \label{l55}
Let $ \Ct\label{;44c} =e^{2c_{\fontsize{6.6}{6.6}\selectfont{\ref{;44}}}} c_{\scriptsize\ref{;43}}(N_0)$.
Then for $\mukNm$-a.e. $\mathsf{s} $, it holds that
%
%
\begin{equation}
\label{52c}
c_{\scriptsize\ref{;44c}}^{-1} \le\sigma\ABN(
\mathsf{x}) \le c_{\scriptsize\ref{;44c}}
\end{equation}
for all $ \mathsf{x} \in\SSS_r^{m} , r<s \in\N\mbox{ and
}N_0\le N\in\N$.
\end{lem}
\begin{pf}
By (\ref{42h}) and (\ref{52a}), we see that
%
\begin{equation}
\label{52e} \frac{\sigma\ABN(\mathsf{x})}{\sigma\ABN(\mathsf{x}')} = e^{-\PsiNrs(\mathsf{x},\mathsf{s}) +
\PsiNrs(\mathsf{x}' ,\mathsf{s})} \frac{\tCN}{\tCNN} 
\le e^{c_{\scriptsize\ref{;44}}} \frac{\tCN}{\tCNN} .
\end{equation}
By (\ref{42i}), we have for $\mukNm$-a.e. $\mathsf{s} $,
%
\begin{eqnarray}
\label{52d} && \mathop{\sup_{N\in\N, r< s\in\N}}_{\mathsf{x}, \mathsf{x}' \in\SSS
_r^{m} } \biggl\{ \frac{ \tCN}{\tCNN} \biggr\}
\nonumber\\
&&\quad= \mathop{\sup_{N\in\N, r< s\in\N}}_{\mathsf{x}, \mathsf{x}' \in\SSS
_r^{m} } \biggl\{  \frac{\int_{\SSS}1_{\Hrk}(\pi_{\Srs}(\mathsf{s})+
\mathsf{z})
e^{ -\PsiNrsi(\mathsf{x},\mathsf{z})- \PsiNssi(\mathsf{s},\mathsf{z})}
\mukNm\circ\pi_{\SSSS_{s\infty} }^{-1} (d\mathsf{z}) } {
\int_{\SSS}1_{\Hrk}(\pi_{\Srs}(\mathsf{s})+
\mathsf{z})
e^{ -\PsiNrsi(\mathsf{x}',\mathsf{z})- \PsiNssi(\mathsf{s},\mathsf{z})}
\mukNm\circ\pi_{\SSSS_{s\infty} }^{-1} (d\mathsf{z})
} \biggr\}
\nonumber
\\[-8pt]
\\[-8pt]
\nonumber
&&\quad= \mathop{\sup_{N\in\N, r< s < t \in\N}}_{\mathsf{x}, \mathsf{x}' \in
\SSS_r^{m} } \biggl\{ \frac{\int_{\SSS}1_{\Hrk}(\pi_{\Srs}(\mathsf{s}) + \mathsf{z})
e^{ -\PsiNrst(\mathsf{x},\mathsf{z})- \PsiN_{rs,st} (\mathsf
{s},\mathsf{z})}
\mukNm\circ\pi_{\SSSS_{s\infty} }^{-1} (d\mathsf{z}) }  {\int
_{\SSS}1_{\Hrk}\bigl(\pi_{\Srs}(\mathsf{s}) +
\mathsf{z}\bigr) e^{ -\PsiNrst(\mathsf{x}',\mathsf{z})- \PsiN_{rs,st} (\mathsf
{s},\mathsf{z})} \mukNm\circ\pi_{\SSSS_{s\infty} }^{-1}
(d\mathsf{z}) } \biggr\}
\\
\nonumber
&&\quad\le e^{c_{\scriptsize\ref{;44}}}\qquad  \mbox{by }(\ref{52a}) .
\end{eqnarray}
Here we used $ \muN(\mathsf{s} (\SSSS)\le n_{N})= 1$ for the third line.
We deduce from~(\ref{52e}) and~(\ref{52d}) that
\[
 \sup_{N\in\N} \sup_{r< s\in\N}\sup_{\mathsf{x}, \mathsf{x}'
\in\SSS_r^{m} } \bigl\{ {\sigma
\ABN(\mathsf{x})}/{\sigma\ABN\bigl(\mathsf{x}'\bigr)}\bigr\} \le
e^{2c_{\scriptsize\ref{;44}}} \qquad\mbox{for $\mukNm$-a.e. $\mathsf{s} $} .
\]
Hence for $\mukNm$-a.e. $\mathsf{s} $,
we see that for all
$ \mathsf{x}, \mathsf{x}' \in\SSS_r^{m} , r<s \in\N
, \mbox{and } N\in\N$,
%
\begin{equation}
\label{52g} e^{-2c_{\scriptsize\ref{;44}}} {\sigma\ABN\bigl(\mathsf{x}'\bigr)}
\le {\sigma\ABN(\mathsf{x})}\le e^{2c_{\scriptsize\ref{;44}}} {\sigma\ABN\bigl(
\mathsf{x}'\bigr)} .
\end{equation}
Multiply (\ref{52g}) by
$ 1_{\SSS_r^{m} }(\mathsf{x}')e^{-\mathcal{H}_{r}^{N}(\mathsf{x}')}$
and integrate w.r.t. $ \Lambda(d\mathsf{x}')$.
Note that by (\ref{42g}) we have
$ \int_{\SSS_r^{m} }{\sigma\ABN(\mathsf{x}')}e^{-\mathcal
{H}_{r}^{N}(\mathsf{x}')}
\Lambda(d\mathsf{x}') =1$.
Then we deduce that for $\mukNm$-a.e. $\mathsf{s} $,
\[
e^{-2c_{\scriptsize\ref{;44}}} \le{\sigma\ABN(\mathsf{x})} \int
_{\SSS_r^{m} }e^{-\mathcal{H}_{r}^{N}(\mathsf{x}')} \Lambda\bigl(d\mathsf{x}'
\bigr) \le e^{2c_{\scriptsize\ref{;44}}}\qquad \mbox{for all }\mathsf{x} \in\SSS_r^{m}
.
\]
This combined with Lemma~\ref{l51} yields (\ref{52c}).
\end{pf}
%

Let\vspace*{-1pt} $\mathcal{H}_{rs}^{N} =
\mathcal{H}_{\Srs}^{\Phi^{N},\PsiN}$ and
$\mathcal{H}_{rs} =
\mathcal{H}_{\Srs}^{\Phi,\Psi}$.
By (\ref{41a}) and (\ref{41b}), we see that
$ \mukNm\circ\pi_{ \Srs}^{-1}$ and
$ \murkm\circ\pi_{ \Srs}^{-1}$ are
absolutely continuous w.r.t. $ e^{-\mathcal{H}_{rs}^{N}}\Lambda$ and
$ e^{-\mathcal{H}_{rs}}\Lambda$, respectively.
Hence, we denote by $ \Delta^{N} $ and $ \Delta$
their Radon--Nikodym densities, respectively.
%
%
\begin{lem} \label{l56}
\textup{(1)} $ \mukNm\circ\pi_{ \Srs}^{-1}$
converges weakly to $ \murkm\circ\pi_{ \Srs}^{-1} $ as $ N\to
\infty$.
\begin{longlist}[(2)]
\textup{(2)}
$ \Delta^{N} e^{-\mathcal{H}_{rs}^{N}} $ converges to
$ \Delta e^{-\mathcal{H}_{rs}}$ in $ \LABone$ as $ N\to\infty$.
\end{longlist}
\end{lem}
\begin{pf}
Let $ E $ be the discontinuity points of $ \pi_{ \Srs}$. Namely
\[
E = \Bigl\{ \mathsf{s}\in\SSS ; \lim_{n\to\infty}\pi_{\Srs}(
\mathsf{s}_n) \not= \pi_{\Srs}(\mathsf{s}) \mbox{ for some }
\{ \mathsf{s}_n \} \mbox{ such that } \limi{n} \mathsf{s}_n
=\mathsf{s}\Bigr\} .
\]
Then by (\hyperlink{A1}{A.1}), we deduce that $ \murkm(E)\le\mu(E)=0 $.
Since $ \mukNm$ converge weakly to $ \murkm$ by Lemma~\ref{l42} and
the discontinuity points of
$ \pi_{ \Srs}^{-1} $ are $ \murkm$-measure zero, we obtain (1).


We proceed with (2). It only remains to prove that
$\{\Delta^{N} e^{-\mathcal{H}_{rs}^{N}}\}_{N\in\N}$ is relatively
compact in $ \LABone$.
Indeed, if this property holds, then their limit points are unique and
equal to $ \Delta e^{-\mathcal{H}_{rs}}$ by (1).

Recall that $ \mathsf{S}^{n}_{rs} =
\{\mathsf{x}\in\SSS; \mathsf{x}(\Srs) = n \}$,
and note that
\[
\Delta^{N} e^{-\mathcal{H}_{rs}^{N}} = \Delta^{N} e^{-\mathcal{H}_{rs}^{N}}
\sum_{n=0}^{\infty}1_{\mathsf{S}^{n}_{rs}} .
\]
We deduce from (1)
that for each $ \epsilon>0 $ there exists an $ n_0 $ such that
%
\begin{equation}
\label{56u} \sup_{N\in\N} \mukNm
\Biggl(\sum
_{n=n_0}^{\infty} \mathsf{S}^{n}_{rs}
\Biggr) < \epsilon ,
\end{equation}
which is equivalent to
%
\begin{equation}
\label{56v} \sup_{N\in\N} \Biggl\| \Delta^{N} e^{-\mathcal{H}_{rs}^{N}} \sum
_{n=n_0}^{\infty}1_{\mathsf{S}^{n}_{rs}}
\Biggr\|_{ \LABone} < \epsilon .
\end{equation}
According to (\ref{56v}), the relative compactness of
$\{\Delta^{N} e^{-\mathcal{H}_{rs}^{N}}\}_{N\in\N}$ in $ \LABone$
follows from that of
$ \{ \Delta^{N} e^{-\mathcal{H}_{rs}^{N}}
1_{\mathsf{S}^{n}_{rs}} \}_{N\in\N}$ for each $ n \in\N$.
Hence, we fix $ n \in\N$ in the rest of the proof.

We set
$ \mu^{N}_{l} = \mukNm(\cdot\cap\mathsf{S}^{n}_{rs}\cap\Hsl)$,
where $ \Hsl$ is as in (\ref{41f}).
Let $ \Delta^{N}_{l} $ be the Radon--Nikodym density of
$ \mu^{N}_{l} \circ\pi_{ \Srs}^{-1}$
w.r.t. $ e^{-\mathcal{H}_{rs}^{N}}\Lambda$.
Since $ \mu^{N}_{l} \le\mukNm$, we see that
$
\Delta^{N}_{l}e^{-\mathcal{H}_{rs}^{N}}\le
\Delta^{N} e^{-\mathcal{H}_{rs}^{N}}
$.
Combining this with (\ref{41g}) yields
%
\begin{eqnarray}
\label{56a}
&&\lim_{l\to\infty} \limsup_{N\in\N} \bigl\|
\Delta^{N} e^{-\mathcal{H}_{rs}^{N}}- \Delta^{N}_{l}e^{-\mathcal{H}_{rs}^{N}}
\bigr\|_{\LABone}
\nonumber
\\[-8pt]
\\[-8pt]
\nonumber
&&\qquad \le \lim_{l\to\infty}\limsup_{N\in\N} \mukNm\bigl(
\Hsl^c\bigr) =0 .
\end{eqnarray}
According to (\ref{56a}), it only remains to prove the relative
compactness of
$\{\Delta^{N}_{l}e^{-\mathcal{H}_{rs}^{N}}\}_{N\in\N}$
in $ \LABone$ for each $ l\in\N$. Hence, we fix $ l \in\N$ in the
rest of the proof.

For $ q \in\N$ we set $ B_r^{q}=\{0< |s-\Sr|<1/q \} $. Let
%
\begin{equation}
\label{56x} \Aq= \bigl\{ \mathsf{s} \in\mathsf{S}^{n}_{rs}
\cap\Hsl ; \mathsf{s}\bigl(B_r^{q}\bigr)=0 \bigr\} .
\end{equation}
By definition, $ \Aq$ is the subset of $ \mathsf{S}^{n}_{rs}\cap
\Hsl$
with no particles in $ B_r^{q}$, where
$ B_r^{q}$ is the intersection of $ \Sr^c$ and
the $ 1/q$-neighborhood of $ \Sr$.
Then the relative compactness of
$\{\Delta^{N}_{l}e^{-\mathcal{H}_{rs}^{N}}\}_{N\in\N}$ follows from
that of
$ \{ \Delta^{N}_{l} e^{-\mathcal{H}_{rs}^{N}} 1_{\Aq} \}_{N\in\N}$
for all sufficiently large $ q\in\N$.
Indeed, by (\ref{41a})--(\ref{41d}),
for each $ \epsilon>0$ there exists a $ q_0\in\N$ such that, for all
$ q \ge q_0 $,
\begin{eqnarray*}
\sup_{N\in\N} \bigl\| \Delta^{N}_{l}e^{-\mathcal{H}_{rs}^{N}} -
\Delta^{N}_{l}e^{-\mathcal{H}_{rs}^{N}}1_{\Aq}
\bigr\|_{\LABone} &\le&\sup_{N\in\N} \mukNm\bigl((\Aq)^c
\bigr)
\\
& \le&\sup_{N\in\N} \int_{ B_r^{q}} \rNone(x) \,dx
\le\epsilon .
\end{eqnarray*}

Let $ \Ct\label{;56b}(q)$ be the constant defined by
%
\begin{equation}
\label{56b} c_{\scriptsize\ref{;56b}}(q) = \sup_{N\in\N} \sup_{ {\mathsf{x}\in\SSS
_r^{m} }}
\sup\biggl\{ \frac
{|\PsiNrs(\mathsf{x},\mathsf{y})-\PsiNrs(\mathsf{x},\mathsf{y}')|
}%
{ d_{\mathsf{S}^{n}_{rs}}\bigl(\mathsf{y},
\mathsf{y}'\bigr)} ; \mathsf{y} \not= \mathsf{y}' \in
\Aq\biggr\} .
\end{equation}
Then we have $ c_{\scriptsize\ref{;56b}} (q) < \infty$.
Note that $ \pi_{\Srs^c} = \pi_{\Sr} + \pi_{\SSSS_{s\infty} }$.
Hence we write $ \pi_{\Srs^c} (\mathsf{s}) = \mathsf{x} + \mathsf{z}$,
where $ \mathsf{x}\in\pi_{\Sr}(\SSS)$ and $ \mathsf{z}\in\pi_{\SSSS_{s\infty} }(\SSS)$.
With this notation, $ \Delta^{N}_{l}(\mathsf{y}) $ can be written as
\begin{eqnarray*}
 \Delta^{N}_{l}(\mathsf{y}) = c \int_{\SSS}
1_{\Hrk\cap\Hsl}\bigl( \mathsf{x} + \pi_{\Srs}(\mathsf{y})+\mathsf{z}
\bigr) e^{
-\PsiNrs(\mathsf{x},\mathsf{y})
- \PsiNssi(\mathsf{y},\mathsf{z}) } \mu^{N}_{l} \circ
\pi_{\Srs^c}^{-1} (d\mathsf{x}\,d\mathsf{z}) .
\end{eqnarray*}
Here $ c $ is a constant.
Then applying (\ref{56b}) and (\ref{52b}) to
$ \PsiNrs(\mathsf{x},\mathsf{y}) $ and $ \PsiNssi(\mathsf
{y},\mathsf{z})$, respectively,
we deduce that
%
\begin{equation}
\label{56c} \sup_{N\in\N} \sup_{\mathsf{y},\mathsf{y}'\in\Aq} \biggl\{
\frac{\Delta^{N}_{l}(\mathsf{y})}{\Delta^{N}_{l}(\mathsf{y}')} \biggr\}
\le e^{(c_{\scriptsize\ref{;56b}}(q)+l) \,d_{\mathsf{S}^{n}_{rs}}(\mathsf{y},\mathsf{y}')} 
.
\end{equation}
Taking the logarithm of (\ref{56c}) and interchanging the role
of $ \mathsf{y}$ and $ \mathsf{y}'$, we deduce that
%
\begin{equation}
\label{56d}\qquad \sup_{N\in\N} \sup_{\mathsf{y},\mathsf{y}'\in\Aq} \bigl\{ \bigl| \log{
\Delta^{N}_{l}(\mathsf{y})} - \log{\Delta^{N}_{l}
\bigl(\mathsf {y}'\bigr)} \bigr| \bigr\} 
\le {
\bigl(c_{\scriptsize\ref{;56b}}(q)+l\bigr) \,d_{\mathsf{S}^{n}_{rs}}\bigl(\mathsf{y},
\mathsf{y}'\bigr)} 
.
\end{equation}

We deduce from (\ref{56c}) and (\ref{56d}) that
$ \{\Delta^{N}_{l} (\mathsf{y})\}_{N\in\N} $ is equi-continuous
in $ \mathsf{y}$ on~$ \Aq$
for each $ q\in\N$. From the definition of $ \Delta^{N}_{l} $, we
see that
\[
\sup_{N\in\N} \bigl\| \Delta^{N}_{l} e^{-\mathcal{H}_{rs}^{N}}
1_{\Aq} \bigr\|_{ \LABone} < \infty .
\]
We deduce from (\ref{41c}) and (\ref{41d}) that
$ \limi{N} e^{-\mathcal{H}_{rs}^{N}} 1_{\Aq} =
e^{-\mathcal{H}_{rs}}1_{\Aq} $ in $ \LABone$,
and that the limit satisfies
$ \| e^{-\mathcal{H}_{rs}}1_{\Aq} \|_{\LABone} > 0 $, which implies
\[
\linfi{N} \bigl\| e^{-\mathcal{H}_{rs}^{N}}1_{\Aq} \bigr\|_{\LABone}> 0 .
\]
These allow us to deduce that
\[
\lsupi{N} \bigl\| \Delta^{N}_{l}1_{\Aq}
\bigr\|_{ L^{\infty}(\SSS,\Lambda)
} < \infty .
\]
We therefore apply the Ascoli--Arzel\'{a} theorem to $ \Delta^{N}_{l}1_{\Aq} $ to
deduce that $ \{ \Delta^{N}_{l}1_{\Aq} \} $ is relatively compact in
$ C_b(\Aq)$ with uniform norm.
Because of the uniform boundedness of
$\{ e^{-\mathcal{H}_{rs}^{N}} 1_{\Aq} \}_{N\in\N} $,
we see that $\{ \Delta^{N}_{l}e^{-\mathcal{H}_{rs}^{N}} 1_{\Aq} \}_{N\in\N}$
is relatively compact in $ \LABone$ for each $ q $.
Hence $\{ \Delta^{N}e^{-\mathcal{H}_{rs}^{N}} \}_{N\in\N}$
is relatively compact in $ \LABone$ because of (\ref{56a}).
Therefore, we complete the proof.
\end{pf}
%

\subsection{Proof of the second step} \label{s52}

%
\begin{lem} \label{l57}
Let $ \muABs= \murkm(\pi_{ \Sr}(\mathsf{s})
\in d \mathsf{x} | \pi_{ \Srs}(\mathsf{s}))$.
Then we have the following:
\begin{longlist}[(1)]
\item[(1)]
$ \muAB$ is absolutely continuous w.r.t. $ e^{- \mathcal
{H}_{r}(\mathsf{x})}\mA$ for \8.

\item[(2)]
For each $ r,m,k \in\N$,
the Radon--Nikodym densities $ \skABs$ of
$ \muABs$ in~\textup{(1)} satisfy for \8
and all $ s \in\N$ such that $ r<s $
%
%
\begin{equation}
\label{57j} c_{\scriptsize\ref{;44c}}^{-1}\le\skABs(\mathsf{x})\le
c_{\scriptsize\ref{;44c}} \qquad\mbox{for $ \muABs$-a.e. $ \mathsf{x}$} .
\end{equation}
\end{longlist}
\end{lem}
\begin{pf}
Similarly to Lemma~\ref{l56}(1), we see that
$ \mukNm\circ(\pi_{\Sr},\pi_{\Srs})^{-1} $
converge weakly to
$ \murkm\circ(\pi_{\Sr},\pi_{\Srs})^{-1} $
as $ N\to\infty$. Hence, for
$ \mathsf{f}$, $ \mathsf{g} \in C_b(\SSS) $,
we have
%
\begin{equation}
\label{57z} \int_{\SSS} \mathsf{f}\bigl(\pi_{\Sr}(
\mathsf{s})\bigr) \mathsf{g}\bigl(\pi_{\Srs}(\mathsf{s})\bigr) \,d\murkm
 = \limi{N} \int_{\SSS} \mathsf{f}\bigl(\pi_{\Sr}(
\mathsf{s})\bigr) \mathsf{g}\bigl(\pi_{\Srs}(\mathsf{s})\bigr) \,d\mukNm
.
\end{equation}

By Lemma~\ref{l55} and the diagonal argument,
there exist subsequences of $\{\sigma\ABN\}_N $,
denoted by the same symbol, with a limit
$\skAB$ such that for all $ k,m, r<s\in\N$,
%
\begin{eqnarray}
\label{57a} \qquad\limi{N} \skABNs\bigl(\pi_{ \Sr}(\mathsf{s})\bigr) =
\skABs\bigl(\pi_{ \Sr}(\mathsf{s})\bigr) \qquad\mbox{$ * $-weakly in $\LAB$} .
\end{eqnarray}
Here $ \skABs$ is a function such that
$ \skABs(\mathsf{x}) = \sigma\ABss(\pi_{ \Sr}(\mathsf{x}))$.
Let
%
\begin{eqnarray}
\label{57h} \mathsf{F}^{N}(\mathsf{s})& =& \mathsf{f}\bigl(
\pi_{\Sr}(\mathsf{s})\bigr) \mathsf{g}\bigl(\pi_{\Srs}(
\mathsf{s})\bigr) \Delta^{N} (\mathsf{s}) e^{-\mathcal{H}_{r}^{N}(\mathsf{s})} ,
\\
\label{57i} \mathsf{F}(\mathsf{s}) &= &\mathsf{f}\bigl(\pi_{\Sr}(
\mathsf{s})\bigr) \mathsf{g}\bigl(\pi_{\Srs}(\mathsf{s})\bigr) \Delta(
\mathsf{s}) e^{-\mathcal{H}_{r}(\mathsf{s})} .
\end{eqnarray}
Then by Lemma~\ref{l56}(2), we see that $ \mathsf{F}^{N} $
converge to
$ \mathsf{F} $ in $ \LABone$. This combined with (\ref{57a}) implies
%
%
\begin{equation}
\label{57g}  \limi{N}\int_{\SSS} \mathsf{F}^{N}(
\mathsf{s}) \skABNs(\mathsf{s})\,d\Lambda = \int_{\SSS}
\mathsf{F}(\mathsf{s}) \skABs(\mathsf{s}) \,d\Lambda .
\end{equation}

By (\ref{57z}), (\ref{57g}) and
$ \Delta(\mathsf{y}) e^{-\mathcal{H}_{r}(\mathsf{y})} \Lambda
(d\mathsf{y}) =
\murkm\circ\pi_{ \Srs}^{-1}(d\mathsf{y}) $,
we obtain
\[
 \int_{\SSS} \mathsf{f}(\mathsf{x}) \mathsf{g}(\mathsf{y})
\,d\murkm = \int_{\SSS} \mathsf{f}(\mathsf{x} ) \mathsf{g}(
\mathsf{y}) \skABs(\mathsf{x}) e^{-\mathcal{H}_{r}(\mathsf{x})} \mA\murkm\circ
\pi_{ \Srs}^{-1}(d\mathsf{y}) ,
\]
where $ \mathsf{x} = \pi_{\Sr}(\mathsf{s})$ and
$ \mathsf{y} = \pi_{\Srs}(\mathsf{s})$.
Hence, we obtain (1) with density $ \skABs$.

By (\ref{52c}) and (\ref{57a}),
we see that $\skABs$ satisfies (\ref{57j}), which implies (2).
\end{pf}

%
\begin{lem} \label{l58}
Let $ \muA(d\mathsf{x}) $ be as in (\ref{qg4}).
Let $ \skAB$ be as in Lemma~\ref{l57}. Then the following limit exists:
%
%
\begin{equation}
\label{58a} \sigma_{r,k,\mathsf{s}}^{ m }(\mathsf{x}) : = \limi{s} \skAB
(\mathsf{x}) \qquad\mbox{for $ \muA$-a.s. $\mathsf{x}$, for $ \murkm$-a.s. $ \mathsf{s}
$} .
\end{equation}
Moreover, $ \sigma_{r,k,\mathsf{s}}^{ m }$ satisfies for \8,
%
\begin{eqnarray}
\label{58c} c_{\scriptsize\ref{;44c}}^{-1}\le\sigma_{r,k,\mathsf{s}}^{ m }(
\mathsf{x})&\le &c_{\scriptsize\ref{;44c}} \qquad\mbox{for $ \muA$-a.e.
}\mathsf{x},
\\
\label{58b} \sigma_{r,k,\mathsf{s}}^{ m }(\mathsf{x}) e^{-\mathcal
{H}_{r}(\mathsf{x})}
\mA&=& \muA(d\mathsf{x}) .
\end{eqnarray}
\end{lem}
\begin{pf}
Define $ \map{M_s}{\mathsf{S}}{\R} $ by
$ M_s (\mathsf{s}) = \skAB(\mathsf{x})$, where
$ \mathsf{x} = \pi_{\Sr}(\mathsf{s})$.
Recall that
$ \skAB$ is the Radon--Nikodym density of
$ \muAB$ w.r.t. $ e^{- \mathcal{H}_{r}(\mathsf{x})} \mA$ and that
$ \muAB= \mu^{m}_{r,k,\pi_{\Srs}(\mathsf{s}),rs}$ by construction.
Hence,
%
%
\begin{equation}
\label{57} M_s (\mathsf{s})e^{-\mathcal{H}_{r}(\mathsf{x}) }\mA =
\mu^{m}_{r,k,\pi_{\Srs}(\mathsf{s}),rs} (d\mathsf{x}) .
\end{equation}

Let
$ \mathcal{F}_s =
\sigma[\pi_{\Sr}, \pi_{\Srs}] $, where
$ r< s \le\infty$.
Then by (\ref{57}), we see that
$ \{ M_s \}_{s \in[r, \infty)}$
is an $ (\mathcal{F}_s) $-martingale, which implies
$ M_{\infty}(\mathsf{s}): =
\limi{s} M_s (\mathsf{s})$
exists for $ \murkm$-a.e. $\mathsf{s}$.
Since
\[
M_s (\mathsf{s}) = \sigma^m_{r,k,\pi_{\Srs}(\mathsf{s}),rs} (
\mathsf{x})\qquad \mbox{where } \mathsf{x} = \pi_{\Sr}(\mathsf{s}) ,
\]
we write
$ M_{\infty}(\mathsf{s}) = \sigma_{r,k,\mathsf{s}}^{ m }(\mathsf
{x}) $.
By construction,
$ \sigma_{r,k,\mathsf{s}}^{ m }(\mathsf{x}) =
\sigma_{r,k,\pi_{\SSSS_{r\infty}}(\mathsf{s})}^{ m } (\mathsf{x})
=\break
\sigma_{r,k,\pi_{\Sr^c}(\mathsf{s})}^{ m } (\mathsf{x}),
$ and, for $ \murkm$-a.s. $ \mathsf{s} $, we can regard
$ \sigma_{r,k,\mathsf{s}}^{ m }(\mathsf{x})$ as a $ \sigma[ \pi_{\Sr}]$-measurable function in $ \mathsf{x}$.
Hence, through disintegration (\ref{qg6}), we obtain (\ref{58a}).

We immediately obtain (\ref{58c})
from (\ref{57j}) and (\ref{58a}).

We see that
$ \{ M_s \}_{s \in[r, \infty)}$ is uniformly integrable by (\ref{57j}).
Hence we deduce from (\ref{58a}) that $ M_s (\mathsf{s})$ converges\vspace*{1pt} to
$ M_{\infty}(\mathsf{s}) = \sigma_{r,k,\mathsf{s}}^{ m }(\mathsf{x})$
strongly in $ L^1(\SSS_r^{m}, \muA) $,
which combined with (\ref{57}), and the definition
$ M_s (\mathsf{s}) =\break \skAB(\mathsf{x}) $ yields
(\ref{58b}).
\end{pf}

\begin{pf*}{Proof of Theorem~\ref{l41}}
By Lemma~\ref{l42}, we see that $ \{ \murkm\}$ satisfies (\ref{qg1}).
Moreover, by (\ref{58c}) and (\ref{58b}) we deduce
that $ \murky$ satisfies (\ref{qg2}),
which completes the proof of Theorem~\ref{l41}.
\end{pf*}

\section{\texorpdfstring{A sufficient condition of (\protect\hyperlink{A5}{A.5})}{A sufficient condition of (A.5)}}\label{s6}
In this section, we give a sufficient condition of (\hyperlink{A5}{A.5}) when $
\Psi$ is a logarithmic function
and $ d= 1,2$. When $ d=2$, we regard $ \R^2$ as $ \mathbb{C}$.
We assume 
%
\begin{equation}
\label{61a} \Psi(x,y) = - \beta\log|x -y |\qquad (\beta\in\R) .
\end{equation}

We take $ \PsiN$ in two different ways. In the first case we assume $
d=1,2$ and
$ \PsiN= \Psi$ for all $ N$, while in the second case $ \PsiN$
depend on $ N$.
To unify these two cases, we introduce
%
%
\begin{equation}
\label{61s} \PsiN(x,y) = - \beta\log\bigl|\wNx-\wNy\bigr| .
\end{equation}
We set for the first case $ d=1,2$ and
%
\begin{equation}
\label{61t} \wNx= x .
\end{equation}
Next we let $ \IN=(-N,N)$ and $ n_{N}= 2^{4\NN}$.
For the second case,
we set $ d=1$ and define the map $ \map{\wN}{\SSSS}{\mathbb{C}}$ by
%
\begin{equation}
\label{61u}\wNx= \cases{\displaystyle \ii\frac{n_{N}}{2\pi} \bigl(1-e^{2\pi\ii x /n_{N}}\bigr)
,&\quad  $\mbox{for }x \in\IN,$ \vspace*{2pt}
\cr
x, & \quad $\mbox{for } x\in
\INN^c,$ \vspace*{2pt}
\cr
\mbox{linear interpolation}, &\quad $\mbox{for } x
\in\INN\setminus\IN.$}
\end{equation}
By construction, we have $ \wN(0)= 0$,
\[
\Re\bigl[\wNx\bigr]= - \Re\bigl[\wN(-x)\bigr]\quad \mbox{and}\quad  \Im\bigl[\wNx\bigr]= \Im\bigl[
\wN(-x)\bigr].
\]
Here $ \Re[ \cdot ]$ and $ \Im[ \cdot ]$
denote the real and imaginary part of $ \cdot $, respectively.
It is easy to see that $ |\wNx| < |\wN(y)| $ for $ |x|< |y|$.
We note that $ \wN$ maps $ \IN$ into a subset of the circle in $
\mathbb{C}$
centered at $ \ii\frac{n_{N}}{2\pi}$ with radius $ \frac
{n_{N}}{2\pi}$.
We take $ n_{N}= 2^{4\NN} $ such that it is large compared with $ N $,
which converges the trajectory of $ \wN(\R)$ to the real axis rapidly as
$ N \to\infty$.

In the former case, (\ref{61t}) is used for the Ginibre random point
field (Theorem~\ref{l23}). We will use this choice to prove the quasi-Gibbs
property of the Bessel random point field in a forthcoming paper.
In the latter case, (\ref{61u}) is used for Dyson's model
(Theorem~\ref{l22}), where we use circular ensembles, and thus, the
above choice of
$ \wN$ is suitable.

The argument in this section may be generalized to higher dimensions $
d \ge3$.
We restrict ourselves to the case $ d=1,2$. As a result, we obtain a rather
simple expression of the Tayler expansion of $ \Psi(\wNx,\wNy) $.
We remark that $ z/|z|^2 = 1/\bar{z}\in\mathbb{C}$.
%
%
\begin{lem} \label{l61}
Assume (\ref{61a}). Let $ x,y\in\R$ such that $ |\wNx|<|\wN(y)|
$. Then
%
\begin{equation}
\label{61d} \Psi\bigl(\wNx,\wNy\bigr) -\Psi\bigl(0,\wNy\bigr) = \beta\sum
_{\ell= 1}^{\infty} \frac{1}{\ell} \Re\biggl[\biggl(
\frac{\bar{\varpi}_{N}(x)}{\bar{\varpi}_{N}(y)}\biggr)^{\ell}\biggr] .
\end{equation}
Here $ \bar{\varpi}_{N}$ denotes the complex conjugate of $ \wN$.
\end{lem}
\begin{pf}
Let $ r = |\wNx|/|\wNy| $ and $ \theta= \angle(\wNx,\wNy) $.
Then 
\begin{eqnarray*}
\Psi\bigl(\wNx,\wNy\bigr) -\Psi\bigl(0,\wNy\bigr)& =& -\frac{\beta}{2} \log\biggl|
\frac{\wNx}{|\wNy|}- \frac{\wNy}{|\wNy|}\biggr|^2
\\
& =& -\frac{\beta}{2} \log \bigl(1+r^2- 2r \cos\theta
\bigr)
\\
& =& -\frac{\beta}{2} \bigl\{\log\bigl(1-re^{ \ii\theta}\bigr) +
\log\bigl(1-re^{- \ii
\theta}\bigr)\bigr\} .
\end{eqnarray*}
Hence, (\ref{61d}) follows from the Tayler expansion.
\end{pf}
%

%
\begin{rem}\label{r61}
When (\ref{61t}) holds, we easily deduce from Lemma~\ref{l61}
that for $ 0 < |x|<|y| $
%
\begin{eqnarray}
\label{61c}\Psi(x,y) -\Psi(0,y)& =& \beta\sum_{\ell= 1}^{\infty}
\frac{1}{\ell}\biggl(\frac{x}{y}\biggr)^{\ell} \qquad\mbox{if $
\SSSS= \R$},
\\
\label{61b} \Psi(x,y) -\Psi(0,y) &=& \beta\sum_{\ell= 1}^{\infty}
\frac{1}{\ell} \Re\biggl[\biggl(\frac{\bar{x}}{\bar{y}}\biggr)^{\ell}
\biggr] \qquad\mbox{if $\SSSS= \mathbb{C} $ }.
\end{eqnarray}

\end{rem}

Let $ \Srs= \Ss\setminus\Sr= \{ y\in\SSSS ; b_r \le|y|< b_s
\} $ as before, where $ \Sr$ and $ b_r $ are given by (\ref{qg0}).
We set
$ \PsiN_{rs}(x,\mathsf{y}) = \sum_{y_i\in\Srs} \PsiN(x,y_i) $,
where $ \mathsf{y}=\sum_i\delta_{y_i}$. By (\ref{61d}),
\begin{eqnarray*}
\PsiN_{rs}(x,\mathsf{y}) - \PsiN_{rs}(w ,
\mathsf{y}) & =& \beta\sum_{\ell= 1}^{\infty}
\frac{1}{\ell} \4 \Re\biggl[\frac{\bar{\varpi}_{N}(x)^{\ell} - \bar{\varpi}_{N}(w)^{\ell
} }{\bar{\varpi}_{N}(y_i)^{\ell} }\biggr]
\\
& = &\beta\sum_{\ell= 1}^{\infty}
\frac{1}{\ell} \Re\biggl[\bigl(\bar{\varpi}_{N}(x)^{\ell} -
\bar{\varpi}_{N}(w)^{\ell} \bigr)\cdot\4 \frac{1}{\bar{\varpi}_{N}(y_i)^{\ell} }
\biggr] .
\end{eqnarray*}
Then, since $ |\Re[ab]|\le|a||b|$
and
$ |\bar{a}^{\ell}-\bar{b}^{\ell}|=|a^{\ell} - b^{\ell}|$,
we have
%
\begin{eqnarray}
\label{63a} &&\frac{|\PsiN_{rs}(x,\mathsf{y}) - \PsiN_{rs}(w ,\mathsf
{y})|}{|x-w|}
\nonumber
\\[-8pt]
\\[-8pt]
\nonumber
&&\qquad \le |\beta| \sum_{\ell= 1}^{\infty}
\frac{|\wNxl- \wNwl|}{\ell|x-w|} \cdot\biggl| \4 \frac{1}{\bar{\varpi
}_{N}(y_i)^{\ell} }\biggr| .
\end{eqnarray}

Our purpose is to estimate
$ {|\PsiN_{rs}(x,\mathsf{y}) - \PsiN_{rs}(w ,\mathsf{y})|}/{|x-w|}$
for $ x\not=w\in\Sr$. Hence, by (\ref{63a}), the main task is
to control the term of the form
\[
\biggl|\4 \frac{1}{\bar{\varpi}_{N}(y_i)^{\ell} }\biggr| .
\]
Taking this into account, we set for $ r , \ell, {k}\in\N$ ,
%
%
\begin{eqnarray}
\label{A7c} \mathsf{U}_{r,\ell,k}& =& \biggl\{ \mathsf{y}\in\mathsf{S} ;
\sup_{N\in\N} \sup_{r< s\in\N}\biggl|\9 \biggr| \le{k}\biggr\} ,
\\
\label{A7d} \bar{\mathsf{U}}_{r,\ell,k} &=& \biggl\{ \mathsf{y}\in
\mathsf{S} ; \sup_{N\in\N} \biggl\{ \sum_{\yyi\in\Sri}
\frac{1}{|\wNyi|^{\ell}- |\wN(b_r )|^{\ell}} \biggr\} \le{k}\biggr\} .
\end{eqnarray}
%
%
\begin{rem}\label{r62} When (\ref{61t}) holds, definitions
(\ref{A7c}) and (\ref{A7d}) become much simpler.
%
\begin{eqnarray}
\label{63s} \mathsf{U}_{r,\ell,k}& =& \biggl\{ \mathsf{y}\in\mathsf{S} ;
\sup_{r< s\in\N}\biggl|\sum_{\yyi\in\Srs} \frac{1}{y_i^{\ell}} \biggr|
\le{k}\biggr\} ,
\\
\label{63t}  \bar{\mathsf{U}}_{r,\ell,k}& =& \biggl\{ \mathsf{y}\in
\mathsf{S} ; \biggl\{ \sum_{\yyi\in\SSSS_{r\infty}} \frac{1}{|\yyi|^{\ell}- b_r^{\ell}}
\biggr\} \le{k}\biggr\} .
\end{eqnarray}
\end{rem}


(A.6)\hypertarget{A6}\ For each $ r \in\N$, there exists an $ {\ell_0}\in\N
$ such that
%
\begin{eqnarray}
\label{A7b}  \limi{k} \limsup_{ N \to\infty} \muN\bigl(\bar{
\mathsf{U}}_{r,{\ell_0},k}^c \bigr) &=& 0 ,
\\
\label{A7a}  \limi{k} \limsup_{ N \to\infty} \muN\bigl(\mathsf{U}_{r,\ell,k}^c
\bigr) &=& 0 \qquad\mbox{for all $ 1 \le\ell< {\ell_0}$} .
\end{eqnarray}
When $ {\ell_0}= 1 $, according to our interpretation,
(\ref{A7a}) always holds by convention.

We now state the main theorem of this section.
%
\begin{thmm} \label{l63}
Assume (\ref{61a}) and (\ref{61s}).
Suppose (\ref{61t}) or (\ref{61u}).
Then \textup{(\hyperlink{A6}{A.6})} implies \textup{(\hyperlink{A5}{A.5})}.
\end{thmm}

\begin{pf}
Let $ \Ct\scriptsize\label{;6x} $ and $ \Ct\label{;6y} $ be the constants
defined by
%
\begin{eqnarray}
\label{63d} c_{\scriptsize\ref{;6x}} &=& |\beta| \cdot\sup_{N\in\N}
\max_{1 \le\ell<
{\ell_0}} \sup_{\7} \frac{|\wNxl- \wNwl|}{\ell|x-w|} ,
\nonumber
\\[-8pt]
\\[-8pt]
\nonumber
 c_{\scriptsize\ref{;6y}} &=& |\beta| \cdot\sup_{N\in\N}
\sup_{{\ell_0}\le
\ell
} \sup_{\7} \frac{|\wNxl- \wNwl|}{ |\wN(b_r )|^{\ell} \ell|x-w|} .
\end{eqnarray}
Then $ c_{\scriptsize\ref{;6x}}$ and $ c_{\scriptsize\ref{;6y}}$ are finite.
Indeed, $ c_{\scriptsize\ref{;6x}} < \infty$ is clear.
Note that the Lipschitz norm of $ \{ \wN\} $ on $ \R$ is uniformly
bounded in $ N \in\N$.
Moreover, $ |\wNx|/ |\wN(b_r )| < 1$ on $ \Sr$.
Hence the Lipschitz norm of the function $ {\wNxl}/ { \wN(b_r )^{\ell
} \ell}$ on $ \Sr$
is uniformly bounded in $ \ell, N\in\N$. This implies $ c_{\scriptsize\ref{;6y}} < \infty$.

By (\ref{63a}) and $ c_{\scriptsize\ref{;6x}}, c_{\scriptsize\ref{;6y}} < \infty$, we have
%
\begin{eqnarray}
\label{63e} &&\frac{|\PsiN_{rs}(x,\mathsf{y}) -\PsiN_{rs}(w,\mathsf{y}) |}{|x-w|}
\nonumber\\
&&\qquad\le c_{\scriptsize\ref{;6x}} \sum_{\ell= 1}^{{\ell_0}-1}
\biggl|\9 \biggr| + c_{\scriptsize\ref{;6y}} \sum_{\ell= {\ell_0}}^{\infty}
\sum_{\yyi\in\Srs}\frac{ |\wN(b_r )|^{\ell}}{|\wNyi|^{\ell}}
\nonumber
\\[-8pt]
\\[-8pt]
\nonumber
&&\qquad=  c_{\scriptsize\ref{;6x}} \sum_{\ell= 1}^{{\ell_0} -1}
\biggl|\9 \biggr| \\
&&\qquad\quad{}+ c_{\scriptsize\ref{;6y}} \sum_{\yyi\in\Srs}
\frac{ |\wN(b_r )|^{{\ell_0}}}{|\wNyi|^{{\ell_0}} - |\wN(b_r
)|^{{\ell_0}} } .\nonumber
\end{eqnarray}
Here we used the formula
$ \sum_{\ell= {\ell_0}}^{\infty} { a^{\ell} }/{b^{\ell}} =
{a^{{\ell_0}}}/({b^{{\ell_0}} - a^{{\ell_0}} })$ valid for $ 0< a\le
b $.
If $ a=b$, then we interpret $ \sum_{\ell= {\ell_0}}^{\infty} {
a^{\ell} }/{b^{\ell}} = \infty$.
Set $ \Ct\label{;6z} = c_{\scriptsize\ref{;6y}}\sup_{N\in\N} |\wN(b_r
)|^{{\ell_0}}$.
By~(\ref{63e}), we see that
\begin{eqnarray*}
&&
\2 \frac{|\PsiN_{rs}(x,\mathsf{y}) -\PsiN_{rs}(w,\mathsf{y})|} {|x-w|}
\\
& &\qquad\le c_{\scriptsize\ref{;6x}} \sum_{\ell= 1}^{{\ell_0}-1}
\biggl\{ \sup_{N\in\N} \sup_{r< s\in\N}\biggl|\9 \biggr| \biggr\}
\\
&&\qquad\quad{} + c_{\scriptsize\ref{;6z}}\biggl\{ \sup_{N\in\N} \sum
_{\yyi\in\Sri} \frac{1}{|\wNyi|^{{\ell_0}} - |\wN(b_r )|^{{\ell_0}}}\biggr\} .
\end{eqnarray*}
Combining this with (\ref{41f}), (\ref{A7c}) and (\ref{A7d}), we
deduce that
\[
\Hrk \supset \Biggl\{ \bigcap_{\ell= 1}^{{\ell_0}-1}
\mathsf{U}_{r,\ell,k/({\ell_0}c_{\scriptsize\ref{;6x}})} \Biggr\} \cap \bar{\mathsf{U}}_{r,{\ell_0},k/({\ell_0}c_{\scriptsize\ref{;6z}})}
.
\]
Hence, we obtain
%
\begin{equation}
\label{63g} \muN\bigl(\Hrk^c\bigr) \le \Biggl\{\sum
_{\ell= 1}^{{\ell_0}-1} \muN\bigl(\mathsf{U}_{r,\ell,k/({\ell_0}c_{\scriptsize\ref{;6x}})}^c
\bigr) \Biggr\} + \muN\bigl(\bar{\mathsf{U}}_{r,{\ell_0},k/({\ell_0}c_{\scriptsize\ref{;6z}})}^c\bigr)
.
\end{equation}
This together with (\hyperlink{A6}{A.6}) implies (\ref{41g}), which completes
the proof.
\end{pf}

\section{\texorpdfstring{Sufficient conditions of (\protect\hyperlink{A6}{A.6})}{Sufficient conditions of (A.6)}} \label{s7}
In this section, we give sufficient conditions of (\hyperlink{A6}{A.6}).
These conditions are used in the proof of Theorem~\ref{l22} and
Theorem~\ref{l23}.
We begin with (\ref{A7b}), the first condition of (\hyperlink{A6}{A.6}).
%
%
\begin{lem} \label{l71}
Assume \textup{(\hyperlink{A4}{A.4})}, (\ref{61a}) and (\ref{61s}).
Assume (\ref{61t}) or (\ref{61u}).
Then (\ref{A7b}) follows from (\ref{71a}) below.
%
\begin{equation}
\label{71a} \sup_{N\in\N} \biggl\{ \int_{1\le|x |<\infty} \biggl
\{ \sup_{M\in\N} \frac{1 }{|\wM(x) |^{{\ell_0}}} \biggr\} \rNone(x )\,dx \biggr\} <
\infty .
\end{equation}
In particular, if (\ref{61t}) is satisfied, then (\ref{A7b})
follows from a simpler condition~(\ref{71z}),
%
\begin{equation}
\label{71z} \sup_{N\in\N} \biggl\{ \int_{1\le|x |<\infty}
\frac{1}{|x|^{{\ell_0}}} \rNone(x )\,dx \biggr\} < \infty .
\end{equation}
\end{lem}
\begin{pf}
Let $ b_r $ be as in (\ref{qg0}). We divide the set
$ \Sri= \{ b_r \le|x | < \infty\} $ in~(\ref{A7d}) into two parts,
$ \SSSS_{r(r+1)} = \{ b_r \le|x | < b_{r+1}\} $ and
$ \SSSS_{(r+1)\infty}= \{ b_{r+1}\le|x | < \infty\} $.
Let $ \xx=\sum_{i}\delta_{x _i}$. We set
\begin{eqnarray*}
 \mathsf{V}_{1,k}& =& \biggl\{ \xx\in\mathsf{S} ; \biggl\{
\sup_{N\in\N} \sum_{\xxi\in\SSSS_{r(r+1)} } \frac{1}{|\wNxi|^{{\ell_0}}- |\wN(b_r )|^{{\ell_0}}}
\biggr\} \le\frac{k}{2}\biggr\}, 
\\
 \mathsf{V}_{2,k} &=& \biggl\{ \xx\in\mathsf{S} ; \biggl\{
\sup_{N\in\N} \sum_{\xxi\in\SSSS_{(r+1)\infty} } \frac
{1}{|\wNxi|^{{\ell_0}}- |\wN(b_r )|^{{\ell_0}}}
\biggr\} \le\frac{k}{2}\biggr\} .
\end{eqnarray*}
Then clearly
$ \bar{\mathsf{U}}_{r,{\ell_0},k} \supset\mathsf{V}_{1,k}\cap
\mathsf{V}_{2,k}$.
To estimate $ \mathsf{V}_{1,k} $, we observe that
\begin{eqnarray*}
&&
\sum_{\xxi\in\SSSS_{r(r+1)} } \frac{1}{|\wNxi|^{{\ell_0}}- |\wN(b_r )|^{{\ell_0}}}
\\
&&\qquad\le \biggl\{ \sup_{\xxi\in\SSSS_{r(r+1)} } \frac{1}{|\wNxi|^{{\ell_0}}- |\wN(b_r )|^{{\ell_0}}} \biggr\}\cdot\xx
(\SSSS_{r(r+1)} ) .
\end{eqnarray*}
Here $ \xx(\SSSS_{r(r+1)} )$ is the number of points $ \xxi$ in $
\SSSS_{r(r+1)} $.
Taking this into account, we set
\begin{eqnarray*}
\mathsf{V}_{3,k} &=& \biggl\{ \xx\in\mathsf{S} ; \sup_{N\in\N}
\sup_{\xxi\in\SSSS_{r(r+1)} } \frac{1}{|\wNxi|^{{\ell_0}}- |\wN(b_r )|^{{\ell_0}}} \le\sqrt{{{k}}/{2}} \biggr\},
\\
\mathsf{V}_{4,k} &=& \bigl\{ \xx\in\mathsf{S} ; \xx(
\SSSS_{r(r+1)} ) \le\sqrt{{{k}}/{2}} \bigr\}.
\end{eqnarray*}
Then we have $ \mathsf{V}_{1,k} \supset\mathsf{V}_{3,k} \cap
\mathsf{V}_{4,k} $.
We therefore obtain
$\bar{\mathsf{U}}_{r,{\ell_0},k} \supset\mathsf{V}_{2,k}\cap
\mathsf{V}_{3,k}\cap\mathsf{V}_{4,k}$ by combining these two
inclusions. Hence we deduce (\ref{A7b}) from
%
\begin{equation}
\label{71g} \limi{k} \limsup_{ N \to\infty} \muN\bigl(\mathsf{V}_{l,k}^c
\bigr) = 0\qquad \mbox{for all }l=2,3,4 .
\end{equation}
We will check (\ref{71g}) for each $ l =2,3,4 $.

As for (\ref{71g}) with $ l =2 $, according to the Chebyshev
inequality, we have
%
\begin{eqnarray}\label{71h}
&& \muN\bigl(\mathsf{V}_{2,k}^c\bigr)
\nonumber\\
&&\qquad\le \frac{2}{k} E^{\muN}\biggl[ \sup_{M\in\N}
\sum_{\xxi\in\SSSS_{(r+1)\infty} } \frac{1}{|\wMxi|^{{\ell_0}}- |\wM(b_r )|^{{\ell_0}}} \biggr]
\nonumber
\\
&&\qquad= \frac{2}{k} \int_{ \SSSS_{(r+1)\infty}}
\sup_{M\in\N} \biggl\{ \frac{1}{ |\wM(x) |^{{\ell_0}}-
|\wM(b_r )|^{{\ell_0}}} \biggr\} \rNone \,dx
\nonumber
\\[-8pt]
\\[-8pt]
\nonumber
&&\qquad=  \frac{2}{k} \int_{ \SSSS_{(r+1)\infty}}
\sup_{M\in\N} \biggl\{ \frac{{|\wM(x) |^{{\ell_0}}}} {
|\wM(x) |^{{\ell_0}} - |\wM(b_r )|^{{\ell_0}}} \frac{1}{{|\wM(x) |^{{\ell_0}}}} \biggr\}
\rNone \,dx
\\
\nonumber
&&\qquad\le\frac{2}{k} \sup_{M\in\N} \biggl\{
\frac{|\wM(b_{r+1})|^{{\ell_0}}} {
|\wM(b_{r+1})|^{{\ell_0}}- |\wM(b_r )|^{{\ell_0}}} \biggr\} 
\\
&&\quad\qquad{}\times \int_{\SSSS_{(r+1)\infty}}
\sup_{M\in\N} \biggl\{ \frac{1 }{|\wM(x) |^{{\ell_0}}} \biggr\} \rNone \,dx
.\nonumber
\end{eqnarray}
Here we used the fact that $ |\wM(x) | < |\wM(y) |$
for $ |x| < |y| $, which implies
\begin{eqnarray*}
\sup_{x \in\SSSS_{(r+1)\infty}} \frac
{|\wM(x) |^{{\ell_0}}} {
|\wM(x) |^{{\ell_0}}- |\wM(b_r )|^{{\ell_0}}} \le\frac
{|\wM(b_{r+1})|^{{\ell_0}}} {
|\wM(b_{r+1})|^{{\ell_0}}- |\wM(b_r )|^{{\ell_0}}}.
\end{eqnarray*}
By (\ref{71a}) and (\ref{71h}), we obtain
(\ref{71g}) with $ l=2 $.

We next consider (\ref{71g}) with $ l =3 $. Let
\[
U_{k} = \bigcup_{\NN\in\N} \bigl\{ x \in
\SSSS_{r(r+1)} ; \bigl|\wN(b_r )\bigr|^{{\ell_0}} \le|x
|^{{\ell
_0}} < \bigl|\wN(b_r )\bigr|^{{\ell_0}} + {\sqrt{2/{k}}}\bigr
\} .
\]
It is not difficult to see that $ U_{k} $ is nonincreasing, and
$ \limi{k} U_{k} = \varnothing$. We note that
%
\begin{eqnarray}
\label{71j}
 \mathsf{V}_{3,k}^c &=& \Bigl\{ \xx\in\mathsf{S}
; \inf_{N\in\N}\inf_{\xxi\in\SSSS_{r(r+1)} }\bigl\{ {\bigl|\wNxi\bigr|^{{\ell_0}}- \bigl|
\wN(b_r )\bigr|^{{\ell_0}}} \bigr\} < \sqrt{2/{k}} \Bigr\}
\nonumber
\\[-8pt]
\\[-8pt]
\nonumber
& =& \bigl\{ \xx\in\mathsf{S} ; 1 \le\xx( U_{k} ) \bigr\} .
\end{eqnarray}
Here we use a convention such that $ \inf\varnothing= \infty$;
that is, we interpret $ \xx\notin\mathsf{V}_{3,k}^c $ when $ \xx
(\SSSS_{r(r+1)} )=0$.
Let $ c_{\scriptsize\ref{;71k}}= \sup\{ \rNone(x ); N\in\N, x \in\SSSS_{r(r+1)} \} $.
Then by (\ref{41b}), we have $ \Ct\label{;71k} < \infty$.
From the second equality in (\ref{71j}) and the Chebyshev inequality,
we obtain
%
\begin{equation}
\label{71k} \muN\bigl(\mathsf{V}_{3,k}^c\bigr) \le
E^{\muN}\bigl[\xx( U_{k} )\bigr] = \int_{ U_{k} }
\rNone(x )\,dx \le c_{\scriptsize\ref{;71k}} \int_{ U_{k} }\,dx .
\end{equation}
Hence, we deduce (\ref{71g}) with $ l =3 $ from (\ref{71k}) and $
\limi{k} U_{k} = \varnothing$.

We finally consider (\ref{71g}) with $ l =4 $. From the Chebyshev
inequality we obtain
\begin{eqnarray*}
\muN\bigl(\mathsf{V}_{4,k}^c\bigr)
\le\sqrt{\frac{2}{{k}}} E^{\muN}\bigl[\xx (\SSSS_{r(r+1)} )
\bigr] = \sqrt{\frac{2}{{k}}} \int_{\SSSS_{r(r+1)} } \rNone(x )\,dx
\le\sqrt{\frac{2}{{k}}} c_{\scriptsize\ref{;71k}} \int_{\SSSS_{r(r+1)} } \,dx.
\end{eqnarray*}
This deduces (\ref{71g}) with $ l =4 $ immediately.
\end{pf}

We proceed with (\ref{A7a}), the second condition of (\hyperlink{A6}{A.6}).

Let $ \Dr= \{ s \in\SSSS; |s|< r \}$ and $ \Drs= \Ds\setminus
\Dr$.
Let $ \map{\mathsf{v}_{\ell,rs}^N }{\SSS}{\mathbb{C}}$ such that
%
\begin{eqnarray}
\label{72a}\mathsf{v}_{\ell,rs}^N (\xx)= \sum
_{\xxi\in\Drs} \frac{1}{\bar
{\varpi}_{N}(x_i)^{\ell}}\qquad \mbox{for }1 \le r < s \le\infty .
\end{eqnarray}
Here we write $ \xx=\sum_{i}\delta_{\xxi}$, as usual.
Note that the sum in (\ref{72a}) makes sense for $ \muN$-a.s. $ \xx$
even if $ s=\infty$.
Indeed, by (2) of (\hyperlink{A4}{A.4}), the total number of particles
has the deterministic bound $ n_{N}$ under $ \muN$.
Hence, $ \mathsf{v}_{\ell,rs}^N (\xx)$ is well defined and finite
for $ \muN$-a.s. $ \xx$ for all $ N\in\N$.

%
\begin{lem} \label{l72}
Under the same assumptions as Lemma~\ref{l71}, (\ref{A7a}) follows from~(\ref{72e}) below.
%
\begin{eqnarray}
\label{72e} \limi{r} \sup_{N\in\N} \Bigl\| \sup_{M\in\N} \bigl|
\mathsf{v}_{\ell
,r\infty}^M\bigr| \Bigr\|_{\LmNone} = 0 \qquad\mbox{for all $
1 \le\ell< {\ell_0}$} .
\end{eqnarray}
\end{lem}

\begin{pf}
By (\ref{72e}), we can and do choose $ \{b_r \}$ and
$ \Ct\label{;61} > 0 $ in such a way that
%
\begin{equation}
\label{72c}\sup_{N\in\N} \Bigl\| \sup_{M\in\N} \bigl|\mathsf{v}_{\ell,b_r \infty
}^M\bigr|
\Bigr\|_{\LmNone} \le c_{\scriptsize\ref{;61}}3^{- r } \qquad\mbox{for all } r \in
\N .
\end{equation}
We note that
$ \sum_{\yyi\in\Srs} {1}/{\bar{\varpi}_{N}(y_i)^{\ell}} =
\mathsf{v}_{\ell,b_r \infty}^N (\xx)- \mathsf{v}_{\ell,b_s \infty
}^N (\xx)$.
Then by (\ref{A7c}), we see that
%
%
\begin{eqnarray}
\label{72d} \muN \bigl(\{ \mathsf{U}_{r,\ell,k} \}^c\bigr)& =&
\muN\Bigl( \sup_{M\in\N} \sup_{r< s\in\N}\bigl| \mathsf{v}_{\ell,b_r
\infty}^M-
\mathsf{v}_{\ell,b_s \infty}^M\bigr|> k \Bigr)\nonumber
\\
&\le& \muN\Bigl( \sup_{M\in\N} \bigl| \mathsf{v}_{\ell,b_r \infty}^M\bigr|>
{k}/{2}\Bigr) + \muN\Bigl( \sup_{M\in\N} \sup_{r< s\in\N}\bigl|
\mathsf{v}_{\ell,b_s
\infty}^M\bigr|> {k}/{2}\Bigr)
\nonumber
\\[-8pt]
\\[-8pt]
\nonumber
&\le& \muN\Bigl( \sup_{M\in\N} \bigl| \mathsf{v}_{\ell,b_r \infty}^M\bigr|>
{k}/{2}\Bigr) + \sum_{s = r+1}^{\infty} \muN\Bigl(
\sup_{M\in\N} \bigl|\mathsf{v}_{\ell
,b_s \infty}^M\bigr|> {k}/{2}\Bigr)
\\
\nonumber
&\le& \frac{2}{k}\cdot\Biggl\{ \sum_{s = r}^{\infty}
\Bigl\|\sup_{M\in\N} \bigl|\mathsf{v}_{\ell,b_s \infty}^M\bigr|
\Bigr\|_{\LmNone}\Biggr\} .
\end{eqnarray}
Here we used Chebyshev's inequality in the last line.
By (\ref{72c}) and (\ref{72d}) we have
\[
\sup_{N\in\N} \muN\bigl(\{\mathsf{U}_{r,\ell,k}\}^c
\bigr) \le \frac{2}{k}\cdot\frac{c_{\scriptsize\ref{;61}} 3^{-r}}{1-3^{-1}} .
\]
Hence, $ \limi{k}\sup_{N\in\N} \muN(\{\mathsf{U}_{r,\ell,k}\}^c) =0$,
which implies (\ref{A7a}).
\end{pf}

We refine Lemma~\ref{l72} in Lemma~\ref{l81}, used in the proof of Theorems
\ref{l22} and \ref{l23} directly.
%
%
\begin{lem} \label{l81}
Let $ \map{ u_{\ell,r}^{N} }{\SSSS}{\mathbb{C}}$ such that
\[
u_{\ell,r}^{N} (x) = 1_{\Sti_{1r} }(x) {\bigl\lceil\bigl| \wNx\bigr|
\bigr\rceil^{\ell} }/{\bar{\varpi }_{N}(x)^{\ell} } .
\]
Here $ \lceil\cdot\rceil$ is the minimal integer greater than or
equal to $ \cdot$.
Let $ \map{\mathsf{u}_{\ell,r}^{N}}{\SSS}{\mathbb{C}}$ such that
$ \mathsf{u}_{\ell,r}^{N}(\xx) = \sum_i u_{\ell,r}^{N} (x_i)$,
where\vadjust{\goodbreak} $ \mathsf{x}=\sum_i \delta_{x_i}$.
Suppose there exists a positive constant $ \Ct\label{;r71}$ such that
%
\begin{equation}
\label{81b} \sup_{r\in\N} r^{c_{\scriptsize\ref{;r71}}-\ell} \sup_{N\in\N} \Bigl\|
\sup_{M\in
\N} \bigl|\mathsf{u}_{\ell,r}^{M}\bigr|
\Bigr\|_{\LmNone} <\infty\qquad \mbox{for all $ 1 \le\ell< {\ell_0}$} .
\end{equation}
In addition, assume the same conditions as for Lemma~\ref{l71}. Then
(\ref{A7a}) holds.
\end{lem}
\begin{pf}
Let $ 1 \le\ell< {\ell_0}$ be fixed. Define $ \map{\mathsf
{w}^{j}_{r } }{\SSS}{\mathbb{C}}$ by
\[
\mathsf{w}^{j}_{r } (\xx) = \sum
_{\xxi\in\Sti_{1r} } \frac{\lceil| \wM(x_i) |\rceil^j }{\bar{\varpi}_M (x_i)^{\ell}} .
\]
Although $ \mathsf{w}^{j}_{r } $ depends on $ M \in\N$,
we omit $ M $ from the notation for simplicity.
Let $ \mathsf{v}_{\ell,1r}^M $ be as in (\ref{72a}).
Then $ \mathsf{w}^{0}_{r } = \mathsf{v}_{\ell,1r}^M $ and
$ \mathsf{w}^{\ell}_{r } = \mathsf{u}_{\ell,r}^{M}$ by definition.
Moreover, we easily deduce that
\[
\mathsf{w}^{j}_{r } = \sum
_{q = 2}^r q \bigl( \mathsf{w}^{j-1}_{q }
- \mathsf{w}^{j-1}_{q-1 }\bigr) \qquad\mbox{for } r\ge2 ,\qquad
\mathsf{w}^{j}_{1} = 0 .
\]
Hence, through a straightforward calculation, we have
%
\begin{eqnarray}
\label{81d} \urr= \frac{\ur}{r} + \sum_{q = 2}^{r-1}
\frac{\up}{q(q+1)}\qquad \mbox{for } r\ge3 , \qquad\mathsf{w}^{j-1}_{2}
=\frac
{1}{2}\mathsf{w}^{j}_{2} .
\end{eqnarray}

By (2) of (\hyperlink{A4}{A.4}), we see that
$ \limi{r}r^{-1}\| \mathsf{w}^{j}_{r } \|_{\LmNone} =0 $ and that
$ \mathsf{w}^{j}_{\infty}:= \limi{r}\mathsf{w}^{j}_{r } $ exists
in $ \LmNone$. Hence, by taking $ r \to\infty$ in (\ref{81d}), we obtain
\[
\mathsf{w}^{j-1}_{\infty} = \sum
_{q = 2}^{\infty} \frac{\up}{q(q+1)} \qquad\mbox{in $L^1\bigl(\SSS, \muN\bigr)$} .
\]
Subtracting (\ref{81d}) from this yields
%
\begin{equation}
\label{81h} \mathsf{w}^{j-1}_{\infty} - \mathsf{w}^{j-1}_{r }
= - \frac{\ur}{r} + \sum_{q = r}^{\infty}
\frac{\up}{q(q+1)} \qquad\mbox{in $L^1\bigl(\SSS, \muN\bigr)$} .
\end{equation}

Take the supremum of the modulus of each terms of
(\ref{81d}) and (\ref{81h}) w.r.t. $ M \in\N$.
Apply Minkowski's inequality to the right-hand sides of
(\ref{81d}) and (\ref{81h}).
Then by taking the supremum w.r.t. $ N\in\N$, we obtain
%
\begin{eqnarray}
\label{81e}\qquad &&\sup_{N\in\N} \Bigl\| \sup_{M\in\N} \bigl|\urr\bigr|
\Bigr\|_{\LmNone}
\nonumber
\\[-4pt]
\\[-12pt]
\nonumber
& &\qquad\le\frac{1}{r} \sup_{N\in\N} \Bigl\| \sup_{M\in\N} \bigl|\ur\bigr|
\Bigr\|_{\LmNone} + \sum_{q = 2}^{r-1}
\frac{1}{q(q+1)}\sup_{N\in\N} \Bigl\| \sup_{M\in\N} \bigl|\up\bigr|
\Bigr\|_{\LmNone} ,
\\
\label{81p}
&& \sup_{N\in\N} \Bigl\|\sup_{M\in\N} \bigl|\mathsf{w}^{j-1}_{\infty}
- \mathsf{w}^{j-1}_{r }\bigr|\Bigr \|_{\LmNone}
\nonumber
\\[-4pt]
\\[-12pt]
\nonumber
&&\qquad \le\frac{1}{r}\sup_{N\in\N} \Bigl\|\sup_{M\in\N} \bigl|\ur\bigr|
\Bigr\|_{\LmNone
} + \sum_{q = r}^{\infty}
\frac{1}{q(q+1)} \sup_{N\in\N} \Bigl\|\sup_{M\in\N} \bigl|\up\bigr|
\Bigr\|_{\LmNone} .
\end{eqnarray}

For each $ j = 1,\ldots,\ell$, there exists
a positive constant $ \Ct\label{;81} = c_{\scriptsize\ref{;81}}( j ) $ such that
%
%
\begin{equation}
\label{81f} \sup_{r\in\N} r^{c_{\scriptsize\ref{;81}}- j } \sup_{N\in\N} \Bigl\|
\sup_{M\in\N} \bigl|\ur\bigr| \Bigr\|_{\LmNone} <\infty
.
\end{equation}
Indeed, when $ j = \ell$, (\ref{81f}) holds by (\ref{81b}) because
$ \mathsf{w}^{\ell}_{r } = \mathsf{u}_{\ell,r}^{M}$.
Suppose (\ref{81f}) holds for some $ 2\le j \le\ell$
with a positive constant $ c_{\scriptsize\ref{;81}}( j )$.
Then by (\ref{81e}), we have~(\ref{81f}) for $ j -1 $ with
a positive constant $ c_{\scriptsize\ref{;81}}( j -1) $. Therefore, through induction,
(\ref{81f}) holds for all $ j = 1,\ldots,\ell$.

Combining (\ref{81p}) and (\ref{81f}), we easily deduce that
for each $ j = 1,\ldots,\ell$,
%
\begin{equation}
\label{81q} \limi{r} \sup_{N\in\N} \Bigl\|\sup_{M\in\N} \bigl|\mathsf
{w}^{j-1}_{\infty} - \mathsf{w}^{j-1}_{r }\bigr|
\Bigr\|_{\LmNone} = 0 .
\end{equation}
Recalling $ \mathsf{w}^{0}_{r } = \mathsf{v}_{\ell,1r}^M $, we have
$ \mathsf{w}^{0}_{\infty} - \mathsf{w}^{0}_{r }= \mathsf{v}_{\ell
,r\infty}^M$.
Hence, by taking $ j=1$ in (\ref{81q}), we obtain
\[
 \limi{r} \sup_{N\in\N} \Bigl\|\sup_{M\in\N} \bigl| \mathsf{v}_{\ell
,r\infty}^M
\bigr| \Bigr\|_{\LmNone} = 0 .
\]
This allows us to deduce (\ref{72e}) in Lemma~\ref{l72}. We therefore
obtain (\ref{A7a}) by Lemma~\ref{l72}.
\end{pf}

\section{Translation invariant periodic measures}\label{s8}
In this section, we make preparations for a proof of Theorem~\ref{l22}.

Let $\SSSS= \Rd$. Let $ \map{\tau_x}{\SSS}{\SSS} $ be the
translation defined by
$ \tau_x (\mathsf{s}) = \sum_i \delta_{x + s_i}$
for $ \mathsf{s} = \sum_i \delta_{s_i} $.
We say that a measure $ \nu$ on $ \SSS$
is translation invariant if
$ \nu\circ\tau_x^{-1} = \nu$ for all $ x \in\Rd$.
We say that $ \nu$ is $ L $-periodic if
$ \nu(\tau_{L\mathbf{e}_i}(\mathsf{s}) = \mathsf{s}) = 1 $
for all $ i = 1,\ldots,d $.
Moreover, we say that $ \nu$ is concentrated on $ A $ if
$ \nu(\mathsf{s}(A^c)>0) = 0 $.
A measure $ \nu$ concentrated on $ (-L/2,L/2]^d $ can be extended
naturally to
the $ L $-periodic measure $ \bar{\nu}$ on the configuration space on
$ \Rd$.
We refer to this measure $ \bar{\nu} $ as the $ L $-periodic
extension of $ \nu$.

Let $ \TN= (-n_{N}/2,n_{N}/2]^d $. We assume that $ \nu$ is
concentrated on $ \TN$ and that~$ \nu$ has a periodic extension that
is translation invariant.
Let $ \rN$ be the $ n $-correlation function of $ \nu$.
Then $ \rN(x) = 0 $ for $ x\notin(\TN)^{n_{N}}$ by assumption.
Let $ \T$ be the two-level cluster function of $ \nu$,
%
\begin{equation}
\label{60T} \T(x,y) = \rNone(x) \rNone(y) - \rNtwo(x,y) .
\end{equation}
Then $ \T(x,y) = 0 $ if $ (x,y)\notin(\TN)^2 $. If $ (x,y)\in(\TN
)^2 $,
$ \T(x,y)$ depends only on $ x-y $ modulo
$ N \mathbf{e}_i $ $ (i = 1,\ldots,d) $, where
$ \mathbf{e}_i $ is the $ i $th unit vector.
Therefore, let $ \map{\T}{\Rd}{\R} $ be the $ n_{N}$-periodic function
such that $ \T(x) = \T(x,0)$ for $ x\in\TN$.
We set
%
%
\begin{equation}
\label{63c} \mN(\xi) = \rNone(0) - \mathcal{F}_{N}(\T) (\xi) .
\end{equation}
Here $ \mathcal{F}_{N}(f)(\xi) =
\int_{\Rd} e^{-2\pi\sqrt{-1} \xi\cdot x }
f 1_{\TN}(x )\,dx $
denotes the Fourier transform of $ f 1_{\TN}$.
%
%
\begin{lem} \label{l61a} 
Assume that $ \nu$ is concentrated on
$ \TN$ and that $ \nu$ has a periodic extension that is translation invariant.
Let $ \map{h }{\Rd}{\R}$ be real valued. Set
$\mathsf{h}_N(\mathsf{s}) = \sum_{s_i\in\TN} h (s_i)$,
where $\mathsf{s} = \sum_i \delta_{s_{i}}$.
Then
\begin{eqnarray*}
 \|\mathsf{h}_N \|_{\Ln}^2 = \biggl
\{\rNone(0)\int_{\TN}h(x)\,dx \biggr\}^2 + \3 \bigl|
\mathcal{F}_{N}(h)\bigr|^2(\xi)\mN(\xi) .
\end{eqnarray*}
\end{lem}
\begin{pf}
From $ \rNone(x) = \rNone(0)1_{\TN}(x)$, we see that
%
\begin{equation}
\label{63f} \int_{\SSS}\mathsf{h}_N \,d\nu= \int
_{\TN} h (x) \rNone(x) \,dx = \rNone(0)\int
_{\TN} h (x) \,dx .
\end{equation}

Let $ \operatorname{Var}^{\nu}[\mathsf{h}_N ]$
be the variance of $ \mathsf{h}_N $
w.r.t. $ \nu$. %
By (\ref{60T}) and the general property of correlation functions,
we see that
\begin{eqnarray*}
\operatorname{Var}^{\nu}[\mathsf{h}_N ] & = &\int
_{\Rd} h^2 (x) \rNone(x)\,dx - \int
_{\Rd\ts\Rd} h(x) h (y) \T(x,y) \,dx \,dy
\\
& =& \rNone(0) \int_{\TN} h^2 (x) \,dx -
\int_{\TN\ts\TN} h(x) h (y) \T(x-y) \,dx \,dy .
\end{eqnarray*}
We used
$ \rNone(x) = \rNone(0)1_{\TN}(x)$ and
$ \T(x,y) = 1_{\TN}(x)1_{\TN}(y)\T(x-y)$ in the second line.
By a direct calculation of the Fourier series, we see that
\[
\int_{\TN} h^2 (x) \,dx = \3 \bigl|
\mathcal{F}_{N}(h) (\xi) \bigr|^2
\]
and
\begin{eqnarray*}
&& \int_{\TN\ts\TN} h(x) h (y)\T(x-y) \,dx \,dy
\\
&&\qquad = \3 \mathcal{F}_{N}(h) (\xi) \overline{
\mathcal{F}_{N}(h* \T) (\xi)}
\\
&&\qquad = \3 \bigl|\mathcal{F}_{N}(h) (\xi) \bigr|^2 \overline{
\mathcal{F}_{N}( \T) (\xi) }
\\
& &\qquad= \3 \bigl|\mathcal{F}_{N}(h) (\xi) \bigr|^2
\mathcal{F}_{N}(\T) (\xi) .
\end{eqnarray*}
Here we used the fact that $ \mathcal{F}_{N}(\T) $ is real valued
because $ \T(x) = \T(-x)$.
Combining these with (\ref{63c}) yields
%
\begin{equation}
\label{63dd} \operatorname{Var}^{\nu}[\mathsf{h}_N ] = \3\bigl |
\mathcal{F}_{N}(h)\bigr|^2(\xi)\mN(\xi) .
\end{equation}

We conclude Lemma~\ref{l61a} from (\ref{63f}) and (\ref{63dd})
immediately.
\end{pf}


\section{\texorpdfstring{Proof of Theorems \protect\ref{l22}}{Proof of Theorems 2.2}} \label{s9}
In this section, we prove Theorem~\ref{l22} using the previous results.
We begin by defining $ \Ksinb$ for $ \beta= 1,4 $.
Let $ \ii= \sqrt{-1}$, as before.
To define $ \Ksinb$, we recall the standard quaternion notation
for $ 2\ts2 $ matrices (see~\cite{mehta}, Chapter~2.4),
\begin{eqnarray}
\nonumber
& \mathbf{1} = \lleft[\matrix{ 1&0\vspace*{2pt}
\cr
0&1 } \rright],\qquad
\mathbf{e}_1 = \lleft[\matrix{ \ii&0\vspace*{2pt}
\cr
0&-\ii} \rright] ,\qquad
\mathbf{e}_2 = \lleft[\matrix{ 0&1\vspace*{2pt}
\cr
-1&0 } \rright],\qquad
\mathbf{e}_3 = \lleft[\matrix{ 0&\ii\vspace*{2pt}
\cr
\ii&0 } \rright] .
\end{eqnarray}

A quaternion\vspace*{1pt} $ q $ is represented as
$ q = q^{(0)}\mathbf{1} +q^{(1)}\mathbf{e}_1 +
q^{(2)}\mathbf{e}_2 +q^{(3)}\mathbf{e}_3 $.
Here the $ q^{(i)} $ are complex numbers.
There is identification between
the $ 2\ts2 $ complex matrices and
the quaternions given by
%
%
\begin{equation}
\label{"91q} \lleft[\matrix{ a&b \vspace*{2pt}
\cr
c&d} \rright] =
\frac{1}{2}(a+d)\mathbf{1} - \frac{\ii}{2}(a-d)\mathbf{e}_1
+ \frac{1}{2}(b-c)\mathbf{e}_2 - \frac{\ii}{2}(b+c)
\mathbf{e}_3
\end{equation}
or equivalently
%
%
\begin{equation}
\label{"91r} \lleft[\matrix{ q^{(0)}+\ii q^{(1)}&
q^{(2)}+\ii q^{(3)} \vspace*{2pt}
\cr
-q^{(2)}+\ii
q^{(3)} & q^{(0)}-\ii q^{(1)} } \rright] =
q^{(0)}\mathbf{1} +q^{(1)}\mathbf{e}_1 +
q^{(2)}\mathbf{e}_2 +q^{(3)}\mathbf{e}_3
.
\end{equation}
We denote by
$ \Theta(q^{(0)}\mathbf{1} +q^{(1)}\mathbf{e}_1 +
q^{(2)}\mathbf{e}_2 +q^{(3)}\mathbf{e}_3 )$
the $ 2\ts2 $ complex matrix defined by
the left-hand side of (\ref{"91r}).
By definition, $ \Upsilon\bigl(\bigl[ {a\enskip b \atop c\enskip d} \bigr]\bigr) $ is
the quaternion on the right-hand side of (\ref{"91q}).
We also remark that these relations can be naturally extended to
the ones between $ (2N)\ts(2N) $ complex matrices and
$ N\ts N $ quaternion matrices.

For a quaternion $ q = q^{(0)}\mathbf{1} +q^{(1)}\mathbf{e}_1 +
q^{(2)}\mathbf{e}_2 +q^{(3)}\mathbf{e}_3 $,
we call $ q^{(0)} $ the scalar part of $ q $.
A quaternion is called scalar if $ q^{(i)} = 0 $ for
$ i = 1,2,3 $. We often identify a scalar quaternion
$ q = q^{(0)}\mathbf{1}$ with the complex number $ q^{(0)} $
by the obvious correspondence.

Let
$ \bar{q} = q^{(0)}\mathbf{1}-
\{ q^{(1)}\mathbf{e}_1 +q^{(2)}\mathbf{e}_2 +q^{(3)}\mathbf{e}_3 \}.
$
A quaternion matrix $ A = [a_{ij}]$ is called self-dual if $ a_{ij} =
\bar{a}_{ji} $ for all $ i,j $.
For a self-dual $ n\ts n $ quaternion matrix $ A = [a_{ij}]$, we set
%
%
\begin{eqnarray}
\label{"91p}\det A = \sum_{\sigma\in\mathfrak{S}_n }
\operatorname{sign} [\sigma] \prod_{i = 1}^{L(\sigma)}
[a_{\sigma_i(1)\sigma_i(2)} 
\cdots a_{\sigma_i(\ell-1)\sigma_i(\ell)} a_{\sigma_i(\ell)\sigma_i(1)}
]^{(0)} .
\end{eqnarray}
Here $ \sigma= \sigma_1\cdots\sigma_{L(\sigma)}$
is a decomposition of $ \sigma$ to products of
the cyclic permutations $ \{\sigma_i\} $
with disjoint indices. We write
$ \sigma_i = (\sigma_i(1),\sigma_i(2),\ldots,
\sigma_i(\ell)) $, where $ \ell$
is the length of the cyclic permutation $ \sigma_i $.
The decomposition is unique up to the order of
$ \{ \sigma_i \} $.
As before, $[\cdot]^{(0)}$ means the scalar part of the quaternion
$\cdot$.
It is known that the right-hand side is well defined.
See Section~5.1 in \cite{mehta} for details.

For a self-dual $ N\ts N $ quaternion matrix $ A = [a_{ij}]$,
it holds that \cite{mehta}, (5.1.15)
%
\begin{equation}
\label{"91x} \det\Theta( A ) = (\det A)^2 .
\end{equation}
Here $\Theta( A )$ is the $ (2N)\ts(2N)$ complex matrix
given by the relation (\ref{"91q}).
We note that the determinant on the left-hand side of (\ref{"91x})
is of the $ (2N)\ts(2N)$ matrix with complex elements, while
that on the left-hand side of (\ref{"91x})
is of the $ N\ts N $ matrix with quaternion elements.

We are now ready to introduce $ \Ksinb$.
Let
$ S(x) = \sin(\pi x)/ \pi x $,
$ D(x) = \frac{d S }{d x}(x) $ and
$ I(x) = \int_0^x S(y)\,dx $.
Let $ \varepsilon(t) = -1/2 $ $ (t>0)$,
$ \varepsilon(t) = 0 $ $ (t=0)$ and
$ \varepsilon(t) = 1/2 $ $ (t<0)$.
%
%
\begin{eqnarray}
\label{"91s} \Ksinx(x,y)& =& \Upsilon\biggl( \lleft[\matrix{ S(x-y)& D(x-y)
\vspace*{2pt}
\cr
I(x-y) -\varepsilon(x-y) & S(x-y) } \rright]\biggr),
\\
\label{"91u}\Ksiny(x,y) &= &S(x-y),
\\
\label{"91t} \Ksinz(x,y)& =& \Upsilon\biggl( \lleft[\matrix{ S\bigl(2(x-y)
\bigr)& D\bigl(2(x-y)\bigr) \vspace*{2pt}
\cr
I\bigl(2(x-y)\bigr) & S\bigl(2(x-y)
\bigr) } \rright]\biggr) .
\end{eqnarray}
We thus clarify the meaning of (\ref{224}).

It is known that the matrices
$ [\Ksinb(x_i,x_j)]_{1\le i,j \le n}$ are self-dual
($ \beta= 1,4$), and that there exist unique random point fields
$ \mu_{\mathrm{dys}, \beta}$ $ (\beta= 1,2,4) $
whose correlation functions $ \{ \rho^n\} $ are given by (\ref{224});
see \cite{mehta}, Chapters 5--8.


%
\begin{lem} \label{l"91}
$ \mu_{\mathrm{dys}, \beta}$ $( \beta= 1,2,4 )$ satisfy \textup{(\hyperlink{A1}{A.1})}.
\end{lem}
\begin{pf}
Since the correlation functions $ \{ \rho^n\} $
have the expression (\ref{224}) and
the kernels $ \Ksinb$ are bounded,
we see that $ \{ \rho^n\}$ satisfy (\hyperlink{A1}{A.1}).
\end{pf}

To prove the quasi-Gibbs property of $ \mu_{\mathrm{dys}, \beta}$,
it is sufficient to check
(\hyperlink{A4}{A.4}) and~(\hyperlink{A5}{A.5}) by Theorem~\ref{l41}.
Therefore, the problem is to construct a finite-particle approximation
$ \{\muN\}$ fulfilling the assumptions in (\hyperlink{A4}{A.4}) and (\hyperlink{A5}{A.5}).
We will take $ \{\muN\}$, whose potentials satisfy
(\ref{61s}) and (\ref{61u}) for $ \beta=1,2,4 $.
Hence, we assume $ n_{N}= 2^{4N} $ and $ \IN= (-N,N)$ as in (\ref{61u}).
We take $ \Phi(x) = 0 $ and $ \PhiN(x) = - \log1_{\IN}(x) $.
$ \Psi$ and $ \PsiN$ are the same as in
(\ref{61a}), (\ref{61s}) and~(\ref{61u}).


To introduce the finite-particle approximation $ \{\muN\}$
we first recall some facts about circular ensembles $ \{\nuN\}$.
Let $ \check{\nu}^{N} $ denote the probability measure on $ \R^{n_{N}} $ defined by
%
%
\begin{equation}
\label{221}\quad d\check{\nu}^{N} =  \frac{1}{Z} \prod
_{i = 1}^{n_{N}} 1_{\TN}( x_i)
\prod_{i, j = 1, i < j }^{n_{N}} \bigl|e^{2\pi\ii x_i/n_{N}}-
e^{2\pi\ii x_j/n_{N}}\bigr|^{\beta} \,dx_1 \cdots \,dx_{n_{N}}
,
\end{equation}
where $ Z $ is the normalization and
$ \TN= (-n_{N}/2,n_{N}/2] $.
It is well known \cite{mehta,for} that
the distribution of
$(e^{2\pi\ii x_i/n_{N}})_{1\le i \le n_{N}}$
under $ \check{\nu}^{N} $ is equal to the distributions
of the spectra of
the circular orthogonal, unitary and symplectic ensembles
for $ \beta= 1 $, $ 2 $ and $ 4 $, respectively.

Let $ \iota$ be a map such that $ \iota((x_i)) = \sum_{i} \delta_{x_i}$.
Set $ \nu^{N}= \check{\nu}^{N} \circ\iota^{-1}$, and let $ \varrho_{N}^{n}$ denote
the $ n $-correlation function of $ \nu^{N}$.
Then by (\ref{221}), we see that $ \varrho_{N}^{n}= 0 $ for $ n >
n_{N}$ and
\begin{eqnarray*}
\varrho_{n_{N}}^{n_{N}}(x_1,
\ldots,x_{n_{N}}) = \frac{n_{N}!}{Z} \prod_{i, j = 1, i < j }^{n_{N}}
1_{\TN}(x_i) \bigl|e^{2\pi\ii x_i/n_{N}}- e^{2\pi\ii x_j/n_{N}}\bigr|^{\beta}
1_{\TN}(x_j) .
\end{eqnarray*}
For each $ n\in\N$, the $ n $-correlation function
$ \varrho_{N}^{n}$ can be written as (see \cite{mehta}, (11.1.10))
%
\begin{equation}
\label{71kk} \varrho_{N}^{n}(x_1,
\ldots,x_n) = \det\bigl[1_{\TN}(x_i)
\KsinNb(x_i-x_j)1_{\TN}(x_j)
\bigr]_{1 \le i,j \le n} ,
\end{equation}
where $ \KsinNb$ is given by (\ref{"91s})--(\ref{"91t}) with the
replacement of\vspace*{1pt}
$ S(x) $, $ D(x) $, and $ I(x) $ by
$ S_N (x) $, $ D_N(x) $, and $ I_N(x) $, respectively.
Here $ S_N $ is defined as
%
%
\begin{equation}
\label{71zz} S_N (x) = \frac{1}{n_{N}}\frac{\sin(\pi x)}{\sin(\pi x/n_{N})} .
\end{equation}
Moreover, we set $ D_N(x) = {dS_N (x)}/{dx} $ and
$ I_N(x) = \int_0^x S_N (y)\,dy $.
One can easily deduce (\ref{71kk}) and (\ref{71zz}) from the
results in \cite{mehta}, Chapter 11, combined with the scaling
$ \theta\mapsto2\pi x /n_{N}$. Indeed, these follow from
(11.1.5), (11.1.6),
(11.3.16), (11.3.22), (11.3.23),
(11.5.6) and
(11.5.13) in \cite{mehta}.\footnote{$ IS_{2N} $ in (11.1.6) of \cite{mehta}
should be $ I_{2N} $. }

We are now ready to introduce the finite-particle approximation $ \{
\muN\} $.
%
\begin{lem} \label{l82}
Let $ \muN= \nu^{N}\circ\pi_{\IN}^{-1} $.
Then $ \mu_{\mathrm{dys}, \beta}$ satisfy \textup{(\hyperlink{A4}{A.4})} with $ \muN$.
Here we take $ \PhiN(x) = -\log1_{\IN}(x) $, and
$ \PsiN$ is given by (\ref{61s}) and (\ref{61u}).
\end{lem}
\begin{pf}
Let $ \rbN$ be the $ n $-correlation function of $ \muN$.
Then by $ \muN= \nu^{N}\circ\pi_{\IN}^{-1}$, we have
$ \rbN(x_1,\ldots,x_n) = \varrho_{N}^{n}(x_1,\ldots,x_n) $ on $ \IN^n $.
Hence, by (\ref{71kk}), we see that $ \rbN$ satisfy
%
\begin{equation}
\label{71ka}  \rbN(x_1,\ldots,x_n) = \det
\bigl[1_{\IN}(x_i) \KsinNb(x_i-x_j)1_{\IN
}(x_j)
\bigr]_{1 \le i,j \le n} .
\end{equation}
By (\ref{71zz}) and (\ref{"91s})--(\ref{"91t}), we deduce that
$ 1_{\IN}(x) \KsinNb(x-y)1_{\IN}(y) $ converge compact uniformly to
$ \KsinNb(x-y)$.
This combined with (\ref{71ka}) yields (\ref{41a}).

Let $ k^{N,n}_i (\mathbf{x}_n)$ be the norm of the $ i $th
row vector of
$ [1_{\IN}(x_i) \KsinNy(x_i-x_j)1_{\IN}(x_j) ]_{1 \le i,j \le n} $, where
$ \mathbf{x}_n=(x_1,\ldots,x_n)$.
Then there exists a constant
$ \Ct\label{;71}$ such that $ | k^{N,n}_i (\mathbf{x}_n) | \le
c_{\scriptsize\ref{;71}} n^{1/2}$
because the kernels $ \KsinNy$ are uniformly bounded.
Hence, we have
%
\begin{eqnarray}
 \label{71b} \qquad\bigl| \det\bigl[1_{\IN}(x_i)
\KsinNy(x_i-x_j)1_{\IN}(x_j)
\bigr]_{1 \le i,j \le
n} \bigr| \le k^{N,n}_i (\mathbf{x}_n)^n
\le c_{\scriptsize\ref{;71}}^n n^{n/2} .
\end{eqnarray}
This combined with (\ref{71ka}) yields (\ref{41b}) with $ \beta=2 $.
We can prove (\ref{41b}) for $ \beta= 1,4$, similarly using identity
(\ref{"91x}) on the quaternion determinant.
We thus obtain (1) of~(\hyperlink{A4}{A.4}).

(2) of (\hyperlink{A4}{A.4}) is clear because $ \muN= \nu^{N}\circ
\pi_{\IN}^{-1} $,
and $ \nuN$ consists of $ n_{N}$ particles.

Let $ \hat{\Phi}^N(x)=-\log1_{\TN}(x)$ and
$ \hat{\Psi}^N(x,y)= -\beta\log| e^{2\pi\ii x/n_{N}}- e^{2\pi
\ii y/n_{N}}|$.
Then by (\ref{221}), we see that $ \nu^{N}$ are
$ ( \hat{\Phi}^N ,\hat{\Psi}^N ) $-canonical Gibbs measures.
Clearly, $ \hat{\Psi}^N(x,y)=\PsiN(x,y)$ for $ x,y\in\IN$.
Hence, $ \muN$ are $(\PhiN,\PsiN)$-canonical Gibbs measures because
$ \muN= \nu^{N}\circ\pi_{\IN}^{-1}$.\vspace*{1pt}

(4) of (\hyperlink{A4}{A.4}) is obvious through construction.
\end{pf}

We next proceed with the proof of (\hyperlink{A5}{A.5}).
For this, it is sufficient to prove (\hyperlink{A6}{A.6}) by Theorem~\ref{l63}.
We note that (\hyperlink{A6}{A.6}) consists of two conditions: (\ref{A7b})
and (\ref{A7a}).
We prove (\ref{A7b}) in the next four lemmas.

Let $ \IN= (-N,N)$ and $ n_{N}= 2^{4N}$, as before.
By (\ref{61u}), we easily see the following: 
%
\begin{eqnarray}
\label{83p} \wNx& =& \frac{n_{N}}{2\pi} \sin\frac{2\pi x}{n_{N}} + \ii
\frac{n_{N}}{2\pi}\biggl( 1 - \cos\frac{2\pi x}{n_{N}}\biggr)
\nonumber
\\[-8pt]
\\[-8pt]
\nonumber
& =& \frac{n_{N}}{\pi} \sin\frac{\pi x}{n_{N}}\cos\frac{\pi x}{n_{N}} +
\ii\frac{n_{N}}{\pi}\sin^2 \frac{\pi x}{n_{N}} \qquad\mbox{for } x\in\IN ,
\\
\label{83f} \bigl|\wNx\bigr| &=& \frac{n_{N}}{\pi} \biggl|\sin\frac{\pi x}{n_{N}}\biggr|\qquad \mbox{for } x
\in\IN .
\end{eqnarray}
Hence, by (\ref{83p}) and (\ref{83f}), we have
%
\begin{equation}
\label{83q} \frac{\wNx}{|\wNx|} = \frac{\sin{\pi x}/{n_{N}} \cos{\pi x}/{n_{N}} } {
|\sin{\pi x}/{n_{N}}|} + \ii \biggl|\sin
\frac{\pi x}{n_{N}}\biggr| \qquad\mbox{for } x \in\IN .
\end{equation}

%
\begin{lem} \label{l82z}
Let $ \wN$ be as in (\ref{61u}). Let $ \Sti_{1\infty} = \{ 1\le|
x |< \infty\}$. 
Then the following holds: %
%
\begin{eqnarray}
\label{82a}  \sup_{N\in\N} \sup_{x\in\Sti_{1\infty} } \frac{|x|}{| \wNx|}
&<&
\infty ,\qquad \sup_{N\in\N} \sup_{x\in\Sti_{1\infty} } \frac{\lceil| \wNx
|\rceil}{| \wNx|} < \infty ,
\\
\label{82b} \qquad \sup_{N\in\N} \sup_{x\in\Sti_{1\infty} } \bigl| \bar{\varpi
}_{N}(x)- x \bigr| &<& \infty,\qquad \sup_{N\in\N} \sup_{x\in\Sti_{1\infty} } \bigl|
\bigl\lceil\bigl|\bar{\varpi }_{N}(x)\bigr| \bigr\rceil- |x| \bigr| < \infty .
\end{eqnarray}
\end{lem}
\begin{pf}
Note that, if $ x \in\IN$, then $ | \wNx|$ is the length of the
segment between the origin and $\wNx$, and that $ |x| $ is the length
of the arc connecting these two points on the circle centered at
$ \ii n_{N}/ 2\pi$ with radius $ n_{N}/ 2\pi$.
This implies $ | \wNx| < |x|$ for $ x\in\Sti_{1\infty} \cap\IN$.
By definition, $ \wN$ is linear on $ \INN\setminus\IN$ and $ \wNx
= x $ on $ \INN^{c} $.
Hence, the maximum of $ \frac{|x|}{| \wNx|}$ over $ \Sti_{1\infty} $
is attained at $ x= \pm1$. Therefore, we have
\[
\sup_{N\in\N} \sup_{x\in\Sti_{1\infty} } \frac{|x|}{| \wNx|} =
\sup_{N\in\N} \frac{1}{| \wN(1) |} = \sup_{N\in\N} \frac{\pi}{n_{N}}
\frac{1}{ \sin {\pi
}/{n_{N}} } = \frac{\pi}{2^4} \frac{1}{\sin{\pi}/{2^4 }} < \infty .
\]
The second inequality in (\ref{82a}) follows from
$ \lceil| \wNx|\rceil< |\wNx| + 1 $ and the first inequality.
We thus obtain (\ref{82a}).

Direct calculation shows that there exists a constant $\Ct\label{;82f}$
independent of $ N $ such that
%
\begin{equation}
\label{82f} \sup_{ x \in\R} \bigl|\bigl\{\wNx- x \bigr\}'\bigr|= \bigl|
\wN'(N) - 1 \bigr|\le c_{\scriptsize\ref
{;82f}}N 2^{-4N} .
\end{equation}
Since $ \wN(0) = 0$ and $ \wNx= x $ for $ |x|\ge N+1$, (\ref{82f}) yields
%
\begin{equation}
\label{82g} \sup_{ x \in\R}\bigl |\wNx- x \bigr| \le c_{\scriptsize\ref{;82f}} N(N+1)
2^{-4N} .
\end{equation}
Because $ |\wNx- x |= |\bar{\varpi}_{N}(x)- x |$, (\ref{82g})
allows us to deduce
the first inequality in (\ref{82b}). The second is clear from the first.
\end{pf}

We set $ \Ct\label{;83}=\sup_{ x\in\Sti_{1r} , N\in\N} {\lceil
| \wNx|\rceil}/{| \wNx|} < \infty$.%
\begin{lem} \label{l83}
Let $ \map{u^N_r }{\R}{\mathbb{C}}$ be such that
\[
u^N_r (x) = { 1_{\Sti_{1r} } (x)\bigl \lceil\bigl| \wNx\bigr|
\bigr\rceil}/{ \bar{\varpi}_{N}(x)} .
\]
Then
%
\begin{eqnarray}
\label{83z} \sup_{N\in\N} \int_{\R}\bigl|u^N_r
\bigr|^2\,dx &\le&2c_{\scriptsize\ref{;83}}^2 r 
,
\\
\label{83b} \sup_{r\in\N} \sup_{N\in\N}\biggl |\int
_{\R} 1_{\IN} u^N_r \,dx\biggr|&<&
\infty .
\end{eqnarray}
\end{lem}

\begin{pf}
From $ |u^N_r | \le c_{\scriptsize\ref{;83}} 1_{\Sti_{1r} }$ and $ \Sti_{1r}
\subset(-r,r)$, (\ref{83z}) is obvious.\vspace*{1pt}

Through construction, we deduce that
$ |\wNx|$ and $ \lceil| \wNx|\rceil$ are even functions.
Moreover, $ \Re[\wNx] $, the real part of $ \wNx$, is an odd
function in $ x\in\R$.
Hence, so is $\Re[1/{ \bar{\varpi}_{N}(x)}] = \Re[\wNx/{ |\wNx|^2
}] $.
Collecting these, we see that
$ \Re[ 1_{\IN} u^N_r] = 1_{\IN} 1_{\Sti_{1r} } \lceil| \wN
|\rceil\Re[ { 1/\bar{\varpi}_{N}}] $ becomes an odd function.
Hence, we have $ \int_{\R}\Re[ 1_{\IN} u^N_r ] \,dx = 0 $. Therefore,
it only remains to estimate
$ \Im[ 1_{\IN} u^N_r]$, the imaginary part of $ 1_{\IN} u^N_r $.

Note that
$ u^N_r (x) = { 1_{\Sti_{1r} } (x) \lceil| \wNx|\rceil} \wNx/{
|\wNx|^2 } $.
We easily see that $ \Im[u^N_r ] \ge0 $ and
$ \Im[\wNx/ |\wNx| ] $ takes its maximum at $ x=\pm N $
according to~(\ref{83p}).
Then by~(\ref{83q}), we have
\[
 \sup_{x\in\R} \bigl|\Im\bigl[ 1_{\IN} u^N_r
(x) \bigr]\bigr| \le c_{\scriptsize\ref{;83}} \sin \frac{\pi N }{n_{N}} \le c_{\scriptsize\ref{;83}}
\frac{\pi N}{2^{4N}} .
\]
Clearly, $\Im[ 1_{\IN}u^N_r ]= 1_{\IN} \Im[u^N_r ]= 0 $ for $ x\not
\in\IN$. Therefore, we deduce that
\[
\sup_{r\in\N}\sup_{N\in\N} \int_{\R} \bigl|\Im
\bigl[1_{\IN} u^N_r (x) \bigr]\bigr| \,dx \le
\sup_{N\in\N} 2N c_{\scriptsize\ref{;83}} \frac{\pi N}{2^{4N}} < \infty .
\]
This implies (\ref{83b}).
\end{pf}

%
\begin{lem} \label{l84}
Let $ u^N_r $ be as in Lemma~\ref{l83}.
Set $ \mathsf{u}^N_r(\mathsf{x}) = \sum_{i}u^N_r (x_i) $.
Then
%
\begin{equation}
\label{aaaa}\lim_{r\to\infty}r^{-3/4 } \sup_{N\in\N} \bigl\|
\mathsf{u}^N_r \bigr\|_{L^2(\SSS,\muN)} = 0 .
\end{equation}
\end{lem}
\begin{pf}
We set
$ \hat{\mathsf{u}}^N_r(\mathsf{x}) = \sum_{i} 1_{\IN} (x_i) u^N_r
(x_i) $.
Then from $ \muN= \nu^{N}\circ\pi_{\IN}^{-1}$ we see that
$ \|\mathsf{u}^N_r \|_{L^2(\SSS,\muN)} = \| \hat{\mathsf{u}}^N_r \|_{L^2(\SSS,\nu^{N})}$.
Hence, (\ref{aaaa}) follows from
%
%
\begin{equation}
\label{84b} \lim_{r\to\infty}r^{-3/4 } \sup_{N\in\N} \bigl\| \hat{
\mathsf {u}}^N_r \bigr\|_{L^2(\SSS,\nu^{N})}= 0 .
\end{equation}

We note that $ \varrho_{N}^{1}(0) = 1 $ according to (\ref{71kk}).
We write $ 1_{\IN}u^N_r = \hat{u}^N_{r,1} + \ii\hat{u}^N_{r,2} $,
where $ \hat{u}^N_{r,m} $ $ (m=1,2)$ are real valued.
We denote by $ \mathcal{F}_N $ the Fourier transform defined before
Lemma~\ref{l61a}. Let $ \mN(\xi) $ be as in (\ref{63c}).
Applying Lemma~\ref{l61a} to $ \hat{\mathsf{u}}^N_r $ and using
$ \varrho_{N}^{1}(0) = 1 $, we have
\begin{eqnarray*}
\bigl\| \hat{\mathsf{u}}^N_r
\bigr\|_{L^2(\SSS,\nu^{N})}^2 = \sum_{m=1}^2
\biggl[ \biggl(\int_{\TN} \hat{u}^N_{r,m}
\,dx \biggr)^2 + \frac{1}{n_{N}}\sum_{\6}
\bigl|\mathcal{F}_{N}\bigl( \hat{u}^N_{r,m} \bigr)
\bigr|^2 (\xi) \mN(\xi) \biggr] .
\end{eqnarray*}
From (\ref{83b}), we deduce that
$ \limi{r}r^{-3/4} \sup_{N\in\N} | \int_{\TN} \hat{u}^N_{r,m} \,dx
| = 0$
for $ m=1,2$.
Hence, it only remains for (\ref{aaaa}) to prove that for $ m=1,2 $,
%
\begin{equation}
\label{84n} \lim_{r\to\infty}r^{-3/2 } \sup_{N\in\N}
\frac{1}{n_{N}}\sum_{\6} \bigl|\mathcal{F}_{N}
\bigl( \hat{u}^N_{r,m} \bigr) \bigr|^2 (\xi) \mN(\xi)
= 0 .
\end{equation}

Let
$ P_{N} = \{ -\frac{n_{N}+1}{2}+p ; 1 \le p \le n_{N}, p \in\N\}$.
Then by an elementary calculation of the triangle series, we have
an expansion of $ S_N (x) $ such that
%
\begin{equation}
\label{71y} S_N (x) = \frac{1}{n_{N}} \sum
_{p\in P_{N}} e^{2\pi x p\ii/n_{N}} .
\end{equation}
This together with $ D_N(x) = {dS_N (x)}/{dx} $ and
$ I_N(x) = \int_0^x S_N (y)\,dy $ yields
%
\begin{eqnarray}
\label{71v}  D_N(x) &= &\frac{2\pi\ii}{n_{N}^2} \sum
_{p\in P_{N}} p e^{2\pi x p\ii/n_{N}},
\\
\label{71v1}  I_N(x) &=& \frac{1 }{2\pi\ii} \sum
_{p\in P_{N}} \frac{1}{p} \bigl( e^{2\pi x p\ii/n_{N}} - 1 \bigr)
= \frac{1 }{2\pi\ii} \sum_{p\in P_{N}}
\frac{1}{p} e^{2\pi x p\ii/n_{N}} .
\end{eqnarray}
For (\ref{71v1}) we use $ 0 \notin P_{N} $,
which follows from $ n_{N}/2 \in\N$.

Let $ \T$ be the two-cluster function of $ \nu^{N}$ defined by (\ref{60T}).
Let $ \mathcal{T}^{N}_{\beta}(x) $ be the $ n_{N}$-periodic function
such that $ \mathcal{T}^{N}_{\beta}(x) = \T(x,0)$ for $ x\in\TN$.
Then through construction [see (\ref{"91p}), (\ref{71kk})],
%
\begin{equation}
\label{71d} \mathcal{T}^{N}_{\beta}(x) = \bigl[\KsinNb(x)
\KsinNb(-x)\bigr]^{(0)}\qquad \mbox{for $ x \in\TN$} .
\end{equation}
Let $ P_{N,1} = P_{N,2} = P_{N}$ and
$ P_{N,4} = \{ p + \frac{1}{2} ; p \in\N,
-n_{N}\le N < n_{N}\} $.
Then (\ref{71y})--(\ref{71d}) combined with
the definition of $ \KsinNb$ yield
%
%
\begin{eqnarray}
\label{71a1} \mathcal{T}^{N}_{2}( x )& =& \bigl|\KsinNy( x
)\bigr|^2 = \frac{1}{n_{N}^2}\biggl|\sum_{p\in P_{N,2}}
e^{2\pi x p\ii/n_{N}}\biggr|^2,
\\
\label{71a3}\mathcal{T}^{N}_{1}( x ) &=&
\frac{1}{n_{N}^2}\biggl|\sum_{p\in P_{N,1}} e^{2\pi x p\ii/n_{N}}\biggr|^2-
\frac{1}{n_{N}^2}\sum_{p,q\in P_{N,1}}\frac{p}{q}
e^{2\pi x (p+q)\ii/n_{N}},
\\
\label{71a4}\mathcal{T}^{N}_{4}( x ) &=&
\frac{1}{n_{N}^2}\biggl|\sum_{p\in P_{N,4}} e^{4\pi x p\ii/n_{N}}\biggr|^2-
\frac{1}{n_{N}^2} \sum_{p,q\in P_{N,4}}\frac{p}{q}
e^{4\pi x (p+q)\ii/n_{N}} .
\end{eqnarray}
For the reader's convenience,
we provide more details of the proof of
(\ref{71a3}) and~(\ref{71a4})
as an Appendix \ref{sA3}.

We now consider the Fourier series
$ \mathcal{F}_{N}(\mathcal{T}^{N}_{\beta})( \xi) =
\int_{\TN} e^{-2\pi\ii\xi\cdot x }
\mathcal{T}^{N}_{\beta}( x )\,dx$.
By~(\ref{71a1})--(\ref{71a4}), we obtain
\[
\sup_{N\in\N} \sup_{\6} \bigl|\mathcal{F}_{N}\bigl(
\mathcal{T}^{N}_{\beta}\bigr) (\xi)\bigr|<\infty .
\]
So $\mN$ defined by (\ref{63c}) for $\nu^{N}$ satisfies
\[
\Ct\label{;63c}:= \sup_{N\in\N} \sup_{\6} \bigl| \mN(\xi)\bigr | <
\infty .
\]
From the isometry of the Fourier series and (\ref{83z}), we have that
for $ m=1,2 $,
\begin{eqnarray*}
\sup_{N\in\N} \biggl\{ \frac{1}{n_{N}}\sum
_{\6} \bigl|\mathcal{F}_{N}\bigl( \hat
{u}^N_r \bigr)\bigr|^2(\xi) \biggr\} =
\sup_{N\in\N} \biggl\{ \int_{\TN} \bigl |
\hat{u}^N_{r,m} \bigr|^2 \,dx \biggr\} \le
2c_{\scriptsize\ref{;83}}^2 r .
\end{eqnarray*}
Combining these two equations, we obtain
\begin{eqnarray*}
\sup_{N\in\N} \biggl\{ \frac{1}{n_{N}}\sum
_{\6} \bigl|\mathcal{F}_{N}\bigl(\hat
{u}^N_{r,m} \bigr) \bigr|^2 (\xi) \mN(\xi) \biggr\}
\le c_{\scriptsize\ref{;63c}} 2c_{\scriptsize\ref{;83}}^2 r\qquad (m=1,2) ,
\end{eqnarray*}
which yields (\ref{84n}). We thus complete the proof.
\end{pf}

%
\begin{lem} \label{l86}
Let $ u^N_r = 1_{\Sti_{1r} } {\lceil| \wN|\rceil}/{ \bar{\varpi
}_{N}}$ be
as in Lemma~\ref{l83}. Then
%
\begin{equation}
\label{86a} \limi{r} r^{-3/4} \sup_{N\in\N} \int
_{\R} \Bigl\{\sup_{M\in\N} \bigl|u^M_r-u^N_r
\bigr| \Bigr\} \rNone \,dx =0 .
\end{equation}
\end{lem}
\begin{pf}
Through straightforward calculation, we have
\begin{eqnarray*}
u^M_r-u^N_r
&=& 1_{\Sti_{1r} } \biggl\{ \frac{\lceil| \wM|\rceil}{\bar{\varpi}_{M}} - \frac{\lceil| \wN|\rceil}{\bar{\varpi}_{N}} \biggr\}
\\
&=& 1_{\Sti_{1r} } \biggl\{ \frac{\lceil| \wM|\rceil}{ \bar{\varpi
}_{M}\bar{\varpi}_{N}} (\bar{
\varpi}_{N}-\bar{\varpi}_{M}) + \frac{1}{ \bar{\varpi}_{N}} \bigl(
\bigl\lceil| \wM|\bigr\rceil-\bigl\lceil| \wN |\bigr\rceil\bigr) \biggr\}
\\
& =& 1_{\Sti_{1r} } \biggl\{ \frac{\lceil| \wM|\rceil}{ \bar{\varpi
}_{M}\bar{\varpi}_{N}} (\bar{
\varpi}_{N} - x + x -\bar{\varpi}_{M}) \\
&&\qquad{}+ \frac{1}{ \bar{\varpi}_{N}}
\bigl(\bigl\lceil| \wM|\bigr\rceil-|x| + |x|-\bigl\lceil| \wN|\bigr\rceil\bigr) \biggr\} .
\end{eqnarray*}
Applying (\ref{82a}) and (\ref{82b}) to the last line,
we have a constant $ \Ct\label{;86a}$ such that
%
\begin{equation}
\label{86c} \bigl| u^M_r (x) -u^N_r
(x) \bigr| \le c_{\scriptsize\ref{;86a}} 1_{\Sti_{1r} } (x) \frac
{ 1 }{|x|}\qquad \mbox{for
all } x \in\Sti_{1r} , M,N\in\N .
\end{equation}
By definition, $ u^M_r (x) = 0 $ on $ \Sti_{1r}^c $.
Hence, by (\ref{86c}) and $ \rNone(x) \le1 $, we obtain
\[
\sup_{N\in\N} \int_{\R} \Bigl\{
\sup_{M\in\N} \bigl|u^M_r-u^N_r
\bigr| \Bigr\} \rNone \,dx  \le c_{\scriptsize\ref{;86a}} \int_{\Sti_{1r} }
\frac{ 1 }{|x|} \,dx 
= c_{\scriptsize\ref{;86a}} 2 \log r .
\]
This deduces (\ref{86a}).
\end{pf}

%
\begin{lem} \label{l85}
Let $ \mathsf{u}^N_r $ be as in Lemma~\ref{l84}. Then %
%
\begin{equation}
\label{85z} \lim_{r\to\infty} r^{-3/4 } \sup_{N\in\N} \Bigl\|
\sup_{M\in\N} \bigl|\mathsf{u}^M_r \bigr|
\Bigr\|_{L^1(\SSS,\muN)} = 0 .
\end{equation}
\end{lem}
\begin{pf}
We note that
$\sup_{M\in\N} |\mathsf{u}^M_r | \le
\{\sup_{M\in\N} |\mathsf{u}^M_r - \mathsf{u}^N_r|\} + |\mathsf{u}^N_r|$.
Hence,
\begin{eqnarray*}
\sup_{N\in\N} \Bigl\| \sup_{M\in\N} \bigl|\mathsf{u}^M_r
\bigr| \Bigr\|_{L^1(\SSS
,\muN)} \le 
\sup_{N\in\N} \Bigl\|\sup_{M\in\N}
\bigl| \mathsf{u}^M_r - \mathsf{u}^N_r
\bigr| \Bigr\|_{L^1(\SSS,\muN)} + \sup_{N\in\N} \bigl\|\mathsf{u}^N_r
\bigr\|_{L^1(\SSS,\muN)} .
\end{eqnarray*}
By Lemma~\ref{l84} and H\"{o}lder's inequality, we have
\[
\limi{r} r^{-3/4 } \sup_{N\in\N} \bigl\|\mathsf{u}^N_r
\bigr\|_{L^1(\SSS
,\muN)} = 0 .
\]
Hence, it only remains to prove
%
%
\begin{equation}
\label{85} \lim_{r\to\infty}r^{-3/4 } \sup_{N\in\N}\Bigl \|
\sup_{M\in\N}\bigl |\mathsf{u}^M_r -
\mathsf{u}^N_r \bigr| \Bigr\|_{L^1(\SSS,\muN)} = 0 .
\end{equation}

We write $ \mathsf{x}=\sum_i \delta_{x_i}$. It is then obvious that
\begin{eqnarray*}
\sup_{M\in\N} \bigl|\mathsf{u}^M_r(\mathsf{x}) -
\mathsf{u}^N_r (\mathsf{x}) \bigr| &= &\sup_{M\in\N} \biggl|
\sum_{i}\bigl\{ u^M_r
(x_i)- u^N_r (x_i)\bigr\} \biggr|
\\
& \le& \sum_{i} \sup_{M\in\N}
\bigl|u^M_r (x_i)- u^N_r
(x_i)\bigr| .
\end{eqnarray*}
Taking the expectation of both sides w.r.t. $ \muN$, we deduce that
%
\begin{equation}
\label{85p}
\Bigl\|\sup_{M\in\N} \bigl|\mathsf{u}^M_r -
\mathsf{u}^N_r \bigr|\Bigr\|_{L^1(\SSS
,\muN)}  \le \int
_{\R} \Bigl\{ \sup_{M\in\N} \bigl|u^M_r-u^N_r
\bigr| \Bigr\} \rNone \,dx .
\end{equation}
Combining (\ref{85p}) with (\ref{86a}), we obtain (\ref{85}),
which completes the proof of Lemma~\ref{l85}.
\end{pf}

\begin{pf*}{Proof of Theorem~\ref{l22}}
According to Theorem~\ref{l41}, it is enough for (\hyperlink{A2}{A.2}) to check
(\hyperlink{A4}{A.4}) and (\hyperlink{A5}{A.5}).
We have already checked (\hyperlink{A4}{A.4}) by Lemma~\ref{l82}.
By Theorem~\ref{l63}, it is sufficient for (\hyperlink{A5}{A.5}) to prove
(\hyperlink{A6}{A.6}).
(\hyperlink{A6}{A.6}) consists of two conditions: (\ref{A7b}) and (\ref{A7a}).

According to (\ref{82a}), there exists a constant $ \Ct\label
{;88}$ such that
${1 }/{|\wM(x) |^{2}}\le c_{\scriptsize\ref{;88}} /{x^{2}} $ for all $ M\in\N$
and $ x\in\Sti_{1r} $.
This combined with $ \rNone(x )\le1 $ yields
%
\begin{eqnarray}
\label{88} \sup_{N\in\N} \int_{\Sti_{1r} } \biggl\{
\sup_{M\in\N} \frac{1 }{|\wM(x) |^{2}} \biggr\} \rNone(x )\,dx \le
c_{\scriptsize\ref{;88}} \sup_{N\in\N} \int_{\Sti_{1r} }
\frac{1 }{x^{2}}\,dx < \infty .
\end{eqnarray}
Hence, (\ref{71a}) is satisfied with $ \ell_0= 2$,
and thus, we conclude (\ref{A7b}) by Lemma~\ref{l71}.
By Lemma~\ref{l85}, we have (\ref{81b})
with $ c_{\scriptsize\ref{;r71}}=1/4$ and $ \ell_0 = 2$, which yields (\ref{A7a})
by Lemma~\ref{l81}.
\end{pf*}

\section{\texorpdfstring{Proof of Theorem~\protect\ref{l23}}{Proof of Theorem 2.3}} \label{s10}
In this section, we prove the quasi-Gibbs property~(\hyperlink{A2}{A.2})
of the Ginibre random point field $ \mug$ (Theorem~\ref{l23}).
Therefore, we set $ \SSSS= \mathbb{C}$, $ \Phi(z) = |z|^2 $ and
$ \Psi(z_1,z_2) = - 2 \log|z_1-z_2| $.
From Theorems~\ref{l41} and~\ref{l63}, we deduce (\hyperlink{A2}{A.2}) from
(\hyperlink{A4}{A.4}) and (\hyperlink{A6}{A.6}).
Therefore, our task is to check these two assumptions.
We begin with the finite-particle approximation $ \mugN$.

Let $ \muNg$ be the determinantal random point field with kernel $
\kgN$
given by
%
\begin{equation}
 \label{0-1a} \kgN(z_1,z_2) = \frac{1}{\pi} \exp
\biggl\{ - \frac{|z_1|^2}{2}-\frac{|z_2|^2}{2}\biggr\} \Biggl\{\sum
_{k = 0}^{N-1} \frac{1}{k!}(z_1
\cdot\bar{z}_2)^k\Biggr\} .
\end{equation}
Then, by definition, its $ n $-point correlation function $ \rgN$ is
given by
%
\begin{equation}
 \label{0-1b} \rgN(z_1,\ldots,z_n) = \det\bigl[
\kgN(z_i,z_j)\bigr]_{1 \le i,j \le n} .
\end{equation}
It is well known (see, e.g., page 943 in \cite{so-}) that
%
\begin{equation}
\label{0-1c} \muNg\bigl(\mathsf{s} (\mathbb{C} )=N\bigr) = 1 .
\end{equation}
Let $ \muNgcheck$ be the probability measure on $ \mathbb{C}^N $
associated with $ \muNg$. By definition, $ \muNgcheck$ is the
symmetric measure satisfying $ d\muNg= \muNgcheck\circ\iota^{-1} $, where
$ \iota((z_1,\ldots,z_n)) = \sum_{i=1}^{N} \delta_{z_i}$.
It is well known (see, e.g., page 943 in \cite{so-}) that
%
%
\begin{equation}
\label{0-1d}\muNgcheck= \frac{1}{Z}e^{-\sum_{i = 1}^{N}|z_i|^2} \prod
_{1\le i<j\le N} |z_i-z_j|^2
\,dz_1\cdots \,dz_{N} .
\end{equation}
%

%
\begin{lem} \label{l0-2}
$\{ \muNg\}_{N\in\N}$ satisfy \textup{(\hyperlink{A4}{A.4})}.
\end{lem}
\begin{pf}
It is clear that the kernels $ \kgN$ converge to $ \kg$ compact
uniformly as $ N \to\infty$.
Hence, (\ref{41a}) follows from (\ref{0-1b}).
Let $ k^{N,n}_i (z_1,\ldots,z_n)$ be the norm of the $ i $th
row vector of the matrix $ [ \kgN(z_i,z_j)]_{1 \le i,j \le n}$.
We see that $ k^{N,n}_i (z_1,\ldots,z_n)\le n^{1/2}/\pi$ because
$ |\kgN(z_1,z_2)|\le1/\pi$ by (\ref{0-1a}). Hence, we obtain %
%
\begin{equation}
\label{0-0a} \bigl|\det\bigl[\kgN(z_i,z_j)
\bigr]_{1 \le i,j \le n}\bigr| \le \prod_{i=1}^{n}k^{N,n}_i
(z_1,\ldots,z_n) \le \frac{n^{n/2}}{\pi^{n} } .
\end{equation}
Therefore, we deduce (\ref{41b}) from (\ref{0-1b}) and (\ref{0-0a}).
We thus see that (1) of (\hyperlink{A4}{A.4}) is satisfied.

By (\ref{0-1c}), we see that
(2) of (\hyperlink{A4}{A.4}) is satisfied with $ n_{N}= N $.

By (\ref{0-1c}) and (\ref{0-1d}), we see that $ \muNg$ is a
$ (|z|^2 ,-2 \log|z|) $-canonical Gibbs measure.
Therefore, (3) of (\hyperlink{A4}{A.4}) holds with
$ \PhiN(z) = |z|^2 $ and $ \Psi(z) = -2 \log|z|$.

(4) of (\hyperlink{A4}{A.4}) is obvious with the above choice of
$ \PhiN$ and $ \Psi$, which completes the proof.
\end{pf}

We proceed with (\hyperlink{A6}{A.6}). For this, we prepare Lemma~\ref{l0-3}.
We denote $ \langle\mathsf{s},f \rangle= \sum_{i}f(s_i)$ for
$ \mathsf{s} = \sum_{i}\delta_{s_i} $.
We set $ \Dr= \{ z \in\mathbb{C} ; |z| < r \}$.
Let $ \arg z $ be the angle of $ z \in\mathbb{C}$; that is,
$ z = |z|e^{\ii\arg z }$.
We write $ f (r) = O (g(r)) \mbox{ as } r \to\infty$ if
$ \limsup_{r \to\infty}|f (r)|/|g(r)| < \infty$.

%
%
\begin{lem} \label{l0-3}
Let $ h_r (z) = 1_{\Dr}(z) e^{ \ii\ell\arg z } $,
where $ \ell\in\mathbb{Z} $. Let $ \map{f}{\mathbb{C}}{\mathbb
{C}} $ be
a bounded, measurable function such that
$ \sup_{|z| = r }|f(z)- z_0 | = O (r^{- 1 }) $ as
$ r \to\infty$ for some $ z_0 \in\mathbb{C}$.
We then have
%
\begin{equation}
\label{0-3b}  \sup_{N} \operatorname{Var}^{\mugN}\bigl(
\langle\mathsf{s}, h_r f \rangle\bigr) = O(r) \qquad\mbox{as $ r \to
\infty$.}
\end{equation}
\end{lem}
We remark that, if we replace $ \muNg$ by the Poisson random point
field whose intensity is the Lebesgue measure, then the right-hand side
of (\ref{0-3b}) becomes $ O(r^2)$.
Therefore Lemma~\ref{l0-3} implies the fluctuation of $ \{ \mugN\} $ is
uniformly small compared with that of the Poisson random point field.
Indeed, Lemma~\ref{l0-3} is the key to the proof of the quasi-Gibbs property.
Shirai \cite{shirai} initiated this kind of small fluctuation property
for the Ginibre random point field $ \mug$ with $ f=1 $.
In \cite{o-s} Shirai's result was generalized to functions $ f $ as above.
Lemma~\ref{l0-3} is its $ N$-particle version, which will be proved in
Section~\ref{sA2} below.

\begin{pf*}{Proof of Theorem~\ref{l23}}
Applying Theorems \ref{l41} and \ref{l63}, we deduce
Theorem~\ref{l23} from (\hyperlink{A4}{A.4}) and (\hyperlink{A6}{A.6}).
We note that (\hyperlink{A4}{A.4}) follows from Lemma~\ref{l0-2}.
Therefore, it only remains to prove two assumptions, (\ref{A7b}) and
(\ref{A7a}) of (\hyperlink{A6}{A.6}).
We check (\ref{A7b}) and (\ref{A7a}) for ${\ell_0}= 3 $.

By (\ref{0-1a}) and (\ref{0-1b}), we have $ \rgNx(z) \le\rgx(z)
= 1/\pi$. Therefore, we have
%
\begin{equation}
\label{0-4a} \int_{|z|\ge1} \frac{1}{|z|^3} \rgNx(z)\,dz \le
\frac{1}{\pi}\int_{|z|\ge1}\frac{1}{|z|^3} \,dz < \infty
.
\end{equation}
This implies (\ref{71z}). Hence, by Lemma~\ref{l71}, we obtain
(\ref{A7b}) with ${\ell_0}= 3 $.

We finally prove (\ref{A7a}).
Let $ u_{\ell,r}^{N} $ and $ \mathsf{u}_{\ell,r}^{N}$ be as in
Lemma~\ref{l81}. It is then easy to see that
\[
u_{\ell,r}^{N} (z)=\biggl(\frac{\lceil|z|\rceil}{|z|}
\biggr)^{\ell} 1_{\Sti
_{1r} }(z) e^{\ii\ell\arg z }.
\]
Hence, $ u_{\ell,r}^{N} $ satisfies the assumption of Lemma~\ref{l0-3}
with $ z_0 = 1 $
and $ f(z) = (\frac{\lceil|z|\rceil}{|z|} )^{\ell}$.
Therefore, by Lemma~\ref{l0-3}, we obtain
%
%
\begin{eqnarray}
\label{0-6}\qquad\quad \lim_{r\to\infty}r^{2c_{\scriptsize\ref{;81}}-2 \ell} \sup_{N}
\operatorname{Var}^{\mugN}\bigl[ \mathsf{u}_{\ell,r}^{N}
\bigr] = 0 \qquad\mbox{with $c_{\scriptsize\ref{;81}} = 1/4 $, say, for $ \ell= 1,2 $} .
\end{eqnarray}
Since $ E^{\mugN}[\mathsf{u}_{\ell,r}^{N}] = 0 $, (\ref{0-6}) implies
$ \lim_{r\to\infty}r^{2c_{\scriptsize\ref{;81}}-2 \ell}
\sup_{N\in\N} E^{\mugN}[|\mathsf{u}_{\ell,r}^{N}|^2]= 0 $, which
allows us to deduce (\ref{81b}).
Hence, by Lemma~\ref{l81}, we obtain (\ref{A7a}) with ${\ell_0}= 3 $.
We thus complete the proof.
\end{pf*}

\subsection{\texorpdfstring{Proof of Lemma~\protect\ref{l0-3}}{Proof of Lemma 10.2}} \label{sA2}
The purpose of this subsection is to prove Lemma~\ref{l0-3}.

Let $ \mathsf{g}(dz) = \frac{1}{\pi} \exp\{ - |z|^2 \} \,dz $ be
the standard complex Gaussian measure.
Let $ \{ \rN\}_{n \in\N}$ be
the correlation function of $ \mugN$ w.r.t. $ \mathsf{g}$.
Then $ \{ \rN\}_{n \in\N}$ is given by
%
\begin{equation}
\label{a11}  \rN(z_1,\ldots,z_n) = \det\bigl[
\KN(z_i,z_j)\bigr]_{i,j = 1,\ldots,n} ,
\end{equation}
where
$ \KN(z_1,z_2) = \sum_{k = 0}^{N-1} {\{ z_1\bar{z}_2 \}^k}/{k!} $.
We see that
$ \rN= \rgN\pi^{n}e^{|z_1|^2+\cdots+|z_n|^2}$ and
$ \KN(w,z) = \pi e^{|w|^2/2}\kgN(w,z)e^{|z|^2/2}$
by construction. Let
%
\begin{eqnarray}
\label{a12a} K(z_1,z_2) = \sum
_{k = 0}^{\infty} \frac{\{ z_1\bar{z}_2 \}^k}{k!} ,\qquad
\KNstar(z_1,z_2) = \sum_{k = N}^{\infty}
\frac{\{ z_1\bar{z}_2 \}^k}{k!} .
\end{eqnarray}
Then $ K = \KN+\KNstar$ by definition. Let
\begin{eqnarray*}
M^N_r = \int h_r ( w )
\overline{ h_r ( z )} \bigl\{ \bigl|K(w , z )\bigr|^2 - \bigl|\KN(w , z
)\bigr|^2 - \bigl|\KNstar(w , z )\bigr|^2 \bigr\} \mathsf{g}(dw )
\mathsf{g}(dz ) .
\end{eqnarray*}

%
\begin{lem} \label{la32}
Let $ e_N^s = \sum_{k = 0}^{N}s^k/k! $. Then
\[
\bigl|M^N_r\bigr| \le 2\bigl\{1-e^{-r^2}e_{N-1}^{r^2}
\bigr\} \bigl\{1-e^{-r^2}e_{N}^{r^2} \bigr\} .
\]
\end{lem}
\begin{pf}
From $ |K|^2 = |\KN|^2+|\KNstar|^2+\KN\KNstarover
+\KNstar\KNover$, we have
\begin{eqnarray*}
\bigl|M^N_r\bigr|& =& \biggl|\int h_r (w)\overline{
h_r }(z)\bigl\{ \KN\KNstarover +\KNstar\KNover\bigr\}\mathsf{g}(dw )
\mathsf{g}(dz )\biggr|
\\
& =& \frac{2}{(N-1)!N!}\biggl\{\int_{\Dr}|w|^{2N-1}
\mathsf{g}(dw )\biggr\}^2
\\
&\le& 2\biggl\{\frac{1}{(N-1)!}\int_{\Dr}|w|^{2N-2}
\mathsf{g}(dw )\biggr\} \biggl\{\frac{1}{N!}\int_{\Dr}|w|^{2N}
\mathsf{g}(dw )\biggr\}
\\
& =& 2\bigl\{1-e^{-r^2}e_{N-1}^{r^2} \bigr
\} \bigl\{1-e^{-r^2}e_{N}^{r^2} \bigr\} .
\end{eqnarray*}
This completes the proof.
\end{pf}

The kernel $ \KNstar$ also generates the
determinantal random point field denoted by~$ \mugNstar$.
%
\begin{lem} \label{la21} (1)
Let $ f $ be a bounded measurable function with compact support. Then
%
\begin{equation}
\label{a21c} \V\bigl(\langle\mathsf{s} , f \rangle\bigr) \le
\frac{2}{\pi} \int_{\mathbb{C}}\bigl | f (z)\bigr|^2 \,dz .
\end{equation}
(2) (\ref{a21c}) also hold for $ \mugNstar$ and $ \mug$.
\end{lem}

\begin{pf}
Since $ \KN(w,z)$ consists of a sum of pairs of
orthonormal functions w.r.t. $ \mathsf{g}(dz) $,
we have the equality
%
\begin{equation}
\label{a21d} \KN(z,z) = \int_{\mathbb{C}} \bigl|\KN(z,w)\bigr|^2
\mathsf{g}(dw) .
\end{equation}
By the standard calculation of correlation functions, we have
\begin{eqnarray*}
&& \V\bigl(\langle\mathsf{s} , f \rangle\bigr)
\\
&&\qquad =\int_{\mathbb{C}} \bigl| f (z)\bigr|^2 \KN(z,z)
\mathsf{g}(dz ) - \int_{\mathbb{C}^2} f (w)\overline{f (z)} \bigl|
\KN(w,z)\bigr|^2 \mathsf{g}(dw) \mathsf{g}(dz) .
\end{eqnarray*}
Combining these two equalities, and then using the inequalities
$ |a-b|^2\le2(|a|^2+|b|^2)$ and
$ |\KN(w,z)|^2\le\KN(w,w) \KN(z,z)$, we obtain
%
\begin{eqnarray}
\label{a21b}
\qquad \V\bigl(\langle\mathsf{s} ,  f \rangle\bigr)& =& \frac{1}{2}
\int_{\mathbb{C}^2} \bigl| f (w)- f (z)\bigr|^2 \bigl|
\KN(w,z)\bigr|^2 \mathsf{g}(dw) \mathsf{g}(dz)
\nonumber
\\[-8pt]
\\[-8pt]
\nonumber
& \le &\int_{\mathbb{C}^2} \bigl\{\bigl| f (w)\bigr|^2+\bigl| f
(z)\bigr|^2\bigr\} \KN(w,w) \KN(z,z) \mathsf{g}(dw) \mathsf{g}(dz) .
\end{eqnarray}
This, combined with the estimates $ 0 \le\KN(z,z)(1/\pi)
e^{-|z|^2}\le1/\pi$,
allows us to conclude (\ref{a21c}).
The proof of (2) is the same as that of (1).
\end{pf}
%
%
\begin{lem}[(Theorem 1.3 in \cite{o-s})] \label{l3}
$ \sup_{1\le r} {r}^{-1} \operatorname{Var}^{\mug} [ \langle\mathsf{s},
h_r f \rangle]
< \infty$.
\end{lem}
\begin{pf}
This lemma is a special case of Theorem 1.3 in \cite{o-s}.
\end{pf}
%
%
\begin{lem} \label{la25}
$ \sup_{1\le N }\sup_{1\le r} \frac{1}{r}\operatorname{Var}^{\mugN}
(\langle\mathsf{s},
h_r \rangle) < \infty$.
\end{lem}
\begin{pf}
By $ K = \KN+\KNstar$ and Lemma~\ref{la21}, we have
\[
 \V\bigl(\langle\mathsf{s}, h_r \rangle\bigr) =\Vz\bigl(\langle
\mathsf{s}, h_r \rangle\bigr)- M^N_r -\VV
\bigl(\langle\mathsf {s}, h_r \rangle\bigr) .
\]
By Lemma~\ref{la32}, we have $|M^N_r|\le2\{1-e^{-r^2}e_{N-1}^{r^2}\}
\{ 1-e^{-r^2}e_{N}^{r^2} \}$.
These, together with Lemma~\ref{l3}, complete the proof.
\end{pf}
\begin{pf*}{Proof of Lemma~\ref{l0-3}}
By $ h_r f = z_0 h_r + h_r ( f - z_0 ) $ and (\ref{a21c}), we have
\begin{eqnarray*}
\V\bigl(\langle\mathsf{s}, h_r f \rangle\bigr) & \le&2
z_0^2 \V\bigl(\langle\mathsf{s}, h_r
\rangle\bigr) + 2 \V\bigl(\bigl\langle\mathsf{s}, h_r ( f -
z_0 )\bigr\rangle\bigr)
\\
& \le&2 z_0^2 \V\bigl(\langle\mathsf{s} ,
h_r \rangle\bigr) + \frac{4}{\pi} \int_{ \Dr}
\bigl| h_r ( f - z_0 ) \bigr|^2 \,dz .
\end{eqnarray*}
Hence, from Lemma~\ref{la25} and the assumption
$ \sup_{|z| = r }|f(z)- z_0 | = O (r^{- 1 }) $, we complete the proof.
\end{pf*}

\begin{appendix}

\section*{Appendix}\label{sA}
\subsection{\texorpdfstring{Proof of Lemmas~\protect\ref{l34} and~\protect\ref{l35}}{Proof of Lemmas 3.4 and 3.5}}\label{sA1}
\mbox{}
\begin{pf*}{Proof of Lemma~\ref{l34}}
Let $ \{ f_{\mathsf{p}} \} $ be a $ \Ermk$-Cauchy sequence
in $ \domai$ such that
$ \lim\| f_{\mathsf{p}} \|_{\Lmk} = 0 $.
Then from (\ref{32a}) and (\ref{32b}), we see that $ \{ f_{\mathsf
{p}} \}$ satisfies
%
%
\setcounter{equation}{0}
\begin{eqnarray}
\label{33} \limi{\mathsf{p},\mathsf{q}} \int_{\SSS}
\Ermky(f_{\mathsf{p}}-f_{\mathsf{q}},f_{\mathsf{p}}-f_{\mathsf{q}})
\murkm(d\mathsf{s}) &= &0,
\\
\label{34}  \limi{\mathsf{p}} \int_{\SSS} \|
f_{\mathsf{p}} \|^2_{L^2(\SSS
_r^{m},\murmy)} \murkm(d\mathsf{s})& =& 0 .
\end{eqnarray}

We prove that $ \limi{\mathsf{p}} \Ermk(f_{\mathsf{p}},f_{\mathsf
{p}}) = 0 $.
For this purpose, it is enough to show that,
for any subsequence $ \{ f_{1,\mathsf{p}} \}$ of
$ \{ f_{\mathsf{p}} \} $, we can choose a subsequence
$ \{ f_{2,\mathsf{p}} \}$ of $ \{ f_{1,\mathsf{p}} \}$ such that
%
\begin{equation}
\label{35} \limi{\mathsf{p}} \Ermk(f_{2,\mathsf{p}} ,f_{2,\mathsf{p}}) = 0 .
\end{equation}
Therefore, let $ \{ f_{1,\mathsf{p}} \} $ be any subsequence of $ \{
f_{\mathsf{p}} \} $.
Then by (\ref{33}) and (\ref{34}),
we can choose a subsequence $ \{ f_{2,\mathsf{p}} \}$
such that $ \murkm(\mathsf{A}_{\mathsf{p}}) \le2^{-k} $ and
$ \murkm(\mathsf{B}_{\mathsf{p}}) \le2^{-k} $, where
\begin{eqnarray*}
\mathsf{A}_{\mathsf{p}} &=& \bigl\{ \mathsf{s} ; \Ermky(f_{2,\mathsf{p}}-f_{2,\mathsf{p}+1},f_{2,\mathsf
{p}}-f_{2,\mathsf{p}+1})
\ge2^{-2k} \bigr\},
\\
\mathsf{B}_{\mathsf{p}}& = &\bigl\{ \mathsf{s} ; \| f_{\mathsf{p}}
\|^2_{L^2(\SSS_r^{m},\murmy)} \ge2^{-2k} \bigr\} .
\end{eqnarray*}
Hence, from the Borel--Cantelli lemma, we see that
\[
\murkm(\limsup\mathsf{A}_{\mathsf{p}}) = \murkm(\limsup\mathsf{B}_{\mathsf{p}})
= 0 .
\]
This means that,
for $ \murkm$-a.s. $ \mathsf{s} $, the sequence
$ \{ f_{2,\mathsf{p}} \} $ is an $ \Ermky$-Cauchy sequence
converging to $ 0 $ in $ L^2(\SSS_r^{m},\murmy) $ as
$ \mathsf{p}\to\infty$. Therefore, by assumption, we have
%
\begin{equation}
\label{36} \limi{\mathsf{p}} \Ermky(f_{2,\mathsf{p}},f_{2,\mathsf{p}}) = 0
\qquad\mbox{for $ \murkm$-a.s. $ \mathsf{s} $. }
\end{equation}

Let $ \murmch$ be the symmetric measure on $ \Sr^m $
such that $ \murmch\circ\iota^{-1} = \murmy$.
For $ f_{2,\mathsf{p}} $, there exists a function
$ \map{f_{2,\mathsf{p}}^{r,m}}{\Sr^m \ts\SSS}{\R} $ such that
$ f_{2,\mathsf{p}}^{r,m} (\mathbf{x}, \mathsf{s})$
is symmetric in
$ \mathbf{x} = (x_1,\ldots,x_m) $
for each $ \mathsf{s} \in\SSS$ and that
$ f_{2,\mathsf{p}}^{r,m}(\mathbf{x}, \mathsf{s}) =
f_{2,\mathsf{p}}(\mathsf{s})$
for $ \mathsf{s} \in\SSS_r^{m}$ decomposed as
$ \mathsf{s} = \iota(\mathbf{x}) +
\pi_{\Sr^c} (\mathsf{s}) $.
Let $ x_l = (x_{l1},\ldots,x_{ld}) \in\Rd$. Then
\begin{eqnarray*}
&& \int_{\SSS_r^{m}} \Ermky (f_{2,\mathsf{p}}-f_{2,\mathsf{p}+1},f_{2,\mathsf{p}}-f_{2,\mathsf{p}+1})
\murkm(d\mathsf{s})
\\
&&\qquad =\int_{\Sr^m \ts\SSS} \frac{1}{2} \sum
_{l = 1}^m \sum_{i,j = 1}^{ d }
a_{ij} (\mathsf{s} , x_l ) \PD{(f_{2,\mathsf{p}}^{r,m}-f_{2,\mathsf{p}+1}^{r,m}
)} {x_{li}}\\
&&\hspace*{107pt}{}\times \PD{(f_{2,\mathsf{p}}^{r,m}-f_{2,\mathsf{p}+1}^{r,m}
)} {x_{lj}} \murmch(d\mathbf{x}) \murkm(d\mathsf{s}) .
\end{eqnarray*}
Hence, by (\ref{33}), we see that the vector-valued function
$ \map{(\nabla_{x_l}f_{2,\mathsf{p}}^{r,m})_{l = 1,\ldots,m}}\break
{\Sr^m \ts\SSS}{(\Rd)^m} $ is a Cauchy sequence in
$ L^2(\Sr^m \ts\SSS\to(\Rd)^m , \murmch) $,
where we equip
$ L^2(\Sr^m \ts\SSS\to(\Rd)^m , \murmch) $
with the inner product
\[
(\mathbf{f},\mathbf{g}) = \int_{\Sr^m \ts\SSS} \sum
_{l = 1}^m \bigl\{ f_l(\mathbf{x} ,
\mathsf{s}) g_l(\mathbf{x} , \mathsf{s}) %
a_0 \bigl(\iota(\mathbf{x}) + \pi_{\Sr^c}(\mathsf{s}) ,
x_l \bigr) \bigr\} \murmch(d\mathbf{x}) \murkm(d\mathsf{s}) .
\]
Here $ \mathbf{f} = (f_1,\ldots,f_m) $, and
$ a_0 $ is the function in (\ref{A33}).
Combining this with (\ref{33}) and (\ref{36}),
we obtain (\ref{35}), which completes the proof.
\end{pf*}

\begin{pf*}{Proof of Lemma~\ref{l35}}
By (\ref{A33}), we deduce the closability of
$ (\Ermky, \break \domai)$ on $ L^2(\SSS_r^{m}, \murmy)$ from that of
$ (\Eaz, \domaz)$ on $ L^2(\SSS_r^{m}, \murmy) $.
Here $ I $ is the $ d\ts d$ unit matrix.

Let $ \murmch$ be as in the proof of Lemma~\ref{l34}.
Then by (\ref{qg2}), $ \murmch$ has a density $ \check{\sigma
}(\mathbf{x})$
w.r.t. $ e^{-\mathcal{H}_{r}(\mathbf{x})}\,d\mathbf{x}$.
Here $ d\mathbf{x}$ denotes the Lebesgue measure on $ \Sr^m $, and
we regard $ e^{-\mathcal{H}_{r}}$ as a symmetric function on $ \Sr^m$
in an obvious manner. In the following we use the same convention for
functions on
the configuration space $ \SSS_r^{m}$.
We note that according to (\ref{qg2}), $ \check{\sigma}$ is
uniformly positive and bounded on
$ \Sr^m $.

Let
$ O_p = \{ \mathbf{x}\in\Sr^m ; p^{-1}< a_0(\mathbf{x})\}
\cap\{ \mathbf{x}\in\Sr^m ; p^{-1}< e^{-\mathcal
{H}_{r}(\mathbf{x})}\}$
$( p\in\N)$.
Recall that $ a_0$ and $ e^{-\mathcal{H}_{r}}$ are lower semicontinuous
(the latter claim follows from the assumption that $ \Phi$ and $ \Psi
$ are upper semicontinuous), which implies that $ O_p $ is an open set.
Moreover, $ \{ O_p\} $ is nondecreasing in $ p$.
Let $ \varepsilon_p $ be the bilinear form on $ \Sr^m $ defined by
\[
\varepsilon_p(f,g) = \int_{O_p}
\mathbb{D}^{a_0I}[f,g] \check{\sigma}e^{-\mathcal{H}_{r}} \,d\mathbf{x} =
\int_{O_p} \mathbb{D}[f,g] a_0\check{
\sigma}e^{-\mathcal{H}_{r}} \,d\mathbf{x} .
\]
Recall that $ a_0\check{\sigma}e^{-\mathcal{H}_{r}}$ and $ \check
{\sigma}e^{-\mathcal{H}_{r}}$ are bounded on $ \Sr^m $ and greater
than or equal to $ p^{-2}$ on $ O_p$.
Hence, $ (\varepsilon_p ,C^{\infty}_b(\Sr^m))$ is closable on
$ L^2(\Sr^m, \check{\sigma}e^{-\mathcal{H}_{r}}\,d\mathbf{x}) =
L^2(\Sr^m, \murmch)$.
Since $ \{ O_p \} $ is nondecreasing, the sequence of closable
bilinear forms $ (\varepsilon_p ,C^{\infty}_b(\Sr^m))$ is nondecreasing.
Hence, the limit bilinear form
$ (\varepsilon_{\infty}, C^{\infty}_b(\Sr^m))$
is also closable on $L^2(\Sr^m, \murmch)$; see \cite{mr},
Proposition 3.7, page 30. We used here
$ \{ f\in C^{\infty}_b(\Sr^m);\varepsilon_{\infty}(f,f)<\infty\}
=C^{\infty}_b(\Sr^m)$.

It is easy to see that the closability of
$ (\varepsilon_{\infty},C^{\infty}_b(\Sr^m)) $ on
$L^2(\Sr^m, \murmch)$ implies the closability of
$ (\Eaz, \domaz)$ on $ L^2(\SSS_r^{m}, \murmy) $.
Hence, we complete the proof.
\end{pf*}

\subsection{\texorpdfstring{The weak convergence of $\{\muN\}$}{The weak convergence of \{mu N\}}}
\label{sA4}
In Section~\ref{s4}, we considered the fact that the measures $ \{
\muN\}
$ in (\hyperlink{A2}{A.2}) converge weakly to $ \mu$. For the sake of
completeness we give a proof of this.
Let $ \Srhat= \{ x \in\SSSS ; |x|< r \} $ and $ \Dr^n = \prod_{m=1}^{n} \{ |x_m|< r \} $
as before.

%
\begin{lemm} \label{l111}
Assume (\ref{41a}) and (\ref{41b}) in \textup{(\hyperlink{A4}{A.4})}.
Then $ \limi{N} \muN= \mu$ weakly.
\end{lemm}
\begin{pf}
A permutation invariant function
$ \map{m_r^n}{\tilde{\SSSS}_r^n }{\R} $ is by definition
the $ n $-density function of $ \mu$ if, for any bounded
$ \sigma[\pi_{\Srhat}] $-measurable function~$ \mathsf{f}$,
\[
\int_{\tilde{\SSS}_r^n } \mathsf{f} \,d\mu = \frac{1}{n!} \int
_{\tilde{\SSSS}_r^n } f_r^n m_r^n
\,dx_1\cdots \,dx_n ,
\]
where
$ \tilde{\SSS}_r^n = \{ \mathsf{x} \in\SSS ; \mathsf{x}
(\Srhat) = n \}$, and
$ \map{f_r^n}{\tilde{\SSSS}_r^n }{\R} $ is the permutation
invariant function such that
$ f_r^n (x_1,\ldots,x_n) = \mathsf{f} (\mathsf{x})$ for
$ \mathsf{x} \in\tilde{\SSS}_r^n $ such that
$ \pi_{\Srhat}(\mathsf{x}) = \sum_i \delta_{x_i} $.

Let $ m^{n}_{N,r} (x_1,\ldots,x_n)$ [resp., $ m^n_r (x_1,\ldots
,x_n)$] be
the $ n $-density function of $ \muN$ (resp., $ \mu$) on $ \Srhat$.
Then by (\ref{41b}), we easily see that
%
\begin{eqnarray}
\label{70c}&& m^{n}_{N,r} (x_1,
\ldots,x_n)
\nonumber
\\[-8pt]
\\[-8pt]
\nonumber
&&\qquad =\sum_{k = 0}^{\infty}
\frac{(-1)^k}{k!} %
\int_{\tilde{\SSSS}_r^n }
\rho_N^{n+k} (x_1,\ldots,x_{n+k})
\,dx_{n+1}\cdots \,dx_{n+k} .
\end{eqnarray}
Combining (\ref{41a}) and (\ref{41b}) with
(\ref{70c}) and the same equality as (\ref{70c})
for $ \mu$ and applying the bounded convergence
theorem, we obtain for each $ r,n \in\N$,
\begin{eqnarray*}
\sup_{N}\sup_{ \tilde{\SSSS}_r^n }\bigl | m^{n}_{N,r}
(x_1,\ldots,x_n)\bigr|& <& \infty,
\\
\lim_{ N \to\infty} m^{n}_{N,r}
(x_1,\ldots,x_n)& =& m^n_r
(x_1,\ldots,x_n) \qquad\mbox{a.e.}
\end{eqnarray*}
From this, we see that the measures satisfy
$\limi{N}\muN\circ\pi_{\Srhat}^{-1} = \mu\circ\pi_{\Srhat}^{-1}$
weakly in $ \pi_{\Srhat}(\SSS) $ for all $ r $.
Hence, it only remains to prove that
the sequence $ \{ \muN\} $ is tight in $ \SSS$.

Now we recall a closed subset $ \SSS_0 $ in $ \SSS$
is compact if and only if there exists an increasing sequence
$ \mathbf{a} = \{ a_r \}_{r\in\N}$
of natural numbers such that
$ \sup_{\mathsf{s}\in\SSS_0 }\mathsf{s} (\Srhat) \le a_r $
for all $ r \in\N$ \cite{resnick}, Section~3.4.
Let $\mathsf{K}(r,a) =
\{\mathsf{s}; \mathsf{s} (\Srhat)\le a \}$.
Set
$\mathsf{K}(\mathbf{a}) = \bigcap_{r \in\N} \mathsf{K}(r,a_r) $
for $\mathbf{a} = \{ a_r \}_{r\in\N}$.
We then see that the set $ \mathsf{K}(\mathbf{a}) $ is compact in
$ \SSS$ because of the equivalence condition given above.

Let $ \epsilon>0 $ be fixed.
Note that $\pi_{\Srhat}(\SSS)$ is also a Polish space because
$ \Srhat$ is Polish~\cite{resnick}, Proposition~3.17.
Since $\{\muN\circ\pi_{\Srhat}^{-1}\}$ is tight
as probability measures in $\pi_{\Srhat}(\SSS)$,
there exists a compact set
$ \mathsf{K}_{r} $ in $ \pi_{\Srhat}(\SSS) $ such that
%
\begin{equation}
\label{70f} \sup_N \muN\circ\pi_{\Srhat}^{-1}
\bigl(\mathsf{K}_{r}^c\bigr)\le\epsilon2^{-r}
.
\end{equation}
Moreover there exists an $ a_r \in\N$ such that
$ \mathsf{K}_{r}\subset\mathsf{K}(r,a_r)$
because $ \mathsf{K}_{r} $ is compact.
We can and do take $ a_r \in\N$ in such a way that $ a_r < a_{r+1}$.
By (\ref{70f}) and $ \mathsf{K}_{r}\subset\mathsf{K}(r,a_r)$,
we have $ \sup_N \muN(\mathsf{K}(r,a_r)^c)\le\epsilon2^{-r}$.
Hence, for $\mathbf{a} = \{ a_r \}_{r\in\N}$, we have
\begin{eqnarray*}
\sup_N \muN\bigl( \mathsf{K}(\mathbf{a})^c \bigr) =
\sup_N \muN\biggl(\bigcup_{r \in\N}
\mathsf{K}(r,a_r)^c \biggr) \le\sup_N \sum
_{r\in\N} \muN\bigl(\mathsf{K}(r,a_r)^c
\bigr) \le\epsilon .
\end{eqnarray*}
This implies $ \{ \muN\} $
is tight, which completes the proof.
\end{pf}

\subsection{\texorpdfstring{Proof of (\protect\ref{71a3}) and (\protect\ref{71a4})}{Proof of (9.30) and (9.31)}} \label{sA3}
In this subsection, we prove (\ref{71a3})\break and~(\ref{71a4}).
Let $ J_N (x) = I_N (x) -\frac{1}{2} \operatorname{sgn}(x)$.
Note that $ S_N $ is an even function, and
$ I_N $, $ D_N $ and $ J_N $ are odd functions.
By (\ref{"91s}) for $ S_N $, $ D_N $ and $ I_N $,
%
%
\begin{eqnarray}
\label{81}
&&
 \KsinNx(x)\KsinNx(-x)
\nonumber\\
&&\qquad = \Theta\biggl( \lleft[\matrix{ S_N (x)&
D_N (x) \vspace*{2pt}
\cr
J_N (x) & S_N (x) }
\rright] \lleft[\matrix{ S_N (-x)& D_N (-x)
\vspace*{2pt}
\cr
J_N (-x) & S_N (-x) } \rright]\biggr)
\\
&&\qquad = \Theta\biggl( \lleft[\matrix{ S_N (x)^2
-D_N (x) J_N (x) & 0 \vspace*{2pt}
\cr
0 & S_N
(x)^2 - D_N (x) J_N (x) } \rright]\biggr) .\nonumber
\end{eqnarray}
Hence, by (\ref{"91q}) and (\ref{71d}), we have
$\mathcal{T}^{N}_{1} = S_N^2 - D_N J_N. $
This, combined with~(\ref{71y})--(\ref{71v1}), yields (\ref{71a3}).
We consider (\ref{71a4}) next. By (\ref{"91t}) for $ S_N $, $ D_N
$ and $ I_N $,
we see that
\begin{eqnarray*}
&& \KsinNz(x)\KsinNz(-x)\nonumber
\\
& &\qquad=  
\Theta\biggl( \lleft[\matrix{ S_N
(2x)& D_N (2x) \vspace*{2pt}
\cr
I_N (2x) &
S_N (2x) } \rright] \lleft[\matrix{ S_N (-2x)&
D_N (-2x) \vspace*{2pt}
\cr
I_N (-2x) & S_N
(-2x) } \rright]\biggr)
\\
\nonumber
&&\qquad=  
\Theta\biggl( \lleft[\matrix{ S_N
(2x)^2 -D_N (2x) I_N (2x) & 0 \vspace*{2pt}
\cr
0 & S_N (2x)^2 - D_N (2x) I_N
(2x) } \rright] \biggr).
\end{eqnarray*}
Hence, by (\ref{"91q}) and (\ref{71d}), we have
$\mathcal{T}^{N}_{4}( x ) = S_N (2x)^2 - D_N (2x) J_N (2x)$.
This, combined with (\ref{71y})--(\ref{71v1}), yields (\ref{71a4}).
\end{appendix}

%

%


\printaddresses

\end{document}